\magnification=1100
\input amstex
\documentstyle{amsppt}
\overfullrule=0pt
\vsize=500pt
\voffset=1.5cm

\document

\def\G{{\Gamma}}

\def\hen{\frak{h}}

\def\op{\operatorname}
\def\ind{{\lim\limits_{\longrightarrow}}}
\def\pro{{\lim\limits_{\longleftarrow}}}
\def\lra{\,\longrightarrow\,}
\def\bda{\big\downarrow}

\def\a{{\alpha}}
\def\o{{\omega}}
\def\l{{\lambda}}
\def\b{{\beta}}
\def\g{{\gamma}}
\def\d{{\delta}}
\def\eps{{\varepsilon}}

\def\1b{{\bold 1}}

\def\Cb{{\bold C}}

\def\Ob{{\bold O}}

\def\Sum{{\ts\sum}}

\def\Prod{{\ts\prod}}

\def\H{{\roman H}}

\def\H{{\roman H}}

\def\top{{\text{top}}}

\def\Spec{\roman{Spec\,}}
\def\Spf{\roman{Spf\,}}
\def\Sch{{\text{Sch}}}

\def\Db{{\bold{D}}}
\def\red{{\text{red}}}

\def\ft{{\text{ft}}}
\def\sqr{{\!\sqrt{\,\,\,}}}

\def\Algb{{\bold{Alg}}}
\def\Affb{{\bold{Aff}}}
\def\Schb{{\bold{Sch}}}
\def\Ischb{{\bold{Isch}}}
\def\Lcb{{\bold{Lc}}}
\def\Shfb{{\bold{Shf}}}
\def\Prob{{\bold{Pro}}}
\def\Funb{{\bold{Fun}}}
\def\Lrsb{{\bold{Lrs}}}
\def\Setsb{{\bold{Sets}}}

\def\Hom{\roman{Hom}}

\def\Id{\roman{Id}}

\def\Res{\text{Res}}

\def\End{\text{End}\,}

\def\Hilb{\text{Hilb}\,}

\def\Ker{\text{Ker}\,}

\def\max{\text{max}\,}

\def\Hilb{\text{Hilb}}

\def\AA{{\Bbb A}}

\def\EE{{\Bbb E}}

\def\NN{{\Bbb N}}

\def\QQ{{\Bbb Q}}

\def\ZZ{{\Bbb Z}}

\def\Vac{{\Vc ac}}
\def\Ac{{\Cal A}}

\def\Cc{{\Cal C}}
\def\Dc{{\Cal D}}
\def\Ec{{\Cal E}}
\def\Fc{{\Cal F}}

\def\Hc{{\Cal H}}

\def\Kc{{\Cal K}}
\def\Lc{{\Cal L}}
\def\Mc{{\Cal M}}
\def\Nc{{\Cal N}}
\def\Oc{{\Cal O}}

\def\Rc{{\Cal R}}
\def\Sc{{\Cal S}}
\def\Tc{{\Cal T}}
\def\Uc{{\Cal U}}
\def\Vc{{\Cal V}}

\def\Xc{{\Cal X}}

\def\and{{\quad\text{and}\quad}}

\def\ts{\textstyle}
\def\ss{\scriptstyle}

\def\qed{\hfill $\sqcap \hskip-6.5pt \sqcup$}        

\def\simto{\,{\buildrel\sim\over\to}\,}
\def\lra{{{\longrightarrow}}}

\def\rra{{\rightrightarrows}}

\def\u1{{\underline 1}}

\newdimen\Squaresize\Squaresize=14pt
\newdimen\Thickness\Thickness=0.5pt
\def\Square#1{\hbox{\vrule width\Thickness
              \alphaox to \Squaresize{\hrule height \Thickness\vss
              \hbox to \Squaresize{\hss#1\hss}
              \vss\hrule height\Thickness}
              \unskip\vrule width \Thickness}
              \kern-\Thickness}
\def\Vsquare#1{\alphaox{\Square{$#1$}}\kern-\Thickness}

%
\title Vertex algebras and the formal loop space\endtitle
\rightheadtext{Vertex algebras and formal loops}
\author M. Kapranov and E. Vasserot\endauthor

\endtopmatter
\document

\head 1. Introduction\endhead

\vskip .5cm

One of the salient mathematical features of string theory
is the importance of vertex algebras. Their role in the
theory can be compared to that of Lie algebras in
the ``ordinary'' physics of point particles.

\vskip .3cm

Mathematically, the approach of string theory
can be cast in terms of analysis on the space of free loops,
i.e., smooth maps $S^1\to X$ where $X$ is a given ``spacetime''
manifold.
Accordingly, one has the folklore principle
that  constructions involving the space of free loops
lead to vertex algebras. One class of  such constructions is provided by
the spaces of highest weight representations of loop groups.
Another is $\Omega_X^{ch}$,  the  chiral de Rham complex
of an algebraic variety $X$, introduced
by Malikov, Schechtman and Vaintrob [MSV].
Heuristically, this complex should be interpreted
in terms of $LX$, the space of free loops
and its subvariety $L^0 X$ consisting of loops extending holomorphically into
the unit disk.  More precisely, $\Omega^{ch}_X$
can be thought of as the semiinfinite de Rham complex
with coefficients in the space of distributions on $LX$  supported
on $L^0X$. This is not, however, the way $\Omega_X^{ch}$ has been
defined mathematically.  The  definition given in [MSV]
is of more computational nature and
proceeds by constructing the action of the group
of diffeomorphisms on  the irreducible module over the Heisenberg algebra.
In that approach it seems miraculous that such an action
exists at all.

\vskip .3cm

The aim of this paper is twofold. First, to give a precise
mathematical theorem underlying the above folklore principle
about vertex algebras. For this, we introduce an algebro-geometric
version of the free loop space $\Lc (X)$ for any scheme $X$ of finite\
type over a field. This is an ind-scheme containing $\Lc^0(X)$,
the scheme of formal germs of curves on $X$ studied  in [DL].
We prove that both $\Lc(X)$ and $\Lc^0(X)$ themselves
possess a non-linear version of the vertex algebra structure
(which makes it clear that any natural linear construction applied to them
should give a vertex algebra in the usual sense).
More precisely, we use the geometric approach to vertex algebras
developed by Beilinson and Drinfeld [BD1] and based on the concepts of chiral
algebras and factorization algebras.
The latter concept has a natural nonlinear version, that of a factorization
monoid.
 What we prove
is that  natural ``global''versions of $\Lc(X)$, $\Lc^0(X)$
have natural structures of factorization monoids.
An earlier known example of a factorization monoid is given by the affine
Grassmannian [G], and this explains why the spaces of
representations of loop groups are vertex algebras. Our construction
is   similar in spirit.

\vskip .3cm

To give a good definition of the algebraic analog of the full loop
space $LX$ one has to overcome a certain subtlety.
Namely, a natural approach would be to try to (ind-)represent
a functor which to any commutative ring $R$ associates
the set of $R((t))$-points of $X$. (This is exactly how one
defines the scheme $\Lc^0(X)$, with $R[[t]]$ instead of
$R((t))$.) If $X$ is affine, this indeed gives a good ind-scheme
which we denote $\tilde\Lc(X)$. But when $X$ is, say, projective,
then (for $R$ a field) there is no difference between
$R[[t]]$-points and $R((t))$-points of $X$ (valuative
criterion
of properness), so it may seem that nothing is gained
by allowing Laurent series. To state this phenomenon
differently, the ind-schemes $\tilde\Lc(U)$ for affine $U\subset X$
do not glue together well. This is in fact understandable
on general grounds: the loop space $LX$ is not the union
of the $LU$ since a loop need not spend all its time in
any given $U$.

To get around this difficulty we adopt the following strategy.
For an affine $X$ we consider $\Lc(X)$, the formal
neighborhood of $\Lc^0(X)$ in $\tilde\Lc(X)$.
So we are dealing with formal loops which
 are ``infinitesimal
in the Laurent direction''.Then,
we prove that the $\Lc(U), U\subset X$, do indeed
possess the right gluing properties. This is
due to the infinitesimal nature of our loops.

The role of nilpotent thickenings in
 Laurent series models for loop spaces
 was first pointed out by C. Contou-Carr\`ere  [CC] who was studying,
in our notation, the group ind-scheme $\tilde\Lc (\Bbb G_m)$
and found that it is a nontrivial formal thickening of $\Lc^0(\Bbb G_m)\times
\Bbb Z$.

\vskip .3cm

Our second goal is to give a direct geometric
 construction
of $\Omega^{ch}_X$ (for smooth $X$) in terms of
our model for the loop space.
By the above, this  construction explains also the
fact that
$\Omega^{ch}_X$ is a sheaf of vertex algebras.
In order to achieve this, we represent $\Lc(X)$ as an ind-pro-object
in the category of schemes of finite type and then
show that the shifted de Rham complexes of the terms
 of this ind-pro-system
arrange naturally into a double inductive system whose
inductive limit is identified with $\Omega^{ch}_X$.

\vskip .3cm

As with the study of formal arcs and motivic integration
[DL], one can view our considerations as algebro-geometric
analogs of the basic constructions of $p$-adic analysis.
The difference between our ind-scheme $\Lc(X)$ and
the more familiar scheme $\Lc^0(X)$ is similar to
the difference between $\QQ_p$ and $\ZZ_p$: while
the latter is a pro-object in the category of
schemes of finite type (resp. finite sets),
the former is an
ind-pro-object.
 Further, our approach to
$\Omega^{ch}_X$ is similar to the construction
of the space of locally constant functions with compact support  on
$\QQ_p = \ind_i \pro_j
p^{-i}\ZZ_p/p^j\ZZ_p$ as the double inductive limit of
the spaces of functions on the finite sets $p^{-i}\ZZ_p/p^j\ZZ_p$,
cf [P]. Notice that the reason  that these spaces of
functions indeed form a double inductive system
(with respect to the maps of inverse image in the $j$-direction
and direct image in the $i$-direction)
is  that the commutative squares in the ind-pro-system
$p^{-i}\ZZ_p/p^j\ZZ_p$ are  Cartesian (so that
 we have base change).
This is an algebraic
counterpart of the property of the local compactness of $\QQ_p$,
see [Kat].
In our situation it is equally important that
 the ind-scheme
$\Lc(X)$ satisfies a certain formal analog of local
compactness.

\vskip .3cm

This work has been done over a period of several years
during which the first author benefitted from
visits to  and the  financial support of the
Universit\'e Cergy-Pontoise,  the IHES, the Ecole
Normale Superieure and the Max-Planck-Institut
f\"ur Mathematik. His research was also partially supported
by grants from NSF and NSERC.
The second author is partially supported by EEC grant
no. ERB. FMRX-CT97-0100. We would like to thank A. Beilinson for
several useful remarks on the earlier version and V. Drinfeld for
pointing out the reference [TT].

\vskip 1cm

\head 2. Construction of the formal loop space\endhead

\vskip .3cm

\noindent {\bf (2.1) Generalities on schemes and ind-schemes.}
If $\Cb$ is a category, then we denote by
$\bold{Ind}(\Cb)$ and
$\Prob(\Cb)$ the categories of ind- and pro-objects of
$\Cb$,
see [AM] [GV] for background. Thus, objects of $\bold{Ind}(\Cb)$ (resp.
$\Prob(\Cb)$) are symbols $``\ind_i" C_i$ (resp. $``\pro_i" C_i$)
where $(C_i)$ is a filtering inductive (resp. projective)
system over $\Cb$, with morphisms defined as in {\it loc. cit.}
Recall that $\bold {Ind}(\Cb)$ can be considered as a full subcategory
in $\Funb^\circ (\Cb, \bold{Sets})$, the category of contravariant functors.

Throughout the paper we fix a field $k$.
We denote by $\Schb \subset \Lrsb$  the categories of
schemes and locally ringed spaces over $k$.
If $R$ is a commutative ring, we will write $\underline{\op{Spec}}\, R$
for the topological space
(the set of prime ideals with the Zariski topology)
 underlying the affine scheme $\op{Spec}\, R$
which is thus the ringed space $(\underline{\op{Spec}}\, R, \, \Cal O_{\op{
Spec}\, R})$.

By an ind-scheme we will mean in this paper an ind-object
of $\Schb$ represented by an inductive system of closed embeddings
of quasi-compact schemes.
The category of ind-schemes
will be denoted by $\Ischb$. Let us make the category $\Schb$
into a Grothendieck site by using Zariski open coverings and let $\Shfb$
be the category of sheaves of sets on $\Schb$.
For any ind-scheme $Y$ the functor $\eta_Y$ on schemes represented
by $Y$ is then a sheaf, so we have the embeddings
$$\Schb\subset\Ischb\subset\Shfb\subset\Funb^\circ(\Schb,
\Setsb).\leqno (2.1.1)$$

Since the category $\Schb$ has finite inverse limits, so
do all the categories in (2.1.1) and the embeddings preserve them.
On the contrary, finite direct limits, such as cokernels (when they exist)
are preserved by the first two of the embeddings but not by the third one:
cokernels in the category of sheaves are not the same as in the
category of all functors (presheaves).

We denote by $\Algb$ the category of $k$-algebras and by $\Affb
\subset \Schb$ the dual category of affine schemes. Note
that $\Ischb$ can be as well realized as a full subcategory
in $\Funb^\circ(\Affb,\Setsb)=\Funb(\Algb,\Setsb)$.

Given two contravariant functors $\phi, \phi':\Schb\to\Setsb$, and a morphism
$F: \phi'\to\phi$, we will say that $F$ is formally smooth
(resp. formally \'etale), if for any nilpotent
extension of affine schemes $S\subset S'$
the natural map
$$\phi'(S')\to\phi(S')\times_{\phi(S)}\phi'(S)$$
is surjective (resp. bijective).

We will say that $F$ is an open embedding, if for any scheme $S$ and any
$u\in \phi(S)$ (which is the same as a morphism $\eta_S\to \phi$)
the fiber product functor
$\eta_S\times_{\phi} \phi'$ is representable by a scheme $S'$ whose
natural morphism to $S$ is an open embedding.

We define formal smoothness and openness for morphisms of ind-schemes
by considering their representable functors.

For a scheme $Z$ we denote by $Z_\red\subset Z$ the corresponding
reduced subscheme.   We extend this notation to ind-schemes by applying
it term by term in inductive systems.

Let $X$ be a $k$-scheme of finite type. We denote by $\Affb_X$
(resp. $\Affb_X^\ft$) the category of schemes affine over $X$
(resp. affine of finite type over $X$). For future use let
us quote the following fact \cite{EGAIV, Corollary 8.13.2}.

\proclaim{(2.1.2) Proposition} The category $\Prob(\Affb_X^\ft)$
is equivalent to $\Affb_X$ via the functor $``\pro_n" S_n\mapsto\pro_n S_n$.
\endproclaim

We also denote $\Schb_X$ the category of all $X$-schemes and
$\Ischb_X$ the category of ind-schemes over $X$.
Thus objects of $\Ischb_X$ are arrows $Y\to X$, $Y\in\Ischb$, or,
equivalently, symbols 
$``\ind_n" Y_n$ where $Y_n\to X$ form
an inductive system of closed embeddings of quasicompact $X$-schemes.

\vskip .3cm

\noindent {\bf (2.2) The scheme of germs of arcs.}

Let $X$ be a scheme. We denote by $\Lc^0(X)$ the scheme of germs of arcs on $X$,
see [DL]  \cite{BLR, Theorem 7.6.4}. It represents the following
covariant functor  $\lambda^0_X$ on the category $\Algb$:
$$\lambda^0_X: R\mapsto \op{Hom}_{\Schb}(\op{Spec}\, R[[t]], X).$$
Here are some of the well-known properties of $\Lc^0(X)$.
Note that if $R$ is a local ring with maximal ideal $M$, then
$R[[t]]$ is a local ring with maximal ideal $M[[t]] + tR[[t]]$.

\proclaim {(2.2.1) Proposition}
(a) For any scheme $S$ the pair $(S, \Oc_S[[t]])$ is a locally ringed
space.

(b) For any scheme $S$ we have
$$\op{Hom}_{\Sch}(S, \Lc^0(X)) = \op{Hom}_{\Lrsb}( (S, \Cal O_S[[t]]), X).$$

(c) The scheme $\Lc^0(X)$ is the projective limit
of the schemes $\Lc^0_n(X)$, $n\in\NN$, representing the functors
$$\l^0_{n,X}: R\mapsto \op{Hom}_{\Schb} \bigl( \op{Spec}\,
R[t]/t^{n+1}, \, X\bigl).$$
 If $X$ is of finite type, than so is each $\Lc^0_n(X)$.

(d) Denote
$\pi_n: \Lc^0_n(X)\to X$, $\pi: \Lc^0(X)\to X$ the natural projections.
They are affine morphisms. For an open subset $U\subset X$
we have $\pi^{-1}(U)=\Lc^0(U)$ and $\pi_n^{-1}(U) = \Lc^0_n(U)$.

(e) If $X$ is smooth, then  so is $\Lc^0_n(X)$ and $\Lc^0(X)$ is
 formally smooth.
\endproclaim

\noindent {\sl Proof:} For (a) and (b) it suffices to assume that
$S= \op{Spec}\, R$ is affine.
 The embedding of
constant series and the evaluation at 0 give ring
homomorphisms
$R\buildrel\alpha\over\to R[[t]]\buildrel\beta\over\to R$
 and hence morphisms of topological spaces
$$\underline{\op{Spec}}\, R \buildrel p\over\longleftarrow
\underline{ \op{Spec}}\, R[[t]]
\buildrel i\over\longleftarrow \underline{\op{Spec}}\, R.\leqno (2.2.2)$$
The statement (a) follows from the next lemma, since
$\Cal O_{\op{Spec}\, R[[t]]}$ is obviously a sheaf of local rings.

\proclaim{(2.2.3) Lemma} We have
$$ i^{-1} \Cal O_{\op{Spec}\, R[[t]]} = \Cal O_{\op{Spec}\, R}[[t]] =
p_*\Cal O_{\op{Spec}\, R[[t]]}.$$
\endproclaim

\noindent {\sl Proof:}
 Let us prove the first equality.
If $\frak p\in\underline{\op{Spec}}\, R$ is a prime ideal in $R$, then
the stalk of $i^{-1} \Cal O_{\op{Spec}\, R[[t]]}$ at $\frak p$
is the localization  of $R[[t]]$ with respect to the multiplicative subset
$\beta^{-1}(R-\frak p)$ while the stalk of
$\Cal O_{\op{Spec}\, R}[[t]]$ is $R[[t]]\otimes_R R_{\frak p}$
where $R_{\frak p}$ denotes, as usual, the localization of $R$ with
respect to $R-\frak p$. Now, to see that the two rings are the same,
 it suffices
to use the following obvious property of formal power series rings:
if $A$ is a commutative ring and $f(t)\in A[[t]]$ is such that $f(0)$
is invertible in $A$, then $f(t)$ is invertible in $A[[t]]$.

The second equality is obvious: the stalk of
$p_*\Cal O_{\op{Spec}\, R[[t]]}$ at $\frak p$ is immediately
seen to coincide with $R[[t]]\otimes_R R_{\frak p}$. \qed

\vskip .1cm

Now,  composing with $i$ defines a map of sets
$$\phi: \op{Hom}_{\Schb}(\op{Spec}\, R[[t]], X)
\to
\op{Hom}_{\Lrsb}\bigl( (\op{Spec}\, R, \Cal O_{\op{Spec}\, R}[[t]]), X
\bigr).$$
A map $\psi$ in the other direction comes from the second equality
in (2.2.3). One verifies easily that $\phi$ and
$\psi$ are mutually inverse. This concludes the proof of (2.2.1)(b).
The rest of (2.2.1) is proved in {\it loc. cit.}

\vskip .2cm

We will also need the following generalization of Proposition 2.2.1(d).

\proclaim{(2.2.4) Proposition} 
Let $\phi: X\to Y$ be an \'etale morphism of schemes. Then

(a)
each morphism $\Lc^0_n(\phi)\,:\,\Lc^0_n(X)\to\Lc^0_n(Y)$ is \'etale,

(b)
the square
$$\matrix 
\Lc^0_n(X)&\to&\Lc^0_n(Y)\cr
\downarrow&&\downarrow\cr
X&\to&Y,
\endmatrix$$
as well as the analogous square for $\Lc^0(X)$, $\Lc^0(Y)$, is Cartesian.

\endproclaim

\noindent {\sl Proof:} We can assume that $X=\Spec(A), Y=\Spec(B)$
are affine.
It is enough to prove that
for each $n\geq 0$ the natural morphism
$$\alpha: \Lc^0_n(X)\to \Lc^0_n(Y) \times_{Y} X $$
is an isomorphism. Let us construct the inverse morphism $\beta$.
  Let $R$ be a $k$-algebra and $f$ be a morphism
$S\to \Lc^0_n(Y) \times_{Y} X$. Thus, $f$ corresponds to
a pair of ring homomorphisms forming the horizontal arrows of
the commutative diagram
$$\matrix
&B&\buildrel u\over\longrightarrow&R[t]/t^{n+1}&\cr
&{\ss\phi^*}\downarrow&&\downarrow{\ss\pi}&\cr
&A&\buildrel v\over\longrightarrow&R.
\endmatrix$$
Here  $\pi$ is the natural projection.
Since $\phi$ is \'etale  and $\pi$ is nilpotent, there is a unique
homomorphism $w: A\to R[t]/t^{n+1}$ such that both resulting triangles
are commutative. Let $g$ be the morphism $S\to \Lc^0_n(X)$
represented by $w$. Then we set $\beta(f)=g$. The verifications
are obvious and left to the reader.
\qed

\vskip .3cm

\noindent {\bf (2.3)  Nil-Laurent series.} Let $R$ be a commutative ring.
A nil-Laurent series is, by definition, a Laurent series $a(t) = \sum_{i\gg
-\infty} a_it^i\in R((t))$ such that all the $a_i$ with $i<0$ are nilpotent.
The set of such series will be denoted $R((t))^{\!\sqrt{\,\,\,}}$.

It is clear that $R((t))^\sqr$ is a subring in $R((t))$, the ring
of all Laurent series. Indeed, let
$\sqrt R$ be the radical of $R$ (the set of all nilpotent elements)
and set $R_\red=R/\sqrt R$.
Consider the homomorphism $\rho\,:\,R((t))\to R_\red((t))$
induced by the projection $R\to R_\red$.
Then $R((t))^\sqr=\rho^{-1}(R_\red[[t]])$.

\proclaim{(2.3.1) Proposition} For $a = \sum a_it^i \in R((t))^\sqr$
the following are equivalent :

(i) the element $a$ is invertible in  $R((t))^\sqr$,

(ii) the element $\rho(a)$ is invertible in $R_\red[[t]]$,

(iii) the element $a_0$ is invertible in $R$.
\endproclaim

\noindent {\sl Proof:} (i)$\Rightarrow$(ii) is obvious. To see that
(ii)$\Leftrightarrow$(iii), note that the invertibility of $\rho(a)$ implies
the invertibility of the image of $a_0$ in $R_\red$ which certainly implies
the invertibility of $a_0$ and the converse is equally obvious. Let us
prove that (ii)$\Rightarrow$(i).
Indeed, if $\rho(a)$ is invertible in $R_\red[[t]]$, we have $ab=1+c$
for some $b\in R((t))^\sqr$, $c\in\Ker\rho=\sqrt{R}((t))$. But every
element $c\in\sqrt{R}((t))$ is topologically nilpotent
(setting $c=c_-+c_+$ with $c_-\in t^{-1}\sqrt{R}[t^{-1}]$ and
$c_+\in \sqrt{R}[[t]]$, we have that $c_-$ is nilpotent while
$c_+$ is topologically nilpotent),
thus $ab$ is invertible in $R((t))^\sqr$, and $a$ is invertible, too.
\qed

\proclaim {(2.3.2) Corollary}
If $R$ is a local ring with maximal ideal $M$, then $R((t))^\sqr$ is
a local ring with maximal ideal
$\rho^{-1}\bigl( M_\red [[t]] + t R_\red[[t]]\bigr)$.
\endproclaim

\vskip .3cm

\noindent {\bf (2.4) The formal loop space.}
We now describe our main construction.
Let $X$ be a scheme of finite type
over $k$. Define a covariant functor $\lambda_X$ from
$\Algb$ to sets as follows:
$$\lambda_X (R) = \op{Hom}_{\Schb} \bigl( \op{Spec}\, R((t))^\sqr, X\bigr).
\leqno (2.4.1)$$

\proclaim {(2.4.2) Theorem} (a) The functor $\lambda_X$
is represented by an ind-scheme
$\Lc(X)$, containing $\Lc^0(X)$.

(b) $\Lc(X)$ is an inductive limit of nilpotent extensions of $\Lc^0(X)$.
In particular, for any open set $Y\subset \Lc^0(X)$ there is a well-defined
ind-scheme $\Lc(X)|_Y$.

(c) If $U$ is an open subset in $X$, then the ind-scheme $\Lc(U)$
is identified with $\Lc(X)|_{\pi^{-1}(U)}$.

(d) If $X$ is smooth, then $\Lc(X)$ is formally smooth.
\endproclaim

Note that Theorem 2.4.2 is closely related to \cite{BD1, Proposition 3.9.3.$(i)$}.
The proof will be finished in the next subsection.

\vskip .1cm

Let $R$ be a commutative ring. Since $R((t))^\sqr$ is a nilpotent extension
of $R[[t]]$, we have $\underline{\op{Spec}}\, R((t))^\sqr
=\underline {\op{Spec}}\, R[[t]]$, so $\op{Spec}\, R((t))^\sqr$
is the ringed space formed by $\underline {\op{Spec}}\, R[[t]]$
and a certain sheaf of rings
$\Cal O_{\op{Spec}\, R((t))^\sqr}$ on it.

\proclaim{(2.4.4) Lemma}(a)  We have, with
respect to the maps in (2.2.2), the equalities
$$ i^{-1} \Cal O_{\op{Spec}\, R((t))^\sqr} = \Cal O_{\op{Spec}\, R}((t))^\sqr =
p_*\Cal O_{\op{Spec}\, R((t))^\sqr}.$$

(b) For any scheme $S$ the sheaf $\Cal O_S((t))^\sqr$ is a sheaf of local
rings

\endproclaim

The proof of (a)
is analogous to that of Lemma 2.2.3. Instead of the property of
$A[[t]]$ quoted there, we use Proposition 2.3.1. Part (b) follows
from (a). \qed

\vskip .1cm

Let $\psi\,:\,\Schb\to\Lrsb$ be the functor such that
$S\mapsto (S, \Cal O_S((t))^\sqr)$.
Let us define a contravariant functor $\lambda'_X$ on the  category
$\Schb$ by
$$\lambda'_X(S)=\op{Hom}_{\Lrsb}\bigl(\psi(S),X\bigr).\leqno(2.4.3)$$

\proclaim{(2.4.5) Proposition}
For an affine scheme $S=\op{Spec}(R)$ we have $\lambda'_X(S) = \lambda_X(R)$.
\endproclaim

\noindent {\sl Proof:} Follows from Lemma 2.4.4 similarly to
Proposition 2.2.1(a).  \qed

\vskip .1cm

In virtue of Proposition 2.4.5, for the proof of Theorem 2.4.2
it suffices to show that
the functor $\l'_X$ on $\Schb$ is ind-representable.
We start by estabilishing some of its properties.

\proclaim {(2.4.6) Proposition} (a) For every scheme $X$ the functor
$\l'_X$ is a sheaf on $\Schb$.

(b) If $U\subset X$ is an open subset, then the induced
morphism of functors $\l'_U\to\l'_X$ is open.

(c) Let $(U_\a)_{\a\in A}$ be an open covering of $X$. Then
$\l'_X$ is equal to the cokernel, in the category $\Shfb$,
of the pair of morphisms
$$\coprod_{\a,\b}\l'_{U_\a\cap U_\b}\rra\coprod_\a\l'_{U_\a},$$
\endproclaim

\noindent{\sl Proof:}
The proposition follows from simple properties of representable
functors on the category  $\Lrsb$.
If $\Cal T = (T, \Oc_\Tc)$
is a   locally ringed space, we will call an open part of $\Cal T$
a ringed space of the form $(U, \Oc_\Tc|_U)$ where $U\subset T$
is an open subset in the usual sense.
An open embedding is, by definition,  a morphism isomorphic
to the inclusion of an open part.
Accordingly, we have
the concept of an open covering of $\Cal T$. This
makes $\Lrsb$ into a Grothendieck site.
For $\Cal X\in\Lrsb$ let $\eta_{\Cal X}$
be the corresponding representable functor on  $\Lrsb$.
As in (2.1), a morphism
$F: \phi'\to\phi$ of functors $\Lrsb\to \bold{Sets}$ will be called
open, if for any $\Cal S\in \Lrsb$ and any element $u\in \phi(\Cal S)$
(i.e., a morphism $\eta_{\Cal S}\to \phi$) the fiber product
functor $\eta_{\Cal S}\times_{\phi} \phi'$ is representable by a
locally ringed space $\Cal S'$ whose
natural morphism to $\Cal S$ is an open embedding.
Let us recall the following basic facts.

\proclaim {(2.4.7) Lemma}
(a) For any $\Cal X\in\Lrsb$ the representable
functor $\eta_{\Cal X}$ is a sheaf on $\Lrsb$.

(b) If $\Cal U\subset\Cal X$ is an open embedding in $\Lrsb$,
then $\eta_{\Cal U}\to \eta_{\Cal X}$ is an open morphism of functors.

(c) Let $\Cal X\in\Lrsb$ and $(\Cal U_\a)_{\a\in A}$
be an open covering of $\Cal X$. Then in the category of sheaves
of sets on $\Lrsb$ we have the equality
$$\eta_{\Cal X}=\op{Coker}\Bigl\{\coprod_{\a,\b}\eta_{\Cal U_\a\cap
\Cal U_\b}\rra\coprod_\a\eta_{\Cal U_\a}\Bigr\}\leqno (2.4.8)$$
or, explicitly, for any $\Cal S\in\Lrsb$,
$$\eta_{\Cal X}(\Cal S) =\lim_{\buildrel\lra\over{\Cal S_\a}}\op{Ker}
\Bigl\{\prod_\a\eta_{\Cal U_\a}(\Cal S_\a)\rra
\prod_{\a\b}\eta_{\Cal U_{\a}\cap \Cal U_\b}(\Cal S_\a\cap \Cal S_\b)\Bigr\},
\leqno(2.4.9)$$
where the limit is taken over the set of
open coverings $(\Cal S_\a)_{\a\in A}$ of $\Cal S$
(the indexing set $A$ being fixed) ordered by simultaneous inclusion.
\endproclaim

\noindent
Proposition 2.4.6 follows from Lemma 2.4.7.
Indeed, we have $\l'_X=\eta_X\circ\psi$,
where $X$ is viewed as a locally ringed space.
Thus, to prove Proposition 2.4.6.$(b)$ using Lemma 2.4.7.$(b)$
it is enough to prove that for any open embedding of schemes
$U\hookrightarrow X$ and any $u\in\l'_X(S)$
there is an open embedding $j\,:\,S'\hookrightarrow S$
such that the following square is Cartesian
$$\matrix
U&\hookrightarrow&X\cr
\uparrow&&\uparrow u\cr
\psi(S')&{\buildrel\psi(j)\over\to}&\psi(S)
\endmatrix$$
(then, $\eta_S\times_{\l'_X}\l'_U=\eta_{S'}$).
This is obvious.
\qed

\vskip3mm

\noindent{\sl Proof of (2.4.7) :}
 Part $(a)$ is obvious, and
Part $(b)$ is proved in \cite{EGA0, (4.5.2)}.
Let us prove Part $(c)$.
For any $\phi\in\eta_\Xc(\Sc)$ and any $\a$,
consider the ringed space
$\Sc_\a=(\phi^{-1}(\Uc_\a ),\Oc_\Sc|_{\phi^{-1}(\Uc_\a)})$.
Clearly, $(\Sc_\a)$ is an open covering of $\Sc$.
Let $\phi_\a$ be the restriction of $\phi$ to $\Sc_\a$.
Then, $(\phi_\a)$ is an element of $\prod_\a\eta_{\Uc_\a}(\Sc_\a)$
lying in the kernel (2.4.9).
Thus we have constructed a map from
$\eta_\Xc(\Sc)$ to the right hand side of (2.4.9).
On the other hand, assume that $(\phi_\a)\in\prod_\a\eta_{\Uc_\a}(\Sc_\a)$
is such that $\phi_\a(\Sc_\a\cap\Sc_\b)\subset\Uc_\a\cap\Uc_\b$ and
${\phi_\a}|_{\Sc_\a\cap\Sc_\b}={\phi_\b}|_{\Sc_\a\cap\Sc_\b}$.
The corresponding maps $\Sc_{\a,\red}\to\Uc_{\a,\red}$ glue together,
giving a continuous map $\phi_\red\,:\,\Sc_\red\to\Xc_\red$.
Moreover, there is an obvious sheaf homomorphism
$\phi^*\,:\,\phi_\red^{-1}\Oc_\Xc\to\Oc_\Sc$ :
the restriction of $\phi^*$ to $\Sc_\a$
is the composition of the chain of maps
$$(\phi_\red^{-1}\Oc_\Xc)|_{\Sc_\a}=(\phi_{\a,\red}^{-1}\Oc_{\Uc_\a})|_{\Sc_\a}
{\buildrel\phi_\a^*\over\to}\Oc_{\Sc_\a}=\Oc_\Sc|_{\Sc_\a}.$$
This establishes (2.4.9).

 Let $C$ be the cokernel as in (2.4.8) but taken in the category
of presheaves. Then the cokernel in the category of sheaves
is just the sheaf $\overline {C}$ associated to the presheaf $C$. By definition,
for $\Sc\in\Lrsb$ the set $C(\Sc)$ consists of pairs $(\a, \phi: S\to
\Uc_\a)$ taken modulo the identifications coming from morphisms
of $\Sc$ into $\Uc_\a\cap\Uc_\b$.  Now,  by the definition of the sheaf
associated to the presheaf,
$$ \overline{C} = \ind_{\Sc = \bigcup_{\gamma\in\Gamma} \Sc_\gamma}
\op{Ker} \biggl\{ \prod_{\gamma\in\Gamma} C(\Sc_\gamma)\to \prod_{
\gamma, \gamma'\in\Gamma} C(\Sc_\gamma\cap\Sc_{\gamma'})\biggr\},$$
where the limit is over all possible open coverings of $\Sc$
(with arbitrary indexing sets $\Gamma$) ordered by refinement.
Notice now that specifying a section of $C$ over $\Sc_\gamma$
includes specifying an index $\alpha$ from the set $A$ indexing
the covering $\{\Uc_\a\}$, so we get a map $p:\Gamma\to A$.
Denoting $\Sc_\a = \bigcup_{p(\gamma)=\a}\Sc_\gamma$, we get
an element of the right hand side of (2.4.9). This establishes
the equivalence of (2.4.8) and (2.4.9). \qed

\vskip3mm

Note that Proposition 2.4.6(c) implies the following.

\proclaim {(2.4.10) Corollary}
We have $\l_X = \lim\limits_{\lra}{}_{
U\subset X \,\,\text{affine}}\l_U$, the limit being taken
in the category $\Shfb$.
\endproclaim

\vskip .3cm

\noindent {\bf (2.5) Proof of Theorem 2.4.2.}
We first assume that $X=\Spec A$ is an affine scheme of finite type.
Consider the larger functor
$\tilde\l_X$ on $\Algb$ defined by
$$\tilde\l_X(R) = \op{Hom}_{\bold{Alg}} \bigl(A,R((t))\bigr).$$

\proclaim {(2.5.1) Proposition} (a)
The functor $\tilde\l_X$ is represented by
an ind-scheme  $\tilde\Lc(X)$, which is an inductive limit of
affine schemes of infinite type.

(b) If $X$ is smooth, then $\tilde\Lc(X)$ is formally smooth.
\endproclaim

\noindent {\sl Proof:} (a) Consider the $k$-ind-scheme
$$\underline{k((t))}=\ind_N\op{Spec}k[a_l;-N\leq l].$$
We can think of the $a_i$ as the coefficients of an indeterminate
Laurent series $\sum a_it^i$. It is clear that $\underline{k((t))}$
represents the functor $\tilde\l_{\AA^1}$. Since $\AA^1$ is
a $k$-algebra object in the category of schemes, multiplication
of Laurent series makes $\underline {k((t))}$ into a $k$-algebra
object in the category of ind-schemes. In particular, each polynomial
$f\in k[x_1,...,x_d]$ defines a morphism of ind-schemes
$\tilde\Lc(f): \underline{k((t))}^d \to \underline{k((t))}$.
Hereafter we will write $k[x]$ instead of $k[x_1,...,x_d]$ to simplify.

Suppose now that $X$ is given in $\AA^d$
by a system of algebraic equations, say
$$f_m(x_1,...,x_d)=0,\qquad m=1,2,...,l.$$
The ind-scheme $\tilde\Lc(X)$ is realized as the closed
sub-ind-scheme in $\underline{k((t))}^d$ defined as the intersection
of the preimages of 0 under the $\tilde\Lc(f_m)$. More explicitly,
replacing $x_i$ by $x_i(t)=\sum_{j\geq -N}a_j^{(i)}t^j$ in the equations
above, we get, for each $N$, a system of algebraic equations in
$k[a_l^{(i)};-N\leq l]$ for each $N$ which defines a subscheme
in $(\Spec k[a_l;-N\leq l])^d$. Our desired ind-scheme
$\tilde\Lc(X)$ is the direct limit of these schemes as $N\to\infty$.

(b) The infinitesimal lifting
condition for $\tilde\l_X$ is formulated for affine schemes
$S = \op{Spec}(R)$.
We need to prove that for any surjection of $k$-algebras $R'\to R$
whose kernel $I$ satisfies $I^n=0$ for some $n$, the map of sets
$\tilde\l_X(R')\to \tilde\l_X(R)$ is surjective.
But the kernel of $R'((t))\to R((t))$ is $I((t))$ which
is also nilpotent of order $n$. So the smoothness of $A$
implies that any morphism $A\to R((t))$ can be lifted to
a morphism $A\to R'((t))$.
\qed

\vskip .2cm

\proclaim {(2.5.2) Corollary} The functor $\l_X$ is represented
by an ind-scheme $\Lc(X)$ which is the inductive limit of
the formal neighborhoods of $\Lc^0(X)$ in the schemes forming an inductive
system for $\tilde\Lc(X)$. If $X$ is smooth, then $\Lc(X)$
is formally smooth.
\endproclaim

This proves parts (a), (b) and (d) of Theorem 2.4.2 for the case of affine $X$.
Part (c) of the theorem follows from Proposition 2.4.6(c).

\vskip .2cm

To prove Theorem 2.4.2 for general $X$, it is enough  to establish
the existence of the limit
$\lim\limits_{\lra}{}_{
U\subset X \,\,\text{affine}}\Lc(U)$ in the category of ind-schemes.
In fact, it is enough to take the limit over a finite set of $U$
consisting of elements of some finite covering and their
intersections.
Indeed, given this, all the other properties follow from the affine case
and from Proposition 2.4.6. But using Proposition 2.4.6(b) again, we see that
for any pair  $U'\subset U$ of affine open subsets in $X$ the
corresponding morphism of ind-schemes $\Lc(U')\to \Lc(U)$ is an
open embedding. Further, the ind-scheme $\Lc(U)_\red$
is actually a scheme, namely $\Lc^0(U)$. So our statement
follows from the next general fact.

\proclaim {(2.5.3) Lemma} Let $(Z_i)_{i\in I}$ be a finite diagram
of ind-schemes in which all the arrows are open embeddings.
Assume that the inductive limit of $Z_{i,\red}$ exists
in the category of ind-schemes. Then so does the inductive limit of the $Z_i$.
\endproclaim

Theorem 2.4.2 is proved. Let us also note the following fact.

\proclaim {(2.5.4) Proposition}
If $\phi: X\to Y$ is an \'etale morphism of schemes of finite type, then
$\Lc(\phi): \Lc(X)\to\Lc(Y)$ is formally \'etale.
\endproclaim

\noindent {\sl Proof:} By the above we can assume that $X=\Spec(A)$,
$Y=\Spec(B)$ are affine.
If $S$ is a scheme and we have two compatible maps
$$\a\,:\,S_\red\to\Lc(X),\qquad\b\,:\,S\to\Lc(Y),$$
we must construct the unique map $S\to\Lc(X)$.
If $S=\Spec R$, then we have a diagram
$$\matrix
A&{\buildrel\a^*\over\to}&R_\red((t))^\sqr&=&R_\red[[t]]\cr
{\ss\phi^*}\uparrow&&\uparrow&&\cr
B&{\buildrel\b^*\over\to}&R((t))^\sqr,&&
\endmatrix$$
where the right arrow, which is the projection $\rho$ in (2.3),
has nilpotent kernel. Therefore there is a unique homomorphism
of rings $A\to R((t))^\sqr$ making the diagram commute.
\qed

\vskip .3cm

\noindent {\bf (2.6) The formal loop space as an ind-object.}
  For future
purposes we construct a certain class of local
realizations of $\Lc(X)$ as  an ind-object in the category
of schemes.

Let $x_1,x_2,...,x_d$ be the coordinates on  the affine space $\AA^d$.
For any $k$-algebra $A$ and any morphism $\phi\,:\,\Spec A\to\AA^d$
let $\phi^*$ be the corresponding map $k[x_1, ..., x_d]\to A$.

\proclaim{(2.6.1) Proposition}
Assume that $X = \Spec A$ is an affine scheme and that $\phi\,:\,X\to\AA^d$
is an \'etale morphism of schemes.

(a) There is a unique morphism $\theta_\phi\,:\,\Lc(X)\to X$
such that
$\theta_\phi(f)(\phi^*x_i)=f(\phi^*x_i)_0$ for any $f\in\l_X(R)$.
The restriction of $\theta_\phi$ to $\Lc_0(X)$
is equal to the projection $p$ from (2.2.1)(b).

(b) We have $\Lc(\AA^d)\times_{\AA^d}X\simeq\Lc(X)$.
\endproclaim

\noindent{\sl Proof:}  Consider the morphism of functors
$\tilde\l_{\AA^d}\to\eta_{\AA^d}$ which maps a homomorphism
$f\in\op{Hom}_{\bold{Alg}}\bigl(k[x_1,...,x_d],R((t))\bigr)$
to the unique morphism in
$\Hom_{\Algb}(k[x_1,...,x_d],R)$ such that $x_i\mapsto f(x_i)_0$.
This morphism of functors can be seen as a morphism of ind-schemes
$\tilde\theta_{\AA^d}\,:\,\tilde\Lc(\AA^d)\to\AA^d$.
It is clear that its restriction to $\Lc^0(\AA^d)$ coincides with
the projection $p$ from (2.2.1)(b).
Let $\theta_{\AA^d}$ be the restriction of $\tilde\theta_{\AA^d}$
to the ind-subscheme $\Lc(\AA^d)\subset\tilde\Lc(\AA^d)$.
Since $\Lc(X)$ is an inductive limit of nilpotent extensions
of the scheme $\Lc^0(X)$ and $\phi$ is formally \'etale,
there is a map $\theta_\phi\,:\,\Lc(X)\to X$ splitting the diagram 
$$\matrix
\Lc^0(X)&\to&X\cr
\downarrow&&\downarrow{\ss\phi}\cr
\Lc(X)&\to&\AA^d
\endmatrix$$
into two commutative triangles.
Here, the lower horizontal  arrow is the composition of maps
$$\theta_{\AA^d}\circ\Lc(\phi)\,:\,\Lc(X)\to\Lc(\AA^d)\to\AA^d.$$
Let $\psi$ be the resulting map
$\Lc(X)\to\Lc(\AA^d)\times_{\AA^d}X$.
We claim that  $\psi$ is an isomorphism.
For this, we construct a map $\chi: \Lc(\AA^d)\times_{\AA^d}X\to\Lc(X)$
inverse to $\psi$.
Let $S=\Spec R$ be a scheme. A morphism
$S\to\Lc(\AA^d)\times_{\AA^d}X$ is a compatible pair
$$\bigl(\a\,:\,S\to\Lc(\AA^d),\,\b\,:\,S\to X\bigr).$$
We need to construct a map $\chi(\alpha,\beta): S\to\Lc(X)$.
First, $\a(S_\red)\subset\Lc^0(\AA^d)$, and $\b(S_\red)\subset X$.
Thus, by Proposition 2.2.4, we have a map $\g\,:\,S_\red\to\Lc^0(X)$. The composition
$\tilde\g$ of $\g$ and the embedding $\Lc^0(X)\subset\Lc(X)$
gives a commutative diagram
$$\matrix
S_\red&{\buildrel\tilde\g\over\to}&\Lc(X)&\cr
\downarrow&&\downarrow{\ss\Lc(\phi)}&\cr
S&{\buildrel\a\over\to}&\Lc(\AA^d).&
\endmatrix\leqno(2.6.2)$$
To complete the construction, we
notice that $\Lc(\phi)$ is formally \'etale by Proposition 2.5.4 and
so we have a map $S\to\Lc(X)$ splitting (2.6.2) into
two commutative triangles. We take this map to be $\chi(\alpha,\beta)$.
The verification that $\chi$ is inverse to $\psi$ is straightforward.
Proposition 2.6.1 is proved.
\qed

\vskip3mm

Let $\EE$ be the set of $\eps = (\eps_{-1}, \eps_{-2}, ...)$,
$\eps_j\in \ZZ_+$ such that $\eps_j=0$ for almost all $j$.
It is equipped with the natural partial order: $\eps\leq \eps'$
if $\eps_j\leq \eps'_j$ for all $j$.
In the remainder of this section we assume that the $k$-scheme $X$ is smooth.
Thus $X$ can be covered by affine open sets
$U=\Spec A$ possessing \'etale maps $\phi\,:\,U\to\AA^d$.
For every such $U, \phi$ we consider the functor
$$\l^\eps_\phi\,:\,R\mapsto\bigl\{f\in\l_U(R)\,|\,
f(\phi^*x_i)_j^{1+\eps_j}=0,\,\forall i\in[1,d],\,\forall j<0\bigr\}.$$

\proclaim{(2.6.3) Proposition}
The functor $\l^\eps_\phi$ is representable
by a closed subscheme $\Lc^\eps(\phi)\subset\Lc(U)$,
such that $\Lc(U)$ is the inductive limit of the schemes $\Lc^\eps(\phi)$.
\endproclaim

\noindent{\sl Proof:} We first consider the case $U=\AA^d, \phi = \op{Id}$.
 Let $N$ be such that $\epsilon_j=0$ for
$j<N$. Define the scheme $\Lc^\eps(\AA^d)$ by
$$\Lc^\eps(\AA^d)=\Spec\bigl(k[a_l^{(i)}; -N\leq l]\bigl/
\bigl((a^{(i)}_l)^{1+\epsilon_l};l<0\bigr)\bigr).
\leqno (2.6.4)$$
It is immediate that this scheme represents
the functor $\l^\eps_\Id$,
and that $\l^\eps_\phi=\l^\eps_\Id\times_{\l_{\AA^d}}\l_U$.
The proof of Proposition 2.6.1 implies that that the map
$$\l^\eps_{\phi}\to\l^\eps_\Id\times_{\eta_{\AA^d}}\eta_U,\quad
f\mapsto\bigl(f\circ\phi^*,\theta_\phi(f)\bigr),\leqno(2.6.5)$$
is an isomorphism of functors.
Thus, the closed subscheme
$$\Lc^\eps(\AA^d)\times_{\AA^d}U\subset\Lc(\AA^d)\times_{\AA^d}U=\Lc(U),
\leqno (2.6.6)$$
see Proposition 2.6.1.$(b)$,
represents the functor $\l^\eps_{\phi}$.
It is clear by the definition
that $\l_U=\ind_\eps\lambda^\eps_\phi$ in the category
of functors $\bold{Alg}\to\bold{Sets}$ and thus we have that
$\Lc(U)=\ind_\eps\Lc^\eps(\phi)$ in the category of ind-schemes.
\qed

\vskip3mm

Let us note the following reformulation of this fact, to be used later.

\proclaim {(2.6.7) Corollary}
If $\phi, \psi$ are two \'etale maps $U\to \AA^d$, then the ind-objects
$``\ind_\eps" \Lc^\eps(\phi)$, $``\ind_\eps"\Lc^\eps(\psi)$ are isomorphic, i.e.,
for any $\eps$ there is $\eps'$ such that $\Lc^\eps(\phi)\subset
\Lc^{\eps'}(\psi)$ and vice versa.
\endproclaim

\noindent{\sl Proof:}
Given $\psi$, any map $S\to\Lc(U)$
with $S$ an affine scheme, factors through some $\Lc^{\eps'}(\psi)$.
Now take $S = \Lc^\eps(\phi)$ which we know to be an affine
scheme by (2.6.4,6).
\qed

\vskip 3mm

\noindent {\bf (2.7) The formal loop space as an ind-pro-object.}
We keep the notations of (2.6). Thus
$U=\Spec(A)$ is affine and $\phi: U\to \AA^d$ is \'etale.
The schemes $\Lc^\eps(\phi)$ are of infinite type.
In fact, each of these schemes is a projective limit of schemes
of finite type, so $\Lc(U)$ can be viewed as an ind-pro-object
in the category of schemes of finite type over $k$. In this subsection we
construct explicit ind-pro-systems for $\Lc(U)$.

Consider the functor
$\l^\eps_{n\phi}\,:\,R\mapsto\l^\eps_\phi(R)/\sim_{n\phi},$ where
$$f\,\sim_{n\phi}\,g\,\iff\,f(\phi^*x_i)-g(\phi^*x_i)\in t^{n+1}R[[t]],\quad
\forall i.$$

\proclaim{(2.7.1) Proposition}
(a) The functor $\l^\eps_{n\phi}$ is representable by a scheme
$\Lc^\eps_n(\phi)$.
The scheme $\Lc^\eps_n(\phi)$ is of finite type and is a nilpotent extension of
$\Lc^0_n(U)=\Lc^0_n(\phi)$.
Moreover, $\Lc^\eps(\phi)=\pro_n\Lc^\eps_n(\phi)$.

(b) The schemes $\Lc^\eps_n(\phi)$ form an ind-pro-system with
Cartesian squares ($n'\geq n,\ \eps'\geq\eps$)
$$\matrix
\Lc^{\eps}_{n'}(\phi)&\hookrightarrow&\Lc^{\eps'}_{n'}(\phi)\cr
\downarrow&&\downarrow\cr
\Lc^{\eps}_{n}(\phi)&\hookrightarrow&\Lc^{\eps'}_{n}(\phi)
\endmatrix$$
where the vertical arrows are smooth affine morphisms.

(c) The ind-pro-object
${\displaystyle``\ind_\eps" ``\pro_n"\Lc^\eps_n(\phi)}$
is independent, up to isomorphism,  on $\phi$.
\endproclaim

\noindent{\sl Proof:}
Claim $(a)$ is entirely similar to Proposition 2.6.3.
We first consider the case $U=\AA^d,\phi=\Id$ and
define the scheme $\Lc^\eps_n(\AA^d)$ by
$$\Lc^\eps_n(\AA^d)=\Spec\bigl(k[a_l^{(i)};-N\leq l\leq n]\bigl/
\bigl((a^{(i)}_l)^{1+\eps_l};l<0\bigr)\bigr), \quad N\gg 0.\leqno(2.7.2)$$
It represents $\lambda^\eps_{n\Id}$. 
The fiber product scheme
$\Lc^\eps_n(\AA^d)\times_{\AA^d}U$
represents the functor 
$\lambda^\eps_{n\Id}\times_{\eta_{\AA^d}}\eta_U$.
The isomorphism of functors (2.6.5) yields an isomorphism of functors
$\lambda^\eps_{n\phi}\to\lambda^\eps_{n\Id}\times_{\eta_{\AA^d}}\eta_U$.
Thus
$$\Lc^\eps_n(\phi)=\Lc^\eps_n(\AA^d)\times_{\AA^d}U.\leqno(2.7.3)$$

Claim $(b)$ is obvious in the case $U=\AA^d$, $\phi=\Id$. 
The general case follows from (2.7.3) 
since the base change of a smooth affine morphism is still smooth affine,
and the base change of a Cartesian square is Cartesian.

Claim $(c)$ follows from Corollary 2.6.7 and Proposition 2.1.2.
\qed

\vskip3mm

Passing to the limit in $\epsilon$ we get the ind-scheme
$\Lc_n(\phi) = \ind_\epsilon \Lc^\epsilon_n(\phi)$. It represents the functor
$R\mapsto \l_U(R)/\sim_{n\phi}$. As in (2.7.2) we get
$$\Lc_n(\AA^d) = \ind_N \op{Spf}\, k[a^{(i)}_l, 0\leq l\leq n]
[[ a^{(i)}_l, -N\leq l<0]],$$
and we see that $\Lc_n(\AA^d)$ is formally smooth. 
Further, we see that
$$\Lc_n(\phi) = \Lc_n(\AA^d)\times_{\AA^d} U,$$
from which we see that $\Lc_n(\phi)$ is also formally smooth.

\vfill\eject

\head 3. The global loop space\endhead

\vskip 3mm

\noindent {\bf (3.1) Localization on a curve.}
Consider the Lie algebra $\op{Der}\, k[[t]]$ and the group scheme
$\op{Aut} \, k[[t]]$ over $k$. They form a Harish-Chandra pair
\cite {BB} \cite {F, \S 6.1}. By construction, we have the action
of this pair on the scheme $\Lc^0(X)$ and on the ind-scheme $\Lc(X)$
constructed in \S 2.

Let $C$ be a smooth curve over $k$. The well known procedure of localization,
see {\it loc. cit.}  and \cite {GKF},  gives then a scheme $\Lc^0(X)_C$
and an ind-scheme $\Lc(X)_C$ over $C$ defined as follows. Let $C^\wedge
\to C$ be the scheme of pairs $(c, t_c)$ where $c\in C$ and $t_c$
is a formal coordinate near $c$. Then $C^\wedge$ has a natural
$(\op{Der} \, k[[t]], \,\,\, \op{Aut}\, k[[t]])$-structure, i.e. a simply
transitive $\op{Der}\, k[[t]]$-action extending the
fiberwise $\op{Aut} \, k[[t]]$-action. We define
$$\Lc^0(X)_C=C^\wedge\times_{\op{Aut}\, k[[t]]} \Lc^0(X), \quad
\Lc(X)_C=C^\wedge\times_{\op{Aut}\, k[[t]]} \Lc(X). \leqno (3.1.1)$$
If $X$ is affine, we define, in a similar way, the ind-scheme $\tilde\Lc(X)_C$
starting from $\tilde\Lc(X)$.
Because of the simple transitivity of the $\op{Der}\, k[[t]]$-action on
$C^\wedge$, the (ind-)schemes thus constructued possess natural
connections along $C$ which are compatible with the embeddings $\Lc^0(X)_C
\subset \Lc(X)_C$ and, for $X$ affine, $\Lc^0(X)_C\subset \tilde\Lc(X)_C$.

Note that $\Lc^0(X)_C$ is nothing but the scheme of infinite jets of morphisms
$C\to X$, so it is easy to describe explicitly the functor represented
by $\Lc^0(X)_C$ (and also  by $\Lc(X)_C$, $\tilde\Lc(X)_C$).
Let us introduce some notations.
For a scheme $S$ and a morphism $f: S\to C$ we denote by $\Gamma(f)
\subset S\times C$ the graph of $f$. Let $f_\red: S_\red\to
C$ be the restriction of $f$ to $S_\red$. We have then the following
sheaves of $k$-algebras on $S\times C$ supported on $\Gamma(f)$:

\itemitem{--}
$\Oc^\wedge_f$, the completion of $\Cal O_{S\times C}$
along $\Gamma(f)$, i.e. the sheaf of functions on the formal neighborhood
of $\Gamma(f)$.

\itemitem{--}
$\Kc^\wedge_f$, the sheaf of functions on the punctured formal neighborhood
of $\Gamma(f)$. Thus $\Kc^\wedge_f$
is obtained from $\Oc^\wedge_f$ by inverting the local equations
of the divisor $\Gamma(f) \subset S\times C$.

\itemitem{--}
$\Kc_f^\sqr\subset\Kc^\wedge_f$, the subsheaf of sections whose
restriction to $S_\red\times C$ lies in
$\Oc^\wedge_{f_\red}\subset\Kc^\wedge_{f_\red}.$

\vskip .2cm

\proclaim{(3.1.2) Proposition} (a) The sheaves
$\Oc^\wedge_f$ and $\Kc_f^\sqr$ are sheaves of local rings.

(b) The scheme $\Lc^0(X)_C$ represents the functor $\l^0_{X,C}: \Schb
\to\bold{Sets}$ which associates to a scheme $S$ the set of pairs
$(f,\rho)$ where $f: S\to C$ is a morphism of schemes and $\rho:
(\Gamma(f),\Oc^\wedge_f)\to X$ is a morphism of locally ringed spaces.

(c) The ind-scheme $\Lc(X)_C$ ind-represents the functor $\l_{X,C}$
defined as in (b) but with $\Oc^\wedge_f$ replaced by $\Kc_f^\sqr$.

(d) Similarly, when $X$ is affine,  the ind-scheme
$\tilde\Lc(X)_C$ ind-represents the functor $\tilde\l_{X,C}$
defined as in (b) but with $\Oc^\wedge_f$ replaced by $\Kc^\wedge_f$.
\endproclaim

\noindent {\sl Proof:} (a) By choosing an etale coordinate
$y$ on $C$ and using the relative etale coordinate $y-f(s)$
on $C\times S$, we reduce to the case when $C=\AA^1$ and $f$
is constant with value $0\in\AA^1$. Then
$$\Oc^\wedge_f = \Oc_S[[t]], \quad \Kc_f = \Oc_S((t)), \quad
\Kc_f^\sqr = \Oc_S((t))^\sqr$$
and our assertion follows from (2.2.1)(a) and (2.4.4)(b).

\vskip .1cm

(b) The projection $\pi: C^\wedge\to C$ induces, for any scheme $S$, a
map of sets
$$\pi_S: \op{Hom}(S, \Lc^0(X)_C)\to \op{Hom}(S, C).$$
It is enough to show that for any $f: S\to C$ the set $\pi_S^{-1}(f)$
is naturally identified with the set of $\rho: (\Gamma(f), \Oc^\wedge_f)
\to X$. Further, since both functors $\lambda^0_{X, C}$ and
$\eta_{\Lc^0(X)_C}$ are sheaves of sets on $\Schb$, it is enough to
construct such an identification Zariski locally on $S$.
But locally on $S$ we have, from the definition (3.1.1):
$$\pi_S^{-1}(f) = \op{Lifts}(f, C^\wedge) \times_{\op{Aut} \, k[[t]]}
\op{Hom}(S, \Lc^0(X)),$$
where $\op{Lifts}(f, C^\wedge)$ is the set of $\tilde f: S\to C^\wedge$
such that $\pi\tilde f=f$. Recall that
$$\op{Hom}(S, \Lc^0(X)) = \op{Hom}_{\Lrsb} ((S, \Oc_S[[t]]), X).$$
This means that
$$\pi_S^{-1}(f) = \op{Hom}_{\Lrsb}((S, \Ac), X),$$
where $\Ac$ is the sheaf of rings on $S$ associated to the presheaf
$$U\mapsto \biggl( \prod_{\tilde f\in\op{Lifts}(f|_U, C^\wedge)}
\Oc_U[[t]]\biggr)^{\op{Aut} k[[t]]}.$$
But it is clear that $\Ac\simeq \Oc^\wedge_f$, whence the statement.
Parts (c) and (d)  are  proved similarly.

\qed

\vskip .3cm

\noindent {\bf (3.2) Factorization monoids.} Let $C$ be a smooth curve, as
before. For any surjection $J\twoheadrightarrow I$ of finite sets
and $i\in I$ we denote by $J_i$ the preimage of $i$. To such a surjection
one associates, in a standard way,
 the  ``diagonal"
embedding $\Delta^{(J/I)}\,:\,C^I\hookrightarrow C^J$.
Let $U^{(J/I)}\subset C^J$ be the locus of $(c_j)_{j\in J}$ such that
$c_j\neq c_{j'}$ whenever the images of $j$ and $j'$ in $I$ are different.
 We denote by $j^{(J/I)}\,:\,U^{(J/I)}\hookrightarrow C^J$  the embedding.

\proclaim {(3.2.1) Definition} Let $Y$ be an ind-scheme over $C$.
A factorization monoid structure on $Y$ is a collection of
ind-schemes $Y_I$ over $C^I$  with a formally integrable connection,
given for each nonempty finite set $I$
such that $Y_{\{1\}} = Y$ and $Y_I$ is formally smooth over $C^I$,
equipped with the following data:

(a) An isomorphism of $C^I$-ind-schemes
$\nu^{(J/I)}\,:\,\Delta^{(J/I)*} Y_J\simto Y_I$
for every surjection $J\twoheadrightarrow I$, compatible with  compositions
of surjections.

(b)  An isomorphism of  $U^{(J/I)}$-ind-schemes
$$\kappa^{(J/I)}\,:\,j^{(J/I)*}(\Prod_{i\in I} Y_{J_i})
\simto j^{(J/I)*} Y_{J}$$
for every $J\twoheadrightarrow I$, such that for $K\twoheadrightarrow J$
the isomorphism $\kappa^{(K/I)}$ coincides with the composition
$\kappa^{(J/I)}\circ(\prod_{i\in I}\kappa^{(K_i/J_i)})$.
We also demand that $\nu$, $\kappa$ are compatible in the following sense :
for every $J\twoheadrightarrow J'\twoheadrightarrow I$ one has
$\nu^{(J/J')}\circ\Delta^{(J/J')*}(\kappa^{(J/I)})=
\kappa^{(J'/I)}\circ(\boxtimes_{i\in I}\nu^{(J_i/J'_i)}).$

\endproclaim

\noindent {\bf (3.2.2) Remarks.} (a) This  is a nonlinear counterpart of the
concept of a factorization algebra due to Beilinson and Drinfeld
\cite {BD1} \cite{G}, see also (6.1) below.
factorization monoids  can be used to construct factorization algebras by
applying ``natural"
linear constructions.

(b) More generally, we can speak about a factorization monoid
structure on any functor $\Cal Y: \Schb\to\bold{Sets}$
(not necessarily representable by an ind-scheme) which is equipped with
a morphism to $C$ (i.e. to the representable functor $\eta_C$).

\vskip .3cm

\noindent {\bf (3.2.3) Example.} Let $G$ be an  affine algebraic group over
$k$. Then $\Lc^0(G)$ is a group scheme and $\tilde\Lc(G)$ is
a  group ind-scheme over $C$. The quotient ind-scheme
${\frak Gr}_G = \tilde\Lc(G)/\Lc^0(G)$ is known as the affine Grassmannian
associated to $G$. The natural family of such Grassmannians over $C$, i.e.
the ind-scheme ${\frak Gr}_G (C) = \tilde\Lc(G)_C/\Lc^0(G)_C$
is known to have a structure of a factorization monoid,
see \cite{G, \S 5.2.2}.

\vskip .3cm

Now, the main result of this section is the following.

\proclaim {(3.2.4) Theorem} Let $X$ be a scheme of finite type.
Then the $C$-scheme $\Lc^0(X)_C$ and the $C$-ind-scheme $\Lc(X)_C$
possess natural structures of factorization monoids so that
the embedding $\Lc^0(X)_C\subset\Lc(X)_C$ is a factorization
homomorphism. Similarly, if $X$ is affine, then the $C$-ind-scheme
$\tilde \Lc(X)_C$ has a natural structure of a factorization monoid
so that the embedding $\Lc^0(X)_C\subset\tilde\Lc(X)_C$ is a factorization
homomorphism.
\endproclaim

\noindent {\bf (3.3) Factorization monoid structure on the functors
represented by $\Lc^0(X)_C,\tilde\Lc(X)_C$ and $\Lc(X)_C$.}
Let $S$ be a scheme and $f_I$  be $I$-uple of morphisms
$f_i\,:\,S\to C$, $i\in I$. We denote by
$\Gamma(f_I)\subset S\times C$
the union of the graphs of the $f_i$'s
and by $f_{I,\red}\subset S_\red\times C$
the union of the graphs of the $f_{i,\red}$.
Let us introduce, similarly to (3.1),
the following sheaves of rings on $S\times C$ with support in
$\Gamma(f_I)$ :

\itemitem{--}
$\Oc^\wedge_{f_I}$, the completion of $\Cal O_{S\times C}$
along $\Gamma(f_I)$, i.e. the sheaf of functions on the formal neighborhood
of $\Gamma(f_I)$.

\itemitem{--}
$\Kc^\wedge_{f_I}$, the sheaf of functions on the punctured formal neighborhood
of $\Gamma(f_I)$.

\itemitem{--}
$\Kc_{f_I}^\sqr\subset\Kc^\wedge_{f_I}$, the subsheaf of sections
whose restriction to $S_\red\times C$ lies in
$\Oc^\wedge_{f_{I,\red}}\subset\Kc^\wedge_{f_{I,\red}}.$

\proclaim{(3.3.1) Proposition}  The sheaves
$\Oc^\wedge_{f_I}$ and $\Kc_{f_I}^\sqr$ are sheaves of local rings.
\endproclaim

\noindent {\sl Proof:}
The case of $\Oc^\wedge_{f_I}$ is immediate because
$(\Gamma(f_I),\Oc^\wedge_{f_I})$ is the formal neighborhood of $\Gamma(f_I)$ in
$S\times C$, hence it is a locally ringed space.
The case of $\Kc_{f_I}^\sqr$ is analog to the proof of Proposition 2.2.1.
Indeed we can assume that $C=\Spec k[t]$ and $S=\Spec R$.
Set $b_i=f_i^*(t)$, $K_{f_I}^\sqr=H^0(\Gamma(f_I),\Kc_{f_I}^\sqr)$, and
$O_{f_I}^\wedge=H^0(\Gamma(f_I),\Oc_{f_I}^\wedge)$.
The ring $K_{f_I}^\sqr$ is identified with the series
$\sum_{l>>-\infty}a_l(t)\prod_i(t-b_i)^l$ where $a_l(t)\in R[t]$ are polynomials
of degree less than $|I|$ with nilpotent coefficients if $l<0$,
and the subring $O_{f_I}^\wedge$ with the series such that $a_l(t)=0$ if $l<0$,
see Example (3.7).
We have $\underline\Spec(K_{f_I}^\sqr)=\underline\Spec(O_{f_I}^\wedge)$ because
$K_{f_I}^\sqr$ is a nilpotent extension of $O_{f_I}^\wedge$.
Let $i\,:\,\Gamma(f_I)\to\underline\Spec(O_{f_I}^\wedge)$ be the natural map.
Then
$$i^{-1}\Oc_{\Spec(K_{f_I}^\sqr)}=\Kc_{f_I}^\sqr.$$
The proposition follows.
\qed

\vskip3mm

\proclaim {(3.3.2) Definition}
Let $I$ be a nonempty finite set.

(a) We define the contravariant functors $\l^0_{X,C^I}$,
$\l_{X,C^I}$ from $\Schb$ to $\bold{Sets}$ as follows.
For a scheme $S$ the set $\l^0_{X,C^I}(S)$ consists of pairs
$(f_I,\rho)$ such that
$$f_I\in\Hom_{\Schb_k}(S,C^I)\and
\rho\in\Hom_{\Lrsb_k}\bigl((\Gamma(f_I),\Oc^\wedge_{f_I}),\,X\bigr).$$
The contravariant functor $\l_{X,C^I}$ is defined
similarly but with $\Oc^\wedge_{f_I}$ replaced by $\Kc_{f_I}^\sqr$.

(b) If $X$ is affine, the functor $\tilde\l_{X,C^I}$
from $\Algb$ to $\bold{Sets}$ is such that
the set $\tilde\l_{X,C^I}(R)$ consists of pairs
$(f_I,\rho)$ with
$$f_I\in\Hom_{\Schb_k}(\Spec R,C^I)\and
\rho\in\Hom_{\Algb_k}\bigl(k[X],\H^0(\Gamma(f_I),\Kc^\wedge_{f_I})\bigr).$$
\endproclaim

\noindent
The embeddings of sheaves of rings
$\Oc^\wedge_{f_I}\hookrightarrow\Kc_{f_I}^\sqr\hookrightarrow\Kc_{f_I}^\wedge$
induce embeddings of functors
$\l^0_{X,C^I}\hookrightarrow\l_{X,C^I}\hookrightarrow\tilde\l_{X,C^I}.$

\proclaim {(3.3.3) Proposition} Let $X$ be a fixed scheme
of finite type. The families of functors $(\l^0_{X,C^I})$,
$(\tilde\l_{X,C^I})$, $(\l_{X,C^I})$ with $I$ running over
nonempty finite sets, form each a factorization monoid in the category
of functors.
\endproclaim

\noindent {\sl Proof:} This is almost obvious by construction. Indeed,
let $J\twoheadrightarrow I$ be a surjection of nonempty sets. Then $\Delta^{(J/I)*}$
of the $J$th functor in any of the three families takes $S$ into the set of
$(f_J, \rho)$ where $f_J$ is a
morphism $S\to C^J$ which in fact lie in the image of
$\Delta^{(J/I)}$. Thus $f_J$ comes from a uniquely defined
$I$-tuple $f_I$. Now, $\Gamma(f_J) = \Gamma(f_I)$ and so each
of the three sheaves $\Oc^\wedge, \Kc^\wedge, \Kc^\sqr$ associated
to them coincide. This gives the datum (a) of Definition 3.2.1.
Similarly,  $j^{(J/I)^*}$ applied to the $J$th functor in any of the families,
takes $S$ into the set of $(f_J, \rho)$ where $f_J: S\to U^{(J/I)}$.
But this means that $\Gamma(f_J) = \coprod_{i\in I} \Gamma(f_{J_i})$
and hence on the level of set of morphisms of $\Gamma(f_J)$ equipped
with any of the three sheaves of rings, we get a direct product.
This gives the datum (b), i.e. the isomorphism $\kappa^{(J/I)}$.
The associativity of these isomorphisms is obvious. \qed

\vskip3mm

Therefore, to establish Theorem 3.2.4,
we need only to prove the representability
of the functors $\l^0_{X,C^I},$ $\l_{X,C^I}$
and, when $X$ is affine, of $\tilde\l_{X,C^I}$, as well as to prove the formal
smoothness of the structure morphisms  of the representing objects
to $C^I$.

\vskip3mm

\noindent {\bf (3.4) The global space of germs of arcs.}

\proclaim{(3.4.1) Proposition}
The functor $\l^0_{X,C^I}$ is representable by a scheme $\Lc^0(X)_{C^I}$
of infinite type over $C^I$.
\endproclaim

\noindent{\sl Proof:} For $n>0$ let $C^{(n)}$ be the $n$th symmetric product
of $C$. As $C$ is a smooth curve, $C^{(n)}$ is identified with
the Hilbert scheme $\Hilb^n(C)$ parametrizing
subschemes in $C$ of finite length $n$. Explicitly, to a point
of $C^{(n)}$, i.e. an effective divisor $D$ on $C$ of degree $n$,
there corresponds the subscheme $Z_D=\Spec\bigl(\Oc_C/\Oc_C(-D)\bigr)$.
The following lemma is well-known.

\proclaim{(3.4.2) Lemma}
Let $T$ be a $k$-scheme, and $X,Z$ be any $T$-schemes.
Assume that the morphism $Z\to T$ is finite,
and that $X\to T$ is of finite type.
The contravariant functor
$$\Schb\to\Setsb,\quad S\mapsto\Hom_{_T}(S\times Z,X),$$
is represented by a $T$-scheme $\underline{Hom}_{_T}(Z,X)$.
\endproclaim

\noindent
Let $Z\subset C\times C^{(n)}$ be the total space of the family
of the schemes $Z_D$, $D\in C^{(n)}$.
We set $Map^n(C,X)=\underline{Hom}_{C^{(n)}}(Z,X\times C^{(n)}).$
Let $u_I^n$ be the composition of the maps
$$C^I\to C^{I\times[1,n+1]}\to C^{((n+1)|I|)} $$
where the first map takes the $I$-uple $(c_i, i\in I)$ to
the $I\times[1,n+1]$-uple $(c_i,c_i,...,c_i, i\in I)$,
each $c_i$ counted $n+1$ times, and the second map is the projection
from the Cartesian product to the symmetric product.
Let $\Lc^0_n(X)_{C^I}$ be the fiber product
$$\Lc^0_n(X)_{C^I}=C^I\times_{C^{((n+1)|I|)}}Map^{(n+1)|I|}(C,X).
\leqno (3.4.3)$$
For any $I$-uple $c\in C^I$ we have inclusions of subschemes of $C$
$$u_I^0(c)\subset u_I^1(c)\subset u_I^2(c)\subset\cdots
\leqno (3.4.4)$$
Thus we get a projective system of schemes
$$\Lc^0_1(X)_{C^I}\leftarrow\Lc^0_2(X)_{C^I}\leftarrow\Lc^0_3(X)_{C^I}
\leftarrow\cdots$$
The morphisms in this projective system are affine because the embeddings
$u_I^n(c)\subset u_I^{n+1}(c)$ are purely nilpotent. Therefore we have
the limit scheme
$$\Lc^0(X)_{C^I}=\pro_n\Lc^0_n(X)_{C^I}.$$
We claim that the scheme $\Lc^0(X)_{C^I}$ represents the functor
$\l^0_{X,C^I}$.
Indeed a morphism from a $k$-scheme $S$ to $\Lc^0_n(X)_{C^I}$ is a pair
$(f_I,\rho)$ where $f_I\,:\,S\to C^I$ and
$\rho\,:\,u_I^n(f_I)\to X$ are morphisms of schemes.
Here $u_I^n(f_I)\subset S\times C$
is the subscheme of relative length $n|I|$ over $S$ corresponding to the
$S$-point $f_I$ of $C^I$ via $u_I^n.$
When we pass to the limit we get
$\pro_n\Oc_{u^n_I(f_I)}=\Oc^\wedge_{f_I}.$
\qed

\vskip3mm

\noindent{\bf (3.4.5) Remark.}
When $X=\Spec A$ is affine, the scheme
$\Lc^0(X)_C$, being the scheme of infinite jets of maps $C\to X$,
is the spectrum of a commutative $\Cal O_C$-algebra
with connection along $C$ (in fact, of the universal such algebra
with connection generated by $\Cal O_C\otimes_k A$). According
to Beilinson and Drinfeld \cite{BD1}, \cite{G}, commutative $\Cal O_C$-algebras
with connection are particular case of chiral algebras which, in their
turn, give factorization algebras. So our construction in this case is
a particular case of theirs.

\proclaim {(3.4.6) Proposition}
Let $\phi: X\to Y$ be an \'etale morphism of schemes of finite type,
and $\pi: C\to D$ be a morphism of smooth curves. 
Then 
 
(a) if $\pi$ is \'etale then each morphism
$\Lc^0_n(\phi)_{\pi^I}\,:\,
\Lc^0_n(X)_{C^I}\to\Lc^0_n(Y)_{D^I}$ is \'etale,

(b) the square
$$\matrix
\Lc^0_n(X)_{C^I}&\to&
\Lc^0_n(Y)_{D^I}\cr
\downarrow&&\downarrow\cr
\Lc^0_0(X)_{C^I}&\to&\Lc^0_0(Y)_{D^I},
\endmatrix $$
as well as the analogous square for 
$\Lc^0(X)_{D^I}$, $\Lc^0(Y)_{D^I}$, is Cartesian. 
\endproclaim

\noindent {\sl Proof:} 
(a)
Let $S$ be any scheme.
Given two compatible morphisms 
$$\alpha=(\alpha_I,\rho_\alpha)\,:\,S_\red\to\Lc^0_n(X)_{C^I},\quad
\beta=(\beta_I,\rho_\beta)\,:\,S\to\Lc^0_n(Y)_{D^I},$$
we must prove that there is a unique morphism 
$\gamma=(\gamma_I,\rho_\gamma)\,:\,S\to\Lc^0_n(X)_{C^I}$ 
which splits the square 
$$\matrix
S_\red&\to&\Lc^0_n(X)_{C^I}\cr
\downarrow&&\downarrow\cr
S&\to&\Lc^0_n(Y)_{D^I},
\endmatrix$$
into two commutative triangles.
We have a commutative square
$$\matrix
S_\red&{\buildrel\alpha_I\over\lra}&C^I\cr
\downarrow&&\downarrow{\ss\pi^I}\cr
S&{\buildrel\beta_I\over\lra}&D^I.
\endmatrix\leqno(3.4.7)$$
Thus, $\pi$ being \'etale, there is a unique morphism 
$\gamma_I\,:\,S\to C^I$ splitting (3.4.7) into two commutative triangles.
Let the subschemes 
$$u_I^n(\alpha_I)\subset S_\red\times C,\quad
u_I^n(\beta_I)\subset S\times D,\quad
u_I^n(\gamma_I)\subset S\times C$$
be as in (3.4).
We have a Cartesian square 
$$\matrix
u_I^n(\alpha_I)&\to&S_\red\times C\cr
\downarrow&&\downarrow\cr
u_I^n(\gamma_I)&\to&S\times C,
\endmatrix$$
yielding a nilpotent extension of schemes
$u_I^n(\alpha_I)\to u_I^n(\gamma_I).$
This map fits into a commutative diagram
$$\matrix
u_I^n(\alpha_I)&{\buildrel\rho_\alpha\over\lra}&X\cr
\downarrow&&\downarrow{\ss\phi}\cr
u_I^n(\gamma_I)&\lra&Y\cr
\endmatrix\leqno(3.4.8)$$
where the lower arrow is the composition of the chain of maps
$$u_I^n(\gamma_I){\buildrel\Id\times\pi\over\lra}
u_I^n(\beta_I){\buildrel\rho_\beta\over\lra}Y.$$
Thus, $\phi$ being \'etale, there is a unique morphism 
$\rho_\gamma\,:\,u_I^n(\gamma_I)\to X$
splitting (3.4.8) into two commutative triangles.
We have proved (a). 

(b)
Let
$$\psi\,:\,\Lc^0_n(X)_{C^I}\to\Lc^0_0(X)_{C^I}\times_{\Lc^0_0(Y)_{D^I}}
\Lc^0_n(Y)_{D^I}$$
be the morphism induced by the diagram.
To prove that $\psi$ is an isomorphism we show that for every scheme $S$
and every two compatible morphisms
$$\alpha\,:\,S\to\Lc^0_0(X)_{C^I},\quad
\beta\,:\,S\to\Lc^0_n(Y)_{D^I}$$
there is a unique $\gamma\,:\,S\to\Lc^0_n(X)_{C^I}$ such that
$\psi(\gamma)=(\beta,\alpha)$.
By definition $\alpha=(\alpha_I,\rho_\alpha)$ with
$\alpha_I\,:\,S\to D^I$, $\rho_\alpha\,:\,u_I^n(\alpha_I)\to Y$ morphisms 
of schemes.
Similarly $\beta=(\beta_I,\rho_\beta)$ with
$\beta_I\,:\,S\to C^I$, $\rho_\beta\,:\,\Gamma(\beta_I)\to X$.
We look for
$\gamma=(\gamma_I,\rho_\gamma)$ with
$\gamma_I\,:\,S\to C^I$, $\rho_\gamma\,:\,u_I^n(\gamma_I)\to X$.

We first prove the existence of $\gamma$.
Take $\gamma_I=\beta_I$.
Next, by compatibility of $\alpha$ and $\beta$ we have
$(\Id\times\pi)(\Gamma(\beta_I))=\Gamma(\alpha_I)$, and this implies
that $\Id\times\pi$ induces a morphism 
$\varpi\,:\,u^n_I(\beta_I)\to u^n_I(\alpha_I)$.
Composing $\varpi$ with $\rho_\alpha$ we get a diagram
$$\matrix
\Gamma(\beta_I)&\to&X\cr
\downarrow&&\downarrow\cr
u^n_I(\beta_I)&\to&Y,
\endmatrix$$
with horizontal arrows $\rho_\beta$, $\rho_\alpha\circ\varpi$ 
and right vertical arrow $\pi$.
This diagram is commutative by compatibility of $\alpha$ and $\beta$.
Now the left vertical arrow is a nilpotent embedding, 
while the right vertical arrow is \'etale.
Therefore there is a unique morphism $\rho_\gamma\,:\,u_I^n(\beta_I)\to X$
splitting the diagram intotwo commutative triangles.

We now check the unicity of $\gamma$.
That $\psi(\gamma)=(\beta,\alpha)$ means that, first, 
$\beta_I=\gamma_I$ and the diagram
(whose horizontal maps are $\rho_\beta$, $\rho_\gamma$ respectively)
$$\matrix
\Gamma(\beta_I)&\to&X\cr
\downarrow&\nearrow&\cr
u^n_I(\beta_I)&&
\endmatrix$$
commutes and, second,
$\alpha_I=\pi^I\circ\gamma_I$ and the diagram
(whose horizontal maps are $\phi\circ\rho_\gamma$, $\rho_\alpha$ respectively)
$$\matrix
u_I^n(\gamma_I)&\to&Y\cr
\downarrow&\nearrow&\cr
u^n_I(\alpha_I)&&
\endmatrix$$
commutes.
The second condition determines $\phi\circ\rho_\gamma$ uniquely as
$\rho_\alpha\circ\varpi$.
So $\rho_\gamma$ splits the following square into two commutative triangles
$$\matrix
\Gamma(\beta_I)&\to&X\cr
\downarrow&&\downarrow\cr
u_I^n(\beta_I)&\to&Y,
\endmatrix$$
so it is unique because $\phi$ is \'etale.
\qed

\vskip .3cm

\noindent {\bf (3.5)
The meromorphic loop space of an affine scheme.} Here we prove the
ind-representability of the functor $\tilde\l_{X,C^I}$ for
$X$ affine.
We first treat the case when $X=\AA^1$. Let $I$ be fixed. To every
point $(c_i)_{i\in I}\in C^I$ we associate the effective divisor
$\sum c_i$ on $C$. For every $m, n\geq 0$ let $\Cal A_{mn}$ be the vector
bundle on $C^I$ whose fiber over $(c_i)$ is the vector space
of global sections of the coherent 0-dimensional sheaf
$$\Oc_C\bigl(m \Sum c_i\bigr)\bigl/\Oc_C\bigl(-n \Sum c_i\bigr)$$
on $C$. Thus $\op{rk}(\Ac_{mn}) = |I|(m+n)$.
Let $A_{mn}$ be the total space of the bundle $\Ac_{mn}$,
considered as an algebraic variety over $C^I$.
When $m,n$ vary, these varieties form
a double inductive-projective system and the following is then
obvious.

\proclaim {(3.5.1) Proposition} (a) For every $m>0$ the limit
$\lim\limits_{\buildrel\longleftarrow\over n} A_{mn}$ exists as a scheme.
The ind-scheme $\tilde\Lc(\AA^1)_{C^I}=\ind_m\pro_n A_{mn}$
represents the functor $\tilde\l_{\AA^1,C^I}$.

(b) $\tilde\Lc(\AA^1)_{C^I}$ has a natural structure of a $k$-algebra
object in the category of ind-schemes over $C^I$.
\endproclaim

\noindent{\sl Proof:} A morphism of $k$-schemes
$\Spec R\to A_{mn}$ is a pair
of a morphism of $k$-schemes $f_I\,:\,S\to C^I$, and
a morphism of sheaves of $k$-algebras
$\Oc_{A_{mn}}\otimes_{\Oc_{C^I}}R\to R$
(here, the right hand side is identified with the constant sheaf on $C^I$).
There is a canonical isomorphism
$$\ind_m\pro_n\Hom\bigl(\Oc_{A_{mn}}\otimes_{\Oc_{C^I}}R,R\bigr)
\simeq\H^0\bigl(\Gamma(f_I),\Kc_{f_I}\bigr)
\simeq\Hom_{\Algb}\bigl(k[x],\H^0(\Gamma(f_I),\Kc_{f_I})\bigr).$$
Claim $(a)$ is proved.
In part $(b)$, the variety $A_{mn}$ is not a ring scheme.
However there is an obvious map $A_{mn}\times A_{m'n'}\to A_{m+m',n+n'}$
which induces a ring structure on the limit of the ind-pro system.
\qed

\proclaim {(3.5.2) Proposition} For any affine $X$ of finite type
the functor $\tilde\l_{X,C^I}$ is representable by
an ind-scheme $\tilde\Lc(X)_{C^I}$ over $C^I$.
\endproclaim

\noindent {\sl Proof:}
By part (b) of Proposition 3.5.1,
for every $f\in k[x_1, ..., x_d]$ we have a morphism of ind-schemes
$$\tilde\Lc(f)_{C^I}: \bigl(\tilde\Lc(\AA^1)_{C^I}\bigr)^d\to
\tilde\Lc(\AA^1)_{C^I}.$$
If now the scheme $X$ is given in some $\AA^d$ be equations
$f_j(x_1, ..., x_d)=0$, then $\tilde\l_{X,C^I}$ is
represented by the ind-scheme $\tilde\Lc(X)_{C^I}$ which
is the intersection, in $\bigl(\tilde\Lc(\AA^1)_{C^I}\bigr)^d$,
of the preimages of 0 under the $\tilde\Lc(f_j)_{C^I}$, i.e.
the inverse limit of an obvious diagram in the category
of ind-schemes. \qed

\vskip .1cm

\proclaim{(3.5.3) Corollary} For an affine scheme $X$ the functor
$\l_{X,C^I}$ is represented by an ind-scheme $\Lc(X)_{C^I}$
which is the inductive limit of the formal neighborhoods
of $\Lc^0(X)_{C^I}$ in the schemes of an inductive system
for $\tilde\Lc(X)_{C^I}$.
\endproclaim

\noindent {\bf (3.6) The global loop space of an arbitrary scheme.}
Let now $X$ be an arbitrary scheme of finite type. The (ind-)representability
of the functor $\l_{X,C^I}$ follows from Corollary 3.5.3 for
the affine case and from the general gluing properties of the functors
summarized in the next proposition.

\proclaim {(3.6.1) Proposition}
(a) The functor $\l_{X,C^I}$ is a sheaf on $\Schb$.

(b)  If $U\subset X$ is an open subset, then the induced
morphism of functors $\l_{U,C^I}\to\l_{X,C^I}$ is open.

(c) Let $\{U_\alpha\}_{\alpha\in A}$ be an open covering of $X$. Then
$\l_{X,C^I}$ is equal to the cokernel, in the category $\Shfb$,
of the pair of morphisms
$$\prod_{\a,\b}\l_{U_\a\cap U_\b, C^I} \rra \prod_\a
\l_{U_\a, C^I}.$$
\endproclaim

\noindent {\sl Proof:} (a) For any $f_I: S\to C^I$ the graph
$\G(f_I)$ is identified with $S$, so $\Cal K^\sqr_{f_I}$ can be 
regarded as a sheaf of local rings on $S$. So our statement
follows from the fact that the representable functor $\eta_X$
on $\Lrsb$ is a sheaf.

\vskip .1cm

(b) Let $S$ be a scheme and $u$ a morphism 
of functors $\eta_S\to\l_{X, C^I}$. WE need to prove that the
 fiber product of $\eta_S$ and $\l_{U, C^I}$ over $\l_{X, C^I}$
is represented by a scheme $S'$ whose natural morphism to $S$ is
an open embedding. To see this, we view $u$ as an element of
$\l_{X, C^I}(S)$, so $u=(f_I, \rho)$ with $f_I: S\to C^I$
and $\rho: (\G(f_I), \Cal K^\sqr_{f_I})\to X$. Notice that
$\G(f_I)\simeq S$, so $\rho$ gives, in particular, a continuous  map of
topological spaces $\bar\rho: S\to X$. It is clear then that the
fiber product mentioned above is represented by the open subset
$S'=\bar\rho^{-1} (U)
\subset S$.

\vskip .1cm

(c) This follows from (b) and from Lemma 2.4.7(c). 
\qed

\vskip .2cm

An immediate corollary of Proposition 3.6.1 and
Corollary 3.5.3 is that the functor
$\l_{X,C^I}$ is representable by an ind-scheme $\Lc(X)_{C^I}$ over $C^I$.
Then, Proposition 3.3.3 implies that
the collections of schemes and ind-schemes
$\Lc^0(X)_{C^I}$ and $\Lc(X)_{C^I}$ are factorization monoids
in the categories of schemes and ind-schemes.
To finish the proof of Theorem 3.2.4 it remains to establish
part (b) of the following.

\proclaim{(3.6.2) Proposition}
(a) If $\phi: X\to Y, \pi: C\to D$ are \'etale morphisms, then the induced
morphism
$\Lc(\phi)_{\pi^I}: \Lc(X)_{C^I}\to \Lc(Y)_{D^I}$
is formally \'etale.

(b) If $X$ is smooth then the morphism
$\Lc(X)_{C^I}\to C^I$ is formally smooth.
\endproclaim

\noindent{\sl Proof:}
To prove Claim $(b)$ it is sufficient to observe that if
$U\subset X$, $V\subset C$ are affine open sets with \'etale maps
$U\to\AA^d$, $V\to\AA^1$ then the composition of maps
$$\Lc(U)_{V^I}\to\Lc(\AA^d)_{\AA^I}\to\AA^I$$
is formally smooth by Claim $(a)$ and  the statement (3.7.5)
of the example below. 

The proof of Claim $(a)$ is similar to that of Proposition 2.5.4.
Let $S$ be a scheme.
We are given $\alpha\,:\,S_\red\to\Lc(X)_{C^I}$, $\beta\,:\,S\to\Lc(Y)_{D^I}$,
and we look for a unique $\gamma$ which splits the square
$$\matrix
S_\red&\to&\Lc(X)_{C^I}\cr
\downarrow&&\downarrow\cr
S&\to&\Lc(Y)_{D^I}.
\endmatrix
$$
The morphism $\alpha$ consists of a pair 
$(\alpha_I\,:\,S_\red\to C^I,\,
\rho_\alpha\,:\,(\Gamma(\alpha_I),\Kc^\sqr_{\alpha_I})\to X)$.
Similarly $\beta$ consists of a pair 
$(\beta_I\,:\,S\to D^I,\,
\rho_\beta\,:\,(\Gamma(\beta_I),\Kc^\sqr_{\beta_I})\to Y)$.
We must construct a pair 
$(\gamma_I\,:\,S\to C^I,\,
\rho_\gamma\,:\,(\Gamma(\gamma_I),\Kc^\sqr_{\gamma_I})\to X)$.
There is a map $\gamma_I$ splitting the square
$$\matrix
S_\red&\to&C^I\cr
\downarrow&&\downarrow\cr
S&\to&D^I
\endmatrix
$$
into two commutative triangles, because $\pi^I$ is \'etale.
We have a Cartesian square
$$\matrix
\Gamma(\alpha_I)&\to&S_\red\times C\cr
\downarrow&&\downarrow\cr
\Gamma(\gamma_I)&\to&S\times C,
\endmatrix$$
which implies that we have a morphism of ringed spaces
$i\,:\,(\Gamma(\alpha_I),\Kc^\sqr_{\alpha_I})
\to(\Gamma(\gamma_I),\Kc^\sqr_{\gamma_I})$
with nilpotent kernel.
This map fits into a diagram of ringed spaces
$$\matrix
(\Gamma(\alpha_I),\Kc^\sqr_{\alpha_I})
&\to&(\Gamma(\gamma_I),\Kc^\sqr_{\gamma_I})\cr
\downarrow&&\downarrow\cr
X&\to&Y,
\endmatrix$$
where 
the left vertical arrow is 
$\rho_\alpha\,:\,(\Gamma(\alpha_I),\Kc^\sqr_{\alpha_I})\to X,$
the right vertical arrow is the composition of
$j\,:\,(\Gamma(\gamma_I),\Kc^\sqr_{\gamma_I})
\to(\Gamma(\beta_I),\Kc^\sqr_{\beta_I})$ and
$\rho_\beta\,:\,(\Gamma(\beta_I),\Kc^\sqr_{\beta_I})\to Y,$
the topological map underlying $j$ is
$j_\top=(\Id_S\times \pi^I)\,:\,\Gamma(\gamma_I)\to\Gamma(\beta_I)$
and the structure morphism 
$j^*_\top\,:\,\Kc^\sqr_{\beta_I}\to\Kc^\sqr_{\gamma_I}$ is an isomorphism.
Note that $i$ yields an isomorphism of underlying topological spaces.
Let $p$ be any point in $\Gamma(\alpha_I)$.
We have a diagram of stalks
$$\matrix
(\Kc^\sqr_{\alpha_I})_p&\leftarrow&(\Kc^\sqr_{\gamma_I})_p\cr
\uparrow&&\uparrow\cr
\Oc_{X,\rho_\alpha(p)}&\leftarrow&
\Oc_{Y,\phi\rho_\alpha(p)}
\endmatrix$$
where the upper horizontal map is $i^*$, hence has a nilpotent kernel, and 
the lower horizontal map is $\phi^*$, hence is \'etale. 
Therefore there is a unique morphism of rings
$\Oc_{X,\rho_\alpha(p)}\to(\Kc^\sqr_{\gamma_I})_p$ for each $p$.
Theses morphisms give a morphism of sheaves of rings
$\rho_{\gamma,\top}^{-1}\Oc_X\to\Kc_{\gamma_I}^\sqr$,
i.e. a desired morphism of ringed spaces $\gamma$.
\qed

\vskip .3cm

\noindent {\bf (3.7) Example : the cases $C=\AA^1$ and $X=\AA^d$.}
Let $C=\Spec k[t]$. If $S=\Spec R$,
a morphism $f_I = (f_i): S\to C^I$ is the same as a collection of elements
$b_i = f_i^*(t)\in R$. Assume these have been fixed.
Then the subscheme $u_I^n(f_I)\subset S\times C$
from the proof of Proposition 3.4.1 is described explicitly :
$$u_I^n(f_I)=\Spec\bigl(R[t]\bigl/\Prod_{i\in I}(t-b_i)^{n+1}\bigr).\leqno
(3.7.1)$$
This implies the following.

\proclaim {(3.7.2) Proposition}
(a) The ring
$$H^0(\Gamma(f_I),\Oc^\wedge_{f_I})=
\pro_n\bigl(R[t]\bigl/\Prod_i(t-b_i)^{n+1}\bigr)$$
is identified with the set of series
$\sum_{l=0}^\infty a_l(t) \prod_i(t-b_i)^l$,
where $a_l(t)\in R[t]$ are polynomials
of degree less than $|I|$ (such polynomials form a set of representatives for
$R[t]/\prod_i(t-b_i)$).

(b) The ring $H^0(\Gamma(f_I), \Kc_{f_I})$ is identified with the set of series
$\sum_{l\gg -\infty}a_l(t)\prod_i(t-b_i)^l$, where $a_l(t)$ are as in (a).

(c) The subring $H^0(\Gamma(f_I), \Kc_{f_I}^\sqr)$ is identified with
the series as in (b) but with the condition that all the coefficients
of the polynomials $a_l(t), l<0$, are nilpotent
elements of $R$.
\endproclaim

\noindent {\bf (3.7.3) Notation.} For fixed $b_i\in R$ and $a\in
H^0(\Gamma(f_I), \Kc_{f_I})$ we denote by $a_l(t)\in R[t]$
the $l$th coefficient of the series corresponding to  $a$ by Proposition
3.7.2(b) and by $a_{l\nu}\in R$, $\nu=0, ..., |I|-1$, the
$\nu$th coefficient of the polynomial $a_l(t)$.

\vskip .2cm

\noindent
Assume moreover that $X=\AA^d$ with coordinates $x_1, ..., x_d$.
In this case we can give a completely explicit description of the
ind-schemes $\tilde\Lc(\AA^d)_{\AA^I}$ and $\Lc(\AA^d)_{\AA^I}$,
using Proposition 3.7.2. Indeed,
$$\tilde\lambda_{\AA^d, \AA^I}(R)=
\bigr\{(b_i,\rho)\,\big|\,(b_i)\in R^I,\,
\rho: k[x]\to H^0(\Gamma(f_I),\Kc_{f_I})\bigr\}.$$
The algebra homomorphism $\rho$ is uniquely determined by the choice of
$$\rho(x_j)=\sum_{l\gg-\infty}\sum_{\nu=0}^{|I|-1} a^{(j)}_{l\nu}
t^\nu \prod_{i\in I}(t-b_i)^l, \quad a^{(j)}_{l\nu}\in R.$$
A choice of $\rho(x_j)$ is the same as a choice
of elements $a^{(j)}_{l\nu}\in R$. Thus the universal case corresponds
to $b_i, a^{(j)}_{l\nu}$, with $i\in I$, $j\in [1,d]$,
$l\in\ZZ$ and $\nu\in[0,|I|-1]$, being independent variables, i.e.
$$\tilde\Lc(\AA^d)_{\AA^I}=\ind_N\Spec\,k\bigl[b_i,a^{(j)}_{l\nu};l\geq-N\bigr].
\leqno (3.7.4)$$
and
$$\Lc(\AA^d)_{\AA^I}=\ind_N \op{Spf}\, k[b_i, a^{(j)}_{l\nu}; l\geq 0]
[[a^{(j)}_{l\nu}; -N\leq l\leq -1]].\leqno (3.7.5)$$
Notice further that
$$\Lc^0_n(\AA^d)_{\AA^I}=\Spec k[b_i, a^{(j)}_{l\nu}; 0\leq l\leq n].$$
In particular, we have a natural morphism
$$\theta_{\AA^d, I}: \Lc(\AA^d)_{\AA^I}\to \Lc^0_0(\AA^d)_{\AA^I}
\leqno (3.7.6)$$
taking all $a^{(j)}_{l\nu}, l\neq 0$, to 0.

\proclaim{(3.7.7) Proposition}
If $X$, $C$ are affine and $\phi\,:\,X\to\AA^d$, $\pi\,:\,C\to\AA^1$ 
are \'etale morphisms, there is a Cartesian diagram
$$\matrix
\Lc(X)_{C^I}&\to&\Lc(\AA^d)_{\AA^I}\cr
\downarrow&&\downarrow\cr
\Lc^0_0(X)_{C^I}&\to&
\Lc^0_0(\AA^d)_{\AA^I}.
\endmatrix$$
\endproclaim

\noindent{\sl Proof :}
The proof is similar to that of Proposition 2.6.1.(b).
Consider the diagram
$$\matrix 
\Lc^0(X)_{C^I}&\to&\Lc^0_0(X)_{C^I}\cr
\downarrow&&\downarrow\cr
\Lc(X)_{C^I}&\to&\Lc^0_0(\AA^d)_{\AA^I}
\endmatrix$$
where the right vertical arrow is $\Lc^0_0(\phi)_{\pi^I}$,
the lower horizontal arrow is the composition of
$\theta_{\AA^d,I}\,:\,\Lc(\AA^d)_{\AA^I}\to\Lc^0_0(\AA^d)_{\AA^I}$
defined in (3.7.6) and
$\Lc(\phi)_{\pi^I}\,:\,\Lc(X)_{C^I}\to\Lc(\AA^d)_{\AA^I}$.
The upper horizontal arrow is the natural projection of $\Lc^0(X)_{C^I}$ to
the 0-th term of the projective system, see (3.4).
Further, the left vertical arrow is an inductive limit of
nilpotent embeddings of schemes, while the right vertical 
arrow is \'etale by (3.4.6). Therefore there exists a unique morphism
$\theta_{\phi,\pi^I}\,:\,\Lc(X)_{C^I}\to\Lc^0_0(X)_{C^I}$ splitting the diagram into two commutative triangles.

Combining $\theta_{\phi,\pi^I}$ with 
$\Lc(\phi)_{\pi^I}\,:\,\Lc(X)_{C^I}\to\Lc(\AA^d)_{\AA^I}$ 
we obtain a morphism
$$\psi\,:\,\Lc(X)_{C^I}\to\Lc^0_0(X)_{C^I}\times_{\Lc_0^0(\AA^d)_{\AA^I}}
\Lc(\AA^d)_{\AA^I}.\leqno(3.7.8)$$
We claim that $\psi$ is an isomorphism 
and to prove this we construct its inverse $\chi$.
Let $S$ be a scheme.
A morphism from $S$ to the RHS of (3.7.8)
is a compatible pair
$$(\alpha\,:\,S\to\Lc(\AA^d)_{\AA^I},\,\beta\,:\,S\to\Lc_0^0(X)_{C^I}).$$
We construct a morphism $\chi(\a,\beta)\,:\,S\to\Lc(X)_{C^I}$.
For this notice that $\alpha(S_\red)\subset\Lc^0(\AA^d)_{\AA^I}$
and $\beta(S_\red)\subset\Lc^0_0(X)_{C^I}$.
By Proposition 3.4.6.(b) we have a map
$\gamma\,:\,S_\red\to\Lc^0(X)_{C^I}$.
The composition $\tilde\gamma$ of $\gamma$ and the embedding
$\Lc^0(X)_{C^I}\to\Lc(X)_{C^I}$ gives a commutative diagram
$$\matrix
S_\red&\to&\Lc(X)_{C^I}\cr
\downarrow&&\downarrow\cr
S&\to&\Lc(\AA^d)_{\AA^I}.
\endmatrix$$
Because $\Lc(\phi)_{\pi^I}$ is formally etale by (3.6.2)(a),
we set $\chi(\a,\b)$ to be the unique splitting of the diagram into 
two exact triangles.
The verification that $\chi$ is inverse to $\psi$ is straightforward.
\qed

\vskip .3cm

\noindent {\bf (3.8) The global loop space as an ind-object.}
We first consider the case $X=\AA^d, C=\AA^1$, employing the notations
of (3.7). Let $\eps = (\eps_{-1}, \eps_{-2}, ...)\in\EE$ be as in
(2.6). Define the scheme
$$\Lc^\eps(\AA^d)_{\AA^I}=\Spec\,\bigl( k[b_i, a^{(j)}_{l\nu}]
\bigl/ ( a^{(j)}_{l, \nu_1}\cdots a^{(j)}_{l,\nu_{1+\eps_l}}; l<0)
\bigl),\leqno (3.8.1)$$
where $j\in [1,d]$, $l\in\ZZ$, and
$\nu_1, ..., \nu_{1+\eps_l}\in [0,|I|-1]$ are arbitrary for $l<0$.
It is clear that
$$\Lc(\AA^d)_{\AA^I}=\ind_\eps \Lc^\eps(\AA^d)_{\AA^I}.$$
next, assume that $X$ is an arbitrary smooth scheme of finite type and
$C$ is an arbitrary smooth curve. Then $X$ can be covered by open $U=\Spec \, A$
possessing an \'etale map $\phi: U\to\AA^d$ and similarly $C$ can be
covered by open $V$ possessing an \'etale $\pi: V\to \AA^1$. 
We set
$$\Lc^\eps(\phi)_{\pi^I}=
\Lc^0_0(U)_{V^I}\times_{\Lc^0_0(\AA^d)_{\AA^I}}\Lc^\eps(\AA^d)_{\AA^I},$$
where the map $\Lc^\eps(\AA^d)_{\AA^I}\to\Lc^0_0(\AA^d)_{\AA^I}$ 
is the restriction of $\theta_{\AA^d,I}$ defined in (3.7.6).
This is an affine scheme.

\proclaim {(3.8.2) Proposition} 
The ind-object $``\ind_\eps" \Lc^\eps(\phi)_{\pi^I}$ in $\bold{Sch}$
is isomorphic to $\Lc(U)_{V^I}$.
\endproclaim

\vskip3mm

\noindent{\sl Proof :}
This follows from the case$X=\AA^d$, $C=\AA^1$, from Proposition 3.7.7, 
and from the fact that fiber products commute with filtering inductive limits.
\qed

\vskip .3cm

\noindent {\bf (3.9) The global loop space as an ind-pro-object.}
We keep the notation of (3.8). For $\eps\in\EE, n\geq 1$ consider the
scheme
$$\Lc^\eps_n(\AA^d)_{\AA^I}=\Spec\bigl(k[b_i,a^{(j)}_{l\nu}; l\leq n]
\bigl/(a^{(j)}_{l\nu_1}\cdots a^{(j)}_{l,\nu_{1+\eps_l}})\bigr).
\leqno (3.9.1)$$
We set
$$\Lc^\eps_n(\phi)_{\pi^I}=
\Lc^0_0(U)_{V^I}\times_{\Lc^0_0(\AA^d)_{\AA^I}}\Lc^\eps_n(\AA^d)_{\AA^I},
\leqno(3.9.2)$$
where the map $\Lc^\eps_n(\AA^d)_{\AA^I}\to\Lc^0_0(\AA^d)_{\AA^I}$ 
is defined as $\theta_{\AA^d,I}$ in (3.7.6).

\proclaim {(3.9.3) Proposition} 
(a) The scheme $\Lc^\eps_n(\phi)_{\pi^I}$ is of finite type 
and is a nilpotent extension of $\Lc^0_n(U)_{V^I}$.
The second projection
$\Lc^\eps_n(\phi)_{\pi^I}\to\Lc^\eps_n(\AA^d)_{\AA^I}$
is \'etale.
Moreover $\Lc^\eps(\phi)_{\pi^I}=\pro_n\Lc^\eps_n(\phi)_{\pi^I}$.

(b) 
The schemes $\Lc^\eps_n(\phi)_{\pi^I}$ form a double ind-pro-system with
Cartesian squares ($n'\geq n,\ \eps\leq\eps'$)
$$\matrix
\Lc^{\eps}_{n'}(\phi)_{\pi^I}&\hookrightarrow&\Lc^{\eps'}_{n'}(\phi)_{\pi^I}\cr
\downarrow&&\downarrow\cr
\Lc^{\eps}_{n}(\phi)_{\pi^I}&\hookrightarrow&\Lc^{\eps'}_{n}(\phi)_{\pi^I}
\endmatrix$$
where the vertical arrows are smooth affine morphisms.

(c) The ind-pro-object
$ ``\ind_\eps" ``\pro_n" \Lc^\eps_n(\phi)_{\pi^I}$
is independent, up to isomorphism, on $\phi, \pi$.
\endproclaim

\noindent{\sl Proof:}
The proof of $(a)$ is similar to that of Proposition 3.8.2.
In order to prove that
$\Lc^\eps_n(\phi)_{\pi^I}\to\Lc^\eps_n(\AA^d)_{\AA^I}$ is \'etale, 
due to (3.9.2) it is sufficient to check that the morphism
$\Lc^0_0(\phi)_{\pi^I}\,:\,\Lc^0_0(U)_{V^I}\to\Lc^0_0(\AA^d)_{\AA^I}$
is \'etale.
This is a consequence of Proposition (3.4.6)(a).

Claim $(b)$ is obvious in the case $\phi=\Id$, $\pi=\Id$. 
The general case follows from (3.9.2)
since the base change of a smooth affine morphism is still smooth affine,
and the base change of a Cartesian square is Cartesian.

(c) The pro-object
$``\pro_n"\Lc^\eps_n(\phi)_{\pi^I}$ in Pro($\Affb^\ft$) can be identified, 
due to Proposition 2.1.2, with the scheme $\Lc^\eps(\phi)_{\pi^I}$.
The ind-object $``\ind_\eps"\Lc^\eps(\phi)_{\pi^I}$ 
in Ind(Pro($\Affb^\ft$)) can then be identified with the ind-scheme
$\Lc(U)_{V^I}$.
\qed

\vskip3mm

\vfill\eject

\head 4. $\Dc$-modules  over ind-schemes\endhead

\vskip.3cm

\noindent {\bf (4.1) Reminder on $\Dc$-modules.} From now on we assume
that $\op{char}(k)=0$.
Fix a $k$-scheme $S$ of finite type.
Let $\Schb_S^\ft$ denote the category of $S$-schemes of finite type.
For any such scheme $X$
let $\Ob_X$ be the category of all quasi-coherent $\Oc_X$-modules.
For a morphism $f: X\to Y$ in $\Schb_S^\ft$ we denote
by  $f_*, f^*$  the functors of the direct
and inverse images on  $\Ob_X$, $\Ob_Y$.
If $f$ is a closed embedding and $\Ec\in\Ob_Y$, let
$$f^!\Ec=f^{-1}\Hc om_{\Oc_Y}(f_*\Oc_X,\Ec)\leqno (4.1.1) $$
be the inverse image of the subsheaf of $\Ec$
consisting of sections supported
scheme-theoretically on $f(X)\subset Y$.

For $X\in \Schb_S^\ft$ let $\Db_{X/S}$ 
be the category of coherent right $\Dc_{X/S}$-modules on $X$.
It is defined as follows,
see \cite{BD2, \S 7.10} or \cite{G, \S 0.2.2}. 
If $X$ is smooth over $S$, then we have the sheaf
of rings $\Dc_{X/S}$ of differential operators from $\Oc_X$ to itself
which are linear over $\Oc_S$.
An object of $\Db_{X/S}$ is then a coherent sheaf of right $\Dc_{X/S}$-modules.
It is quasicoherent over $\Oc_X$.
Next, if  $X$ admits a closed
embedding into a smooth $S$-scheme $Y$, one defines $\Db_{X/S}$
as $\Db_{Y/S,X}$, the full subcategory of $\Db_{Y/S}$ 
consisting of modules supported (as sheaves) on $X$. 
This definition is independent on the choice
of embedding: if $X$ is embedded into two smooth $S$-schemes $Y_1$ and $Y_2$,
then one has an equivalence $\Db_{Y_1/S, X}\to \Db_{Y_2/S, X}$
which is unique up to a unique isomorphism of functors.
Now, given any $X$, an embedding
into a smooth scheme always exists locally on $X$. Therefore we have
an open covering $X=\bigcup U_\a$, the categories $\Db_{U_\a/S}$,
$\Db_{U_{\a\b}/S}$ etc. and the obvious restriction functors among them.
One then  defines an object $\Mc$ of $\Db_{X/S}$
as a collection of objects $\Mc_\a\in \Db_{U_\a/S}$ together
with isomorphisms of their images in $\Db_{U_{\a\b}/S}$
whose images in $\op{Mor}(\Db_{U_{\a\b\gamma}/S})$ satisfy
the obvious compatibility conditions.

Given $X\in \Schb_S^\ft$ and $\Mc\in\Db_{X/S}$, we define a sheaf
$\Mc^{\Oc}\in\Ob_X$ as follows. If $X$ admits a closed embedding
$i: X\to Y$ with $Y$ smooth over $S$ and $\Mc$ is represented by a sheaf
of right $\Dc_{Y/S}$-modules supported on $X$ (which we also denote $\Mc$)
then we set $\Mc^{\Oc} =  i^! \Mc$. This definition is easily seen
to be independent on the choice of $Y$. In the case
of a general $X$ one defines $\Mc^\Oc$ by gluing the
sheaves given by above procedures on open parts of $X$ admitting
embeddings into smooth schemes.

Any smooth morphism $f\,:\,X\to Y$ in $\Schb_S^\ft$ induces the functor
of inverse image $f^\bullet\,:\,\Db_{Y/S}\to\Db_{X/S}$.
If $\omega_{X/Y}$ is the relative canonical bundle, then
$$f^\bullet(\Mc)^\Oc = f^*(\Mc^\Oc)\otimes_{\Oc_X}\omega_{X/Y}.$$
We have then  a canonical embedding
$$\Mc^\Oc\hookrightarrow f_*((f^\bullet\Mc)^\Oc\otimes\omega^{-1}_{X/Y}).
\leqno (4.1.2)$$

If $f\,:\,X\to Y$ is a closed embedding in $\Schb_S^\ft$,
we have an exact functor of direct image
$f_\bullet\,:\,\Db_{X/S}\to\Db_{Y/S}$.
In the particular case when $X,Y$ are smooth over $S$, 
we can view $\Mc\in\Db_{X/S}$ as a sheaf on $X$ and we have
$$f_\bullet(\Mc) = f_*(\Mc\otimes_{\Dc_{X/S}}\Dc_{X\to Y}).$$
Here $\Dc_{X\to Y}$ is the sheaf of differential operators from
$f^{-1}\Oc_Y$ to $\Oc_X$ linear over $\Oc_S$, see [BB], [G].

For a general closed embedding $f$ we have a canonical embedding
$$\Mc^\Oc\simeq f^!f_\bullet(\Mc)\hookrightarrow f^{-1}f_\bullet(\Mc).
\leqno(4.1.3)$$
For example, if $X$ is non-reduced and $i: X_\red\hookrightarrow X$
is the reduced part, then $i_\bullet$  identifies
$\Db_{X_\red/S}$ with $\Db_{X/S}$.

We have the following base change property.

\proclaim{(4.1.4) Lemma}
Suppose that in a Cartesian diagram of $S$-schemes of finite type
$$\matrix &X&\buildrel i\over\hookrightarrow&X'&\cr
&{\ss f}\bda&&\bda{\ss f'}&\cr
&Y&\buildrel j\over\hookrightarrow&Y'&\cr
\endmatrix$$
the morphisms $f, f'$ are smooth, $i,j$ are closed embeddings.
Then for any $\Mc\in\Db_{Y/S}$, we have
$$i_\bullet f^\bullet\Mc\simeq {f'}^\bullet j_\bullet\Mc.$$
\endproclaim

\vskip3mm

Note that if $f$ is a closed  embedding of smooth schemes over $S$,
the inverse image functor $f^\bullet$ is still defined.
Furthermore, the projection formula holds for right $\Dc$-modules.
More precisely, if $f$ is a closed or open embedding of smooth $S$-schemes
and $\Mc\in\Db_{X/S}$, $\Nc\in\Db_{Y/S}$, then there is a canonical isomorphism
$$f_\bullet\bigl(f^\bullet(\Nc)\otimes_{\Oc_X}\Mc\otimes_{\Oc_X}\o_X^{-1}\bigr)
\simeq\Nc\otimes_{\Oc_Y}f_\bullet(\Mc)\otimes_{\Oc_Y}\o_Y^{-1}.$$

\vskip .3cm

\noindent {\bf (4.2) $\Dc$-modules over pro-schemes.}

\proclaim{(4.2.1) Definition} (cf. \cite{Kap}.) Let $A$ be a filtering poset
and $(\Cb_\a)_{\a\in A}$ be an inductive system of categories
labelled by $A$. In other words, for  each $\a\leq\b$ we have a functor
$i_{\a\b}: \Cb_\a\to\Cb_\b$, for each $\a\leq\b\leq\g$ a natural isomorphism
$i_{\b\g}\circ i_{\a\b}\Rightarrow i_{\a\g}$ and these isomorphisms
satisfy the obvious coherence conditions for any $\a\leq\b\leq\g\leq\d$.

 The inductive limit $2\ind_\a\Cb_\a$ is the category whose
objects are pairs $(\a,x_\a)$, $\a\in A, x_\a\in\op{Ob}(\Cb_\a)$ and
$$\op{Hom}((\a,x_\a),(\b,y_\b))=\ind_{\g\geq \a, \b}
\op{Hom}_{\Cb_\a} (i_{\a\g}(x_\a), i_{\b\g}(y_\b)).$$

 \endproclaim

\proclaim{(4.2.2) Definition}
A $S$-scheme $X_\infty$ (possibly of infinite type) is called compact
if it can be represented as $\pro_\a\, X_\a$ where $(X_\a)_{\a\in A}$
is a filtering projective system over
$\bold{Sch}_S^{\op{ft}}$ such that all the maps $p_{\a\b}: X_\b\to X_\a$,
$\a\leq\b$, are affine morphisms.
\endproclaim

\proclaim {(4.2.3) Proposition} A scheme is compact if and only if it is
quasi-compact and quasi-separated.
The category of compact $k$-schemes can be identified with a full
subcategory in $\bold{Pro}(\bold{Sch}_S^{\op{ft}})$, via
$X_\infty=\pro_\a\,X_\a\mapsto``\pro_\a" X_\a$.
\endproclaim

\noindent {\sl Proof:} This follows from  [TT], Appendix C, Theorem C9.

\proclaim{(4.2.4) Definition}
(a)
A compact $S$-scheme $X_\infty$ is called smooth if it can be represented
as $\pro_\a\,X_\a$ for some $(X_\a)$ as in (4.2.2) with the
extra property that each $X_\a$ is a smooth $S$-scheme and each $p_{\a\b}$
is a smooth affine morphism.

(b)
$X_\infty$ is called almost smooth if it can be represented
as $\pro_\a\,X_\a$ for some $(X_\a)$ as in (4.2.2) with the
extra property that each $p_{\a\b}$ is a smooth affine morphism.
\endproclaim

\noindent{\bf (4.2.5) Example.}
For any smooth $X\in\Schb^{\ft}_k$,
the scheme $\Lc^0(X)$ is smooth and compact over $\Spec(k)$, and 
the scheme $\Lc^0(X)_{C^I}$ is smooth and compact over $C^I.$
Moreover,
the schemes $\Lc^\eps(\phi)$, $\Lc^\eps(\phi)_{\pi^I}$ from
(2.6.3) and (3.8) are almost smooth.

\vskip3mm

Let $q\,:\,X_\infty\to S$ be a compact almost smooth $S$-scheme 
and $(q_\a\,:\,X_\a\to S)$ be as in Definition 4.2.4.
We have two inductive systems of categories
$(\Ob_{X_\a}, p_{\a\b}^*)$, $(\Db_{X_\a/S}, p_{\a\b}^\bullet)$.
We set $\Db_{X_\infty/S}=2\ind_\a\Db_{X_\a/S}$.

\proclaim{(4.2.6) Proposition}
The category $\Db_{X_\infty/S}$ is independent, up to canonical equivalence
of categories, on the choice of $(X_\alpha)$ as in (4.2.4).
\endproclaim

\noindent{\sl Proof :}
We first consider the case when $X_\infty$ is smooth, so that each $X_\alpha$ 
is smooth over $S$.
Let $\Db^\ell_{X_\alpha/S}$ be the category of left coherent 
$\Dc_{X_\alpha/S}$-modules. We have an equivalence
$$\Db^\ell_{X_\alpha/S}\to\Db_{X_\alpha/S},\
\Nc\mapsto\Nc\otimes\omega_{X_\alpha/S}.$$
Let also $\Dc_{X_\infty/S}$ be the sheaf of rings of differential operators
on $X_\infty$ linear over $\Oc_S$.
It is equipped with the natural topology, see \cite{KT}, \S 1.7.
Using the pull-back of left $\Dc$-modules
(which is the same as for quasi-coherent sheaves)
we get an inductive system of categories $\Db^\ell_{X_\alpha/S}$.
It is proved in \cite{KT}, \S 1.9, that
$2\ind_\alpha\Db^\ell_{X_\alpha/S}$
is identified with the category of discrete, 
locally finitely generated quasicoherent sheaves of left 
$\Dc_{X_\infty/S}$-modules, 
and thus is independent on the choice of $(X_\alpha)$.
Therefore the category
$2\ind_\alpha\Db_{X_\alpha/S}$,
being equivalent to the previous one,
is also independent.

Now assume that $X_\infty$ is almost smooth.
Fix $\alpha_0\in A$.
Using a covering of $X_{\alpha_0}$ by affine open subsets, 
we reduce to the case when $X_{\alpha_0}$ (and thus $X_\infty$) is affine.
We can also assume that $\alpha_0$ is the minimal element in $A$.
Let us embed $X_{\alpha_0}$ as a closed subscheme 
into a smooth affine $S$-scheme $Y_{\alpha_0}$.
We can then extend each $p_{\alpha_0,\alpha}\,:\,X_\alpha\to X_{\alpha_0}$
to a smooth map $q_{\alpha_0,\alpha}\,:\,Y_\alpha\to Y_{\alpha_0}$.
We get then a smooth compact scheme $Y_\infty=\pro_\alpha Y_\alpha$
containing $X_\infty$ as a closed subscheme.
The category $2\ind_\alpha\Db_{X_\alpha/S}$ is then identified with the category
of sheaves of discrete locally 
finitely generated left $\Dc_{Y_\infty/S}$-modules supported on $X_\infty$.
\qed

\vskip3mm

Informally, an object of $\Db_{X_\infty/S}$ is a ``$\Dc$-module
pulled back from some $X_\a$".
Let $p_\a\,:\,X_\infty\to X_\a$ be the projection.
Let $p_\a^\bullet\,:\,\Db_{X_\a/S}\to\Db_{X_\infty/S},$
$\Mc_\a\mapsto(\a,\Mc_\a)$ be the canonical functor.
If the compact $S$-scheme $X_\infty$ is smooth, there is the functor
$$\Db_{X_\infty/S}\to 2\ind_\a\Ob_{X_\a}\subset\Ob_{X_\infty},\quad
p_\a^\bullet\Mc_\a\mapsto(p_\a^\bullet\Mc_\a)^\Oc:=
\bigl((p^\bullet_{\a\b}\Mc_\a)^\Oc
\otimes_{\Oc_{X_\b}}\o^{-1}_{X_\b}\bigr)_{\b\geq\a}.$$
In particular, we can associate to any $\Mc\in\Db_{X_\infty/S}$
its ``space of global sections", i.e., 
the direct image to $S$ as an $\Oc$-module.
This is a quasi-coherent sheaf on $S$ such that
if $\Mc=p^\bullet_\a\Mc_\a$ then
$$q_*(\Mc):=q_*(\Mc^\Oc)=
\ind_{\b\geq\a}(q_\b)_*\bigl(X_\b,(p^\bullet_{\a\b}
\Mc_\a)^\Oc\otimes_{\Oc_{X_\b}}\o^{-1}_{X_\b}\bigr).\leqno(4.2.7)$$
When $S=\Spec k$, we write $\Gamma(X_\infty,\Mc)$ for $q_*(\Mc)$.

\vskip .3cm

\noindent {\bf (4.3) $\Dc$-modules over ind-schemes.}

\proclaim{(4.3.1) Definition}
Let $A$ be a filtering poset and $(\Cb_\a)_{\a\in A}$ be a
projective system of categories labelled  by $A$.
In other words, for each $\a\leq\b$ we have
a functor $j_{\a\b}: \Cb_{\b}\to\Cb_\a$
and for any $\a\leq\b\leq\g$ a natural isomorphism
$j_{\a\b}\circ j_{\b\g}\Rightarrow j_{\a\g}$ satisfying the
obvious compatibility conditions.

 The projective
limit $2\pro \Cb_\a$ is the category whose objects
are systems consisting of objects $x_\a\in\Cb_\a$ given
for all $\a\in A$ and isomorphisms $j_{\a\b}(x_\b)\to x_\a$
given for each $\a\leq\b$ and satisfying the
compatibility condition for each $\a\leq\b\leq\g$.
Morphisms are defined in the obvious way.
\endproclaim

\proclaim{(4.3.2) Definition}
Let $X^\infty$ be an ind-$S$-scheme. We say that $X_\infty$ is discrete
over $S$
if it can be represented as $X^\infty = ``\ind_\a" X^\a$ where
$(X^\a)_{\a\in A}$ is a
filtering inductive system over $\bold{Sch}_S^{\op{ft}}$ such that each map
$i_{\a\b}\,:\,X^\a\to X^\b$, $\a\leq\b$, is a closed embedding.
\endproclaim

\noindent{\bf (4.3.3) Example.}
For any $\phi$, $U$ as in (2.7) the
ind-scheme $\Lc_n(\phi)$ is discrete over $\Spec(k)$.

\vskip .2cm

Let $q\,:\,X^\infty\to S$ be a discrete ind-scheme over $S$
and $(q^\a\,:\,X^\a\to S)$ be as in Definition 4.3.2. 
We have then the projective system of categories
$(\Ob_{X^\a}, i_{\a\b}^!)$. We define $\Ob_{X^\infty} = 2\pro \Ob_{X^\a}$.
If $(\Ec^\a, \g_{\a\b}: \Ec^\a\to i_{\a\b}^!\Ec^\b)$
is an object of $\Ob_{X^\infty}$, then
the direct images $q_*^\a(\Ec^\a)$ form an inductive system and
we define
$$q_*\Ec=\ind_\a\,q_*^\a(\Ec^\a).\leqno(4.3.4)$$
When $S=\Spec(k)$ we write $\Gamma(X^\infty,\Ec)$ for $q_*\Ec$.

We will also use the category $\hat\Ob_{X^\infty}$ which is the limit
of the projective system of categories
$(\Ob_{X^\a}, i_{\a\b}^*)$.
There is a functor
$$\Ob_{X^\infty}\times\hat\Ob_{X^\infty}\to\Ob_{X^\infty},\quad
(\Ec,\hat\Fc)\mapsto \Ec\otimes\hat\Fc=
(\Ec^\a\otimes_{\Oc_{X^\a}}\hat\Fc^\a).\leqno(4.3.5)$$
See \cite{BD2, \S 7.11.4} for more details on $\Ob_{X^\infty}$, $\hat\Ob_{
X^\infty}$.

We set also $\Db_{X^\infty/S}=2\ind_\a (\Db_{X^\a/S}, i_{\a\b\bullet}).$
Let $i_\a$ be the embedding $X^\a\hookrightarrow X^\infty$.
Let $i_{\a\bullet}\,:\,\Db_{X^\a/S}\to\Db_{X\infty/S}$,
$\Mc^\a\mapsto(\a,\Mc^\a)$ be the canonical functor.
It is exact. There is also the functor
$$\Db_{X^\infty/S}\to\Ob_{X^\infty},\quad
\Mc=i_{\a\bullet}\Mc^\a\mapsto
\Mc^\Oc=(i_{\a\bullet}\Mc^\a)^\Oc:=
(i_{\a\b\bullet}\Mc^\a)^\Oc_{\b\geq\a},\leqno(4.3.6)$$
see (4.1.3). In particular, to any $\Mc\in\Db_{X^\infty/S}$
we can associate its direct image to $S$ :
if $\Mc$ is represented by
$\Mc^\a\in\Db_{X^\a/S}$, then
$$q_*(\Mc):=
q_*(\Mc^\Oc)=
\ind_{\b\geq\a}q_*^\b\bigl((i_{\a\b\bullet}\Mc^\a)^\Oc\bigr).
\leqno (4.3.7)$$

\vskip3mm

\noindent{\bf (4.3.8) Remark.}
$(a)$ We have an exact functor
$$2\ind_\a\Ob_{X^\a}\to\Ob_{X^\infty},\quad
(\a,\Ec^\a)\mapsto(i_{\a\b*}\Ec^\a)_{\b\geq\a}.$$
The two categories are not equivalent in general.

$(b)$ The category $\Ob_{X^\infty}$ is closed by inductive limits.

$(c)$ If the ind-scheme $X^\infty$ is not discrete anymore,
the $\Oc_{X^\a}$-module $i_{\a\b}^!\Ec^\b$ may not be quasi-coherent.
However the category $\Ob_{X^\infty}$ is still well-defined.

\vskip .2cm

Let $X^\infty = \ind_\a X^\a$ be a formally smooth (over $S$)
discrete ind-$S$-scheme.
Following [BD2, Proposition 7.11.8] we define the tangent sheaf
of $X$ to be the object $\hat\Theta_{X^\infty/S}
\in\hat\Ob_{X^\infty}$ such that
$$\Ec\otimes \hat\Theta_{X^\infty/S} =
\ind_\a \Cal Hom (\Omega^1_{X^\a/S}, \Ec^\a),\quad
\forall \Ec\in\Ob_{X^\infty}.$$
Denoting $i^\a: X^\a\to X^\infty$ the canonical embedding,
we have that $i^{\a*}\hat\Theta_{X^\infty/S}$ is a locally free
$\Oc_{X^\a}$-module (possibly of infinite rank),
see \cite{BD2, Proposition 7.12.13}.

\vskip .3cm

\noindent{\bf (4.4) $\Dc$-modules over ind-pro-schemes.}
The following definition is inspired by the paper of K. Kato [Kat].

\proclaim{(4.4.1) Definition}
An ind-$S$-scheme $X^\infty_\infty$ is called locally compact if
it can be represented as
$$X^\infty_\infty = \ind_{\alpha\in A} \pro_{\beta\in B} X^\alpha_\beta$$
where $(X^\alpha_\b)$ is a bi-filtering
ind-pro-system over $\Schb_S^\ft$ with the following
properties:

(1) For each $\beta\in B$ and $\a\leq \a'\in A$
the structure map $i_\b^{\a\a'}: X^\a_\b\to
X^{\a'}_\b$ is a closed embedding.

(2) For each $\a\in A$ and $\b\leq \b'\in B$ the structure map
$p^\a_{\b\b'}: X^\a_{\b'}\to X^\a_\b$ is an affine morphism.

(3) For each $\a\leq \a'\in A$ and $\b\leq\b'\in B$ the commutative
square
$$\matrix
X^\a_{\b'}&\hookrightarrow&X^{\a'}_{\b'}\cr
\downarrow&&\downarrow\cr
X^\a_{\b}&\hookrightarrow&X^{\a'}_\b\cr
\endmatrix$$
is Cartesian.

We denote by $\Lcb_S\subset\Ischb_S$ the full subcategory
formed by locally compact ind-$S$-schemes.

\endproclaim

Let $(X^\a_\b)$ be an ind-pro-system as above.
The maps $p^\a_{\b\b'}$ being affine, the projective
limit $\pro_\b X^\a_{\b}$ is represented by a compact scheme
 (possibly of infinite type) denoted by $X^\a_\infty$. Similarly, we have the discrete
ind-schemes $X^\infty_\b = \ind_\a X^\a_\b$. By definition
$X^\infty_\infty = \ind_\a X^\a_\infty$.

\proclaim {(4.4.2) Proposition} For
an ind-pro-system satisfying the conditions of (4.4.1) we also
have $X_\infty^\infty = \pro_\beta X^\infty_\beta = \pro_\b\ind_\a X^\a_\b$.

\endproclaim

\noindent {\sl Proof:} By passing to the  functors represented by our (ind-)schemes,
we reduce the statement (a)
to the following lemma whose proof we leave
to the reader.

\proclaim {(4.4.3) Lemma}
Let $(T^\a_\b)_{\a\in A, \b\in B}$ be a bi-filtering ind-pro-system of sets.
Then there is a canonical map
$$c: \ind_\a\pro_\b T^\a_\b \to \pro_\b\ind_\a T^\a_\b.$$
If, moreover, all the squares in $(T^\a_\b)$ are Cartesian,
then $c$ is an isomorphism.
\endproclaim

\proclaim{(4.4.4) Definition} We say that a locally compact
ind-$S$-scheme $X^\infty_\infty$ is smooth (over $S$)
if it admits a presentation as in (4.1.1) where :

(1) All the $p^\a_{\b\b'}$
are smooth morphisms of relative dimension $d_{\b\b'}$
(independent on $\a$).
There is an element $(\a,\b)\in A\times B$
such that $X^\a_\b$ is smooth over $S$.

(2) All the ind-$S$-schemes $X_\b^\infty = \ind_\a X^\a_\b$ are formally
smooth over $S$.
\endproclaim

Let $X^\infty_\infty$ be a locally compact smooth ind-$S$-scheme
and $(X^\a_\b)$ be an ind-pro-system as in (4.4.1), (4.4.4).
By Lemma 4.1.4, we have then a double inductive system of categories
$(\Db_{X^\a_\b/S}, i_{\b \bullet}^{\a\a'}, p^{\a\bullet}_{\b\b'})$
and we define the category of (right) $\Dc$-modules
on $X^\infty_\infty$ to be $\Db_{X^\infty_\infty/S} = 2\ind_{\a,\b}
\Db_{X^\a_\b/S}$.

\proclaim{(4.4.5) Proposition}
The category $\Db_{X^\infty_\infty/S}$ is independent,
up to canonical equivalence, on the choice of $(X^\alpha_\beta)$ as in (4.4.4).
\endproclaim

\noindent{\sl Proof :}
Each $X^\alpha_\infty=\pro_\beta X^\alpha_\beta$ being almost smooth, 
the category $\Db_{X^\alpha_\infty/S}=2\ind_{\beta}\Db_{X^\alpha_\beta/S}$
depends, by Proposition 4.2.6, on $X^\alpha_\infty$ only.
Next, for $\alpha<\alpha'$ the functor
$\Db_{X^\alpha_\infty/S}\to\Db_{X^{\alpha'}_\infty/S}$
depends only on the morphism
$X^\alpha_\infty\to X^{\alpha'}_\infty$.
This is seen by the same argument as in (4.2.6).
Let $(\tilde X^{\tilde\alpha}_{\tilde\beta})$ 
be another ind-pro-system as in (4.4.4) representing $X_\infty^\infty$.
So $X^\infty_\infty=\ind_\alpha X^\alpha_\infty=
\ind_{\tilde\alpha}\tilde X^{\tilde\alpha}_\infty$.
The second equality (of ind-objects) means that each $X^\alpha_\infty$ 
is included into some
$\tilde X^{\tilde\alpha}_\infty$
as a closed subset and vice versa.
This means that $\Db_{X^\alpha_\infty/S}$ 
is identified with a full subcategory
in some $\Db_{\tilde X^{\tilde\alpha}_\infty/S}$ and vice versa.
Therefore their 2-limits are identified.
\qed

\vskip3mm

\noindent {\bf (4.4.6) Remark.} Although we have defined
$\Db_{X^\infty_\infty/S}$ as an
abstract category, it is impossible, in general,   to view its objects
as sheaves in a more conventional sense.
For example, it is impossible to associate to
an object of $\Db_{X^\infty_\infty/S}$ its direct image onto $S$.
Indeed, assume that $S=\Spec(k)$ for simplicity.
If such an object $\Mc$ is represented by some $\Mc^\a_\b\in\Db_
{X^\a_\b}$, then the spaces of global sections of
the coherent sheaves
$$(p^{\a'\bullet}_{\b\b'}i_{\b\bullet}^{\a\a'}\Mc^\a_\b)^\Oc =
\bigl(p^{\a'*}_{\b\b'}(i_{\b\bullet}^{\a\a'}\Mc^\a_\b)^\Oc\bigr)
\otimes_{\Oc_{X^{\a'}_{\b'}}}
\omega_{X^{\a'}_{\b'}/X^{\a'}_\b}\in\Ob_{X^{\a'}_{\b'}}$$
do not form an inductive system because of the twist by the relative
canonical class.
For a compact smooth pro-scheme it is possible
to get around this problem by untwisting by
the absolute canonical classes of the terms of the projective system,
see (4.2.7).
To achieve the same effect in the ind-pro-case one would need
to make sense of the (absolute) canonical class of the ind-scheme $X^\infty_\b$,
i.e. of the determinant of the (possibly infinite-dimensional)
vector bundle $\hat\Theta_{X^\infty_\b}$. The impossibility
of doing this (``the determinantal anomaly'')
is precisely the reason why there is no natural space
$\Gamma(X_\infty^\infty,\Mc)$.

\vskip3mm

\noindent {\bf (4.4.7) Example.} 
(a)
If $X$ is a smooth affine variety
admitting an \'etale map to $\AA^d$, then the ind-scheme $\Lc(X)$
is locally compact and smooth.
Thus, the category $\Db_{\Lc(X)}$ is well-defined.
If $X$ is no longer affine it admits a covering by affine open sets $U_\a$
admitting an \'etale map to $\AA^d$.
Then, an object of $\Db_{\Lc(X)}$ is a sheaf on $\Lc(X)$ whose restriction to
$\Lc(U_\a)$ is an object of $\Db_{\Lc(U_\a)}$.

(b)
If $X$ is smooth then
$\Lc(X)_{C^I}\to C^I$ is a locally compact and smooth ind-scheme over $C^I$.

(c)
If $X$ is a smooth affine variety, then the ind-scheme $\tilde\Lc(X)$ is locally
compact and formally smooth. But we do not know if it is smooth in the sense
of (4.4.4).

\vfill\eject

\head 5.  De Rham complexes on ind-schemes \endhead

\vskip.3cm

As in \S 4, let $S$ be a $k$-scheme of finite type.

\vskip.3cm

\noindent{\bf (5.1) Reminder on the De Rham complexes.}
Let $q\,:\,X\to S$ be a smooth $S$-scheme of finite type 
and $\Mc\in\Db_{X/S}$ be a right $\Dc_{X/S}$-module.
Its de Rham complex $\Dc\Rc(\Mc)$ is given by
$$\Dc\Rc^{i}(\Mc)=\Cal Hom_{\Oc_X}(\Omega_{X/S}^{-i},\Mc)=\Mc\otimes_{\Oc_X}
\bigwedge\nolimits^{-i}\Theta_{X/S},\quad i\leq 0.\leqno (5.1.1)$$
If $(x_j)$ is a relative \'etale coordinate system
on an open part of $X$, then the  differential is given by the formula
$d=\sum \partial_{x_j}\otimes dx_j$ where $dx_j\in\Omega^1_X$ is considered
as the contraction operator $\bigwedge^{-i}\Theta_X\to\bigwedge^{-i-1}
\Theta_X$. We denote by $DR(\Mc)=q_*(\Dc\Rc(\Mc))$ the complex
of direct images. 

Let  $i: X\to Y$ be a closed embedding of smooth $S$-schemes 
and $\Mc\in\Db_{X/S}$.
The embedding (4.1.3)  induces
an embedding of the de Rham complexes
$\Dc\Rc(\Mc) \hookrightarrow i^* \Dc\Rc(i_\bullet\Mc)$
and therefore an embedding of the complexes of direct images 
$$DR(\Mc) \hookrightarrow DR(i_\bullet
(\Mc)).\leqno (5.1.2)$$

Let $p: X\to Y$ be a smooth morphism of smooth $S$-schemes of
relative dimension $d$ and $\Mc\in\Db_{Y/S}$.  
Let $q\,:\,X\to S$, $r\,:\,Y\to S$ be the structure maps.
The embedding (4.1.2) induces
an embedding of de Rham complexes which now involves a shift in the degrees:
$$p^*\Dc \Rc(\Mc) \hookrightarrow \Dc\Rc(p^\bullet(\Mc))[d].$$
It is induced by the map
$$p^*\bigwedge\nolimits^{i}\Theta_{Y/S}\to
\omega_{X/Y}\otimes_{\Oc_X}\bigwedge\nolimits^{i+d}\Theta_{X/S}.$$
In particular, we get an embedding of the complexes of direct images 
$$DR(\Mc)\hookrightarrow DR(p^\bullet\Mc)[d].
\leqno (5.1.3)$$
Note that without passing to the de Rham complexes there is no embedding
of $r_*\Mc$ into $q_*p^\bullet\Mc$.
The map (5.1.3) can be seen as a $\Dc$-module manifestation
of the fact that ``fermions cancel the determinantal anomaly''.

\vskip .2cm

\noindent {\bf (5.1.4) Example.} 
Take $S=\Spec(k)$.
Let $Y=\AA^d$ with coordinates $a_1, ...,
a_d$ and $i: X\hookrightarrow Y$ be the embedding of the
affine subspace $\{ a_1 = ... = a_l=0\}$, $l\leq d$.
The algebra $\Gamma(Y, \Dc_Y)$ is just the Heisenberg (Weyl)
algebra $D_Y$ generated by $a_1, ..., a_d$ and $a_1^*, ..., a_d^*$
subject to the relations
$$[a_m, a_n]= [a_m^*, a_n] = 0, \quad [a_m^*, a_n] = \delta_{mn}.$$
The space of global sections $\Gamma(Y, i_\bullet \omega_X)$ is the right
$D_Y$-module
$$\omega_{XY} = D_Y\big/ (a_m^*, a_n; \, n\leq l < m) D_Y.$$
Let also $C_Y$ be the Clifford algebra generated by odd elements $b_1, ...,
b_d$, $b_1^*, ..., b_d^*$ subject to the relations
$$[b_m, b_n]_+ = [b^*_m, b^*_n]_+ = 0, \quad [b^*_m, b_n]_+ = \delta_{mn}.$$
Denote $CD_Y = C_Y\otimes_k D_Y$ the tensor product algebra.
Then the global de Rham complex of $i_\bullet\omega_X$
is identified with the right $CD_Y$-module
$$DR_{XY} = CD_Y\big/
\bigl(a_m^*, a_n, b_p;\, n\leq l<m,\, p=1, ...,d\bigr) CD_Y.$$

\vskip .3cm

\noindent{\bf (5.2) De Rham complexes for ind-schemes.} Let
$X^\infty = \ind_{\a\in A} X^\a$ be a formally smooth discrete
ind-$S$-scheme with structure maps $q,q^\a$. 
Denote by $i^\a: X^\a\hookrightarrow X^\infty$ the canonical
embedding. The considerations of (5.1) generalize easily to give
the global  de Rham complex of any $\Mc\in\Db_{X^\infty/S}$. Explicitly,
let $\Mc$ have the form $i_{\a\bullet}\Mc^\a$. Then the $i$th term
of its de Rham complex is
$$DR^{i}(\Mc)=\Mc^\Oc\otimes\bigwedge\nolimits^{-i}\hat\Theta_{X^\infty/S}=
\ind_{\a'\geq\a}q^{\a'}_*\Hc om(\Omega^{-i}_{X^{\a'}/S},
(i^{\a\a'}_{\bullet}\Mc^\a)^\Oc).\leqno (5.2.1)$$
Here $\Mc^\Oc$ is defined in (4.3.6).

\proclaim{(5.2.2) Proposition} Let $p^\infty: X^\infty\to Y^\infty$ be
a morphism of formally smooth discrete ind-$S$-schemes
which is smooth of relative dimension $d$.
Then for any $\Mc\in\Db_{Y^\infty/S}$ we have an embedding of
the shifted de Rham complexes
$DR(\Mc)\hookrightarrow DR(p^{\infty\bullet}\Mc)[d].$
\endproclaim

\noindent{\sl Proof:} The conditions on $p^\infty$ in the proposition
are equivalent to the following: we can represent $X^\infty = \ind_{\a\in A}
X^\a$, $Y^\infty = \ind_{\a\in A} Y^\a$ with the same filtering poset $A$,
and we can represent $p^\infty$ by a morphism of inductive systems
$(p^\a: X^\a\to Y^\a)$ of $S$-schemes
such that for each $\a\leq \b$ the arising
commutative square
$$\matrix
X^\a&{\buildrel i^{\a\b}\over\hookrightarrow}&X^\b\cr
{\ss p^\a}\downarrow&&\downarrow{\ss p^\b}\cr
Y^\a&{\buildrel j^{\a\b}\over\hookrightarrow}&Y^\b\cr
\endmatrix$$
is Cartesian, and each $p^\a$ is smooth of relative dimension $d$.
Then $p^\infty$ is formally smooth.
By \cite{BD2, Lemma 7.12.23} we have an exact sequence
$$0\to\hat\Theta_{X^\infty/Y^\infty}\to
\hat\Theta_{X^\infty/S}\to p^{\infty*}\hat\Theta_{Y^\infty/S}\to 0.
\leqno(5.2.3)$$
Moreover
$$(i^\b)^*(\hat\Theta_{X^\infty/Y^\infty})=\Theta_{X^\b/Y^\b}.\leqno(5.2.4)$$
Let $\Mc$ have the form $j_{\a\bullet}\Mc^\a$,
where $j^\a\,:\,Y^\a\hookrightarrow Y^\infty$ is the canonical embedding.
Then $p^{\infty\bullet}\Mc=i_{\a\bullet}p^{\a\bullet}\Mc^\a$.
The base change for Cartesian squares gives
$$\matrix
DR^{i}(p^{\infty\bullet}\Mc)[d]
&=\ind_{\b\geq\a}q^\b_*\bigl((p^{\b\bullet}j^{\a\b}_{\bullet}\Mc^\a)^\Oc
\otimes_{\Oc_{X^\b}}i^{\b*}\bigwedge\nolimits^{-i-d}\hat\Theta_{X^\infty/S}
\bigr)
\hfill\cr
&=\ind_{\b\geq\a}q^\b_*\bigl(
p^{\b*}(j^{\a\b}_{\bullet}\Mc^\a)^\Oc\otimes_{\Oc_{X^\b}}\o_{X^\b/Y^\b}
\otimes_{\Oc_{X^\b}}i^{\b*}\bigwedge\nolimits^{-i-d}\hat\Theta_{X^\infty/S}
\bigr),
\hfill\cr
\endmatrix$$
see Lemma 4.1.4.
Let $r\,:\,Y^\infty\to S$, $r^\a\,:\,Y^\a\to S$ be the structure maps.
By (5.2.3), (5.2.4) there is an embedding
$$p^{\b*}j^{\b*}\bigwedge\nolimits^{-i}\hat\Theta_{Y^\infty/S}=
i^{\b*}p^{\infty*}\bigwedge\nolimits^{-i}\hat\Theta_{Y^\infty/S}
\to\o_{X^\b/Y^\b}
\otimes_{\Oc_{X^\b}}i^{\b*}\bigwedge\nolimits^{-i-d}\hat\Theta_{X^\infty/S}.$$
Hence there is an embedding
$$DR^i(\Mc)=\ind_{\b\geq\a}r^\b_*\bigl(
(j^{\a\b}_{\bullet}\Mc^\a)^\Oc\otimes_{\Oc_{Y^\b}}
j^{\b*}\bigwedge\nolimits^{-i}\hat\Theta_{Y^\infty/S}\bigr)\to$$
$$\to\ind_{\b\geq\a}q^\b_*\bigl(
p^{\b*}\bigl((j^{\a\b}_{\bullet}\Mc^\a)^\Oc\otimes_{\Oc_{Y^\b}}
j^{\b*}\bigwedge\nolimits^{-i}\hat\Theta_{Y^\infty/S}\bigr)\bigr)
\to DR^{i}(p^{\infty\bullet}\Mc)[d].$$
We are done.
\qed

\vskip .3cm

\noindent {\bf (5.3) De Rham complexes for ind-pro-schemes.}
Let $X^\infty_\infty$ be a locally compact smooth ind-$S$-scheme
and $\Mc$ be an object of $\Db_{X^\infty_\infty/S}$.
We fix an ind-pro-system $(X^\a_\b)$ for $X$ as in (4.4.1), (4.4.4).
We have then the formally smooth discrete ind-schemes $X_\beta^\infty$
and projections $p_\b: X^\infty_\infty\to X_\b^\infty$.
We also have schemes $X^\a_\infty$ and embeddings 
$i^\a\,:\,X^\a_\infty\to X^\infty_\infty$.
Let also $i^\a_\b\,:\,X^\a_\b\to X^\infty_\b$,
$p^\a_\b\,:\,X_\infty^\a\to X^\a_\b$ be the natural embeddings and projections.
The category $\Db_{X^\infty_\infty/S}$ being the double direct limit
of $\Db_{X^\a_\b/S}$, we can think of $\Mc$ as being of the form
$p_\b^\bullet i^\a_{\b\bullet}\Mc^\a_\b=i_\bullet^\a p_\b^{\a\bullet}\Mc^\a_\b$
for some $\Mc^\a_\b\in\Db_{X^\a_\b/S}$.
Recall that $d_{\b\b'}$ denotes the relative dimension of the smooth
morphism $p^\a_{\b\b'}: X^\a_{\b'}\to X^\a_\b$, $\b\leq \b'$.
We choose numbers $d_\b, \b\in B$, such that $d_{\b\b'} = d_\b-d_{\b'}$
(this can be done uniquely up to an overall constant).
Set $\Mc_\b=i^\a_{\b\bullet}\Mc^\a_\b\in\Db_{X^\infty_\b/S}$.
Proposition 5.2.2 implies then that the shifted global de Rham complexes
$DR(p_{\b\b'}^\bullet \Mc_\b)[d_{\b'}]$ form an inductive system
of complexes of vector spaces and we define the de Rham complex
of $\Mc$ to be
$$DR(\Mc) = \ind_{\beta'\geq\beta} DR(p_{\b\b'}^\bullet \Mc_\b)[d_{\b'}].
\leqno (5.3.1)$$
Explicitely, by (5.2.1) we have
$$DR^{i}(\Mc)=\ind_{\b'\geq\b}\ind_{\a'\geq\a}
q^{\a'}_{\b'*}\Hc om\bigl(\Omega^{-i-d_{\b'}}_{X^{\a'}_{\b'}/S},
(p^{\a'\bullet}_{\b\b'}i^{\a\a'}_{\b\bullet}\Mc^\a_\b)^\Oc)\bigr),$$
where $q^{\a'}_{\b'}\,:\,X^{\a'}_{\b'}\to S$ is the structure morphism.

\proclaim{(5.3.3) Proposition}
$DR^i(\Mc)$ depends only on $X^\infty_\infty$ 
and $\Mc$ as an object of $\Db_{X^\infty_\infty/S}$,
but not on the choice of a system $(X^\a_\b)$.
\endproclaim

\noindent{\sl Proof:}
Interchanging the two inductive limits and using base change for Cartesian 
squares in the diagram $(X^\a_\b)$, we can write
$$DR^{i}(\Mc)=\ind_{\a'\geq\a}\ind_{\b'\geq\b}
q^{\a'}_{\b'*}\Hc om\bigl(\Omega^{-i-d_{\b'}}_{X^{\a'}_{\b'}/S},
(i^{\a\a'}_{\b'\bullet}p^{\a\bullet}_{\b\b'}\Mc^\a_\b)^\Oc)\bigr).$$
For any $\a'\geq\a$ the limit over $\b'$ depends only on the scheme 
$X^{\a'}_\infty$ and the object
$$\Mc^{\a'}=p_\b^{\a'\bullet}i^{\a\a'}_{\b\bullet}\Mc^\a_\b
\in\Db_{X^{\a'}_\infty/S}.$$
Therefore the limit over $\a'\geq\a$ of the limits above depends only on the
ind-object
$``\ind_{\a'\geq\a}"X^{\a'}_\infty$ (which is $X^\infty_\infty$) and
the object 
$$i_{\a'\bullet}\Mc^{\a'}\in 2\ind_{\a'\geq\a}\Db_{X^{\a'}_\infty/S}=
\Db_{X^\infty_\infty/S},$$
which is $\Mc$.
\qed

\vskip .3cm

\noindent {\bf (5.4) The de Rham complexes on $\Lc(X)$.}
We now specialize to the particular case $X^\infty_\infty = \Lc(X)$
where $X$ is a smooth affine algebraic variety over $k$ admitting an
\'etale map $\phi$ to $\AA^d$. In this case $A= \EE, B=\NN$
with the terms of the ind-pro-systems being $\Lc^\eps_n(\phi)$,
see (2.6.7).
We take $S=\Spec(k)$.
Given an object $\Mc\in\Db_{\Lc(X)}$ we associate to it
its de Rham complex $DR(\Mc)$ as in (5.3).
Note that it is independent on the choice
of an \'etale map to $\AA^d$ because two such
maps lead to isomorphic ind-pro-objects
in the category $\Schb^\ft$, see Proposition 2.7.1.$(c)$.

Let now $X$ be an arbitrary smooth algebraic variety over $k$
and $\Mc$ be an object of $\Db_{\Lc(X)}$.
By covering $X$ with affine open $U$ admitting \'etale maps
to $\AA^d$, we get a complex of sheaves $U\mapsto DR(\Mc|_{\Lc(U)})$
which we denote $\Dc\Rc(\Mc)$. Recall that we have the diagram
$$X\buildrel p\over\longleftarrow \Lc^0(X) \buildrel i\over\hookrightarrow
\Lc(X).\leqno (5.4.1)$$
Thus, every right $\Dc_X$-module $\Nc$ gives an object
$i_\bullet p^\bullet\Nc$ of $\Db_{\Lc(X)}$.
We write $\Cc\Dc\Rc (\Nc)$ for $\Dc\Rc(i_\bullet p^\bullet
\Nc)$ and call it the chiral de Rham complex of $\Nc$.
In particular, we write
$\Cc\Dc\Rc_X$ for $\Dc\Rc(i_\bullet p^\bullet\omega_X)$.
More generally, denoting $p_n: \Lc^0(X)\to \Lc^0_n(X)$ the
projection, we can start with any right $\Dc$-module
$\Nc$ on the algebraic variety $\Lc^0_n(X)$: then
$i_\bullet p_n^\bullet \Nc$ is an object of $\Db_{\Lc(X)}$
and we can form its de Rham complex.
It is a complex of sheaves on $X$.

\vskip .2cm

\noindent {\bf (5.4.2) Example.} 
Let $X=\AA^1$. Then the complex of global sections
of $\Cc\Dc\Rc_X$, i.e., the complex $DR(i_\bullet p^\bullet\omega_{\AA^1})$
can be found explicitly as follows.

Let $V$ be the topological $k$-vector space $k((t))$ and
$V^*$ be its topological dual (over $k$).
Denote by $=\langle l, v\rangle$ the canonical
pairing of $l\in V^*$ and $v\in V$.

Then $V^*$ can be identified with $k((t))dt$,
the space of 1-forms, the pairing between $V$ and $V^*$
being $(f, \o)\to \op{res}(f\cdot\o)$. Let $D$
be the Heisenberg algebra generated by $V^*$ and $V$
with $[l, v] = \langle l,v\rangle$ and $C$ be the
Clifford algebra generated by $V^*, V$ with
$[l,v]_+=\langle l,v\rangle$. 
Denote $CD = C\otimes_k D$.
This is a certain completion of the algebra $\tilde{CD}$ generated
by symbols $a_n,b_n,a_n^*,b_n^*$ for $n\in\ZZ$ subject to the relations
$$[a_m, a_n] = [a^*_m, a^*_n]=0, \quad [a^*_m, a_n] = \delta_{m, -n},$$
$$[b_m, b_n]_+ = [b^*_m, b^*_n]_+ = 0, \quad
[b^*_m, b_n]_+ = \delta_{m, -n},$$
$$[a_m, b_n] = [a^*_m, b_n] = [a_m, b^*_n] =
[a^*_m, b^*_n]=0.$$
More precisely, we write a generic element of $V$ as $\sum_ma_mt^m$, so
$a_m,a_m^*$ are elements of $D$. Similarly, writing a generic element 
of $V^*$ as $\sum_mb_mt^{m-1}dt$ we view $b_m,b_m^*$ as elements of $C$.
Let $\tilde V=k[t,t^{-1}]$, $\tilde V^*=k[t,t^{-1}]dt$.
Then $\tilde{CD}=\tilde C\otimes_k\tilde D$, where $\tilde C$ is the
Clifford algebra generated by $V^*,V$ 
and $\tilde D$ is the Heisenberg algebra generated by $\tilde V,\tilde V^*$.
Let $\tilde{CD}^+\subset\tilde{CD}$ be the
right ideal generated by linear combinations of
$a_n, b_n, a^*_{n+1}, b^*_{n+1}$ with $n\geq 0$,
and $CD^+\subset CD$ be the ideal generated by possibly infinite linear combinations with the above property.
We denote $Vac=CD/CD^+$ and $\tilde{Vac}=\tilde{CD}/\tilde{CD}^+$
the corresponding vacuum modules.

\proclaim{(5.4.3) Proposition}
(a)
The natural morphism $\tilde{Vac}\to Vac$ is an isomorphism.

(b)
The de Rham complex $DR(i_\bullet p^\bullet\omega_{\AA^1})$ 
is identified (as a vector space) with $Vac.$
\endproclaim

\noindent {\bf Remark.} Although this description is similar to
(5.1.4), there is a difference :
here the ideal is generated by $b_{n+1}^*, b_{n}$
for $n<0$ while in (5.1.4) all the $b_n$ are in the ideal. This is because
in our present situation we are dealing with a semiinfinite de
Rham complex obtained as an inductive limit of usual
de Rham complexes with respect to maps shifting the degrees.

\vskip .2cm

\noindent {\sl Proof:}
(a) The quotient $k((t))/k[[t]]$ is identified with $k[t,t^{-1}]/k[t]$.
The ideal $CD^+$ includes the Taylor series part of $V,V^*\subset D$
and $V,V^*\subset C$.
So $CD/CD^+$ is identified with $\tilde{CD}/\tilde{CD}^+$.

(b) To simplify, we write $\Lc=\Lc(\AA^1)$, etc.
We have
$$\Lc^\eps_M=\Spec\bigl(k[a_l;-N_\eps\leq l\leq M]\big/(a_l^{1+\eps_l})\bigr),$$
where, for each $\eps\in\EE$, we set $N_\eps=\max\{l;\eps_{-l}\neq 0\}$.
Hence
$$\Lc_M=\ind_N\Spf\bigl(k[a_l;0\leq l\leq M][[a_l;-N\leq l<0]]\bigr).$$
For any $M\in\NN$, $N\in\NN\cup\{\infty\}$, we put
$$Y^N_M=\Spec\bigl(k[a_l;-N\leq l\leq M]\bigr).$$
Then $\Lc_M$ is just the limit, over $N>0$,
of the formal completions of $Y^N_M$ along $Y^0_M$.
Since the de Rham complex with coefficients
in the $\Dc$-module of distributions along a subvariety
depends only on the completion along this subvariety,
we can write
$$DR(i_\bullet p^\bullet\omega_{\AA^1}) =
\ind_{M, N} DR(i_{M, N,\bullet} \omega_{Y^0_M})[M],
\leqno (5.4.4)$$
where $i_{M, N}: Y^0_M\hookrightarrow Y^N_M$ is the
embedding. Note that $i_{M, N}$ is just the embedding
of an affine subspace, so we are in the situation
of Example 5.1.4.

Let $D^N_M$ be the subalgebra in $CD$
generated by $a_n,$ $-N\leq n\leq M$,
and $a^*_n$, $-M\leq n\leq N$.
It is identified with
the algebra of polynomial
differential operators on functions
of $a_{-N}, ..., a_M$, with $a^*_{-i}$ corresponding
to $\partial/\partial a_i$. 
Similarly let $C^N_M$ be 
generated by $b_n,$ $-N\leq n\leq M$,
and $b^*_n$, $-M\leq n\leq N$.
Denote $CD_M^N=C_M^N\otimes_k D_M^N$.
We see therefore that
the $(M, N)$th term of the inductive
system in (5.4.4) is identified with
$$Vac_M^N := CD_M^N\big/ (a_i, i>0;\, a_i^*, i\geq 0;\,b_i^*) CD_M^N.
\leqno (5.4.5)$$
Denote by $1^N_M$ the generator of this module.
Then for $N\leq N'$ and any $M$ the embedding takes $1_M^N$
to $1_M^{N'}$, while for $M\leq M'$ and any $N$ it takes
$1^N_M$ into $1_{M'}^N b_{M+1}...b_{M'}$.
From the normal form of elements it is clear that $\tilde{CD}=\ind_{M,N}CD^N_M$,
$\tilde{Vac}=\ind_{M,N}Vac^N_M$,
and we are done since $\tilde{Vac}=Vac$.
\qed

\vfill\eject

\head 6. Identification of the chiral de Rham complex
\endhead

\vskip 3mm

In this section we construct, in a geometric way,
the structure of a vertex algebra on the
chiral de Rham complex $\Cc\Dc\Rc_X$ and
compare it with the construction of [MSV].

\vskip .3cm

\noindent{\bf (6.1) Factorization algebras, and
De Rham complexes on $\Lc(X)_{C^I}$.}
Let $C$ be a smooth curve, as before.
For any non-empty finite set $I$ we set
$U^{(I)}=U^{(I/I)}$, $\Delta^{(I)}=\Delta^{(I/\{1\})}$
and $j^{(I)}=j^{(I/I)}$, see (3.2).
Hereafter we write $U,\Delta,j$ instead of
$U^{(I)}$, $\Delta^{(I)},$ $j^{(I)}$ if $I$ has cardinal 2.
We will need the following notation : for any (possibly empty) $I$ let
$\bar I=I\sqcup\heartsuit$ be the corresponding pointed set.
For any surjection 
$I\twoheadrightarrow J$ we denote 
$\bar I\twoheadrightarrow\bar J$ the surjection equal to 
$I\twoheadrightarrow J$ on $I$ and taking $\heartsuit$ to $\heartsuit$. 
Let us recall the definition of a factorization algebra,
see \cite{BD1, \S 3.4}.

\proclaim{(6.1.1) Definition}
(a) Let $\Ec$ be a quasi-coherent sheaf on $C$.
A  structure of a
factorization algebra on $\Ec$ is a collection
of quasi-coherent $\Oc_{C^I}$-modules $\Ec_I$
for each non-empty finite set $I$,
such that $\Ec_I$ is flat along the diagonal strata,
$\Ec_{\{1\}}=\Ec$, and

--  an isomorphism of $\Oc_{C^I}$-modules
$\nu^{(J/I)}\,:\,\Delta^{(J/I)*}\Ec_J\simto\Ec_I$
for every $J\twoheadrightarrow I$, compatible with the compositions
of $J\twoheadrightarrow I$,

-- an isomorphism of $\Oc_{U^{(J/I)}}$-modules
$$\varkappa^{(J/I)}\,:\,j^{(J/I)*}(\boxtimes_I\Ec_{J_i})
\simto j^{(J/I)*}\Ec_{J}$$
for every $J\twoheadrightarrow I$, compatible with the compositions of
$J\twoheadrightarrow I$, and compatible with $\nu$,

--  a global section $1_\Ec\in\H^0(C,\Ec)$ such that for every
$f\in\Ec$ one has
$1_\Ec\boxtimes f\in\Ec_{\{1,2\}}\subset j_*j^*(\Ec\boxtimes\Ec)$
and $\Delta^*(1_\Ec\boxtimes f)=f$.

(b) A module over $\Ec$ is a quasi-coherent sheaf $\Mc$ on $C$ with
a collection of quasi-coherent $\Oc_{C^{\bar I}}$-modules $\Mc_{\bar I}$
for each non-empty finite set $I$,
such that $\Mc_{\bar I}$ is flat along the diagonal strata,
$\Mc_{\{\heartsuit\}}=\Mc$, and

--   an isomorphism of $\Oc_{C^{\bar I}}$-modules
$\nu^{(\bar J/\bar I)}\,:\,\Delta^{(\bar J/\bar I)*}\Mc_{\bar J}\simto\Mc_{\bar I}$
for every $\bar J\twoheadrightarrow\bar I$,
compatible with the compositions of $\bar J\twoheadrightarrow\bar I$,

--   an isomorphism of $\Oc_{U^{(\bar J/\bar I)}}$-modules
$$\varkappa^{(\bar J/\bar I)}\,:\,
j^{(\bar J/\bar I)*}\bigl((\boxtimes_I\Ec_{J_i})
\boxtimes\Mc_{\bar J_\heartsuit}\bigr)\simto j^{(\bar J/\bar I)*}\Mc_{\bar J}$$
for every $\bar J\twoheadrightarrow\bar I$,
compatible with the compositions of $\bar J\twoheadrightarrow\bar I$,
and compatible with $\nu$, such that

-- for any $f\in\Mc$ one has
$1_\Ec\boxtimes f\in\Mc_{\{1,\heartsuit\}}\subset j_*j^*(\Ec\boxtimes\Mc)$
and $\Delta^*(1_\Ec\boxtimes f)=f$.
\endproclaim

\noindent
We have the following immediate global counterpart of (5.4).

\proclaim{(6.1.2) Proposition}
For any right $\Dc_X$-module $\Mc$ there is a unique complex
$\Cc\Dc\Rc(\Mc)_{C^I}$ of sheaves of $\Oc_{C^I}$-modules
on $X\times C^I$ such that:

(a) the fiber of $\Cc\Dc\Rc(\Mc)_C$ at a point of the curve $C$
is isomorphic to the complex $\Cc\Dc\Rc(\Mc)$,

(b) the collection $\bigl(\Cc\Dc\Rc(\o_X)_{C^I}\bigr)$ is a
factorization algebra on the curve $C$.
The collection $\bigl(\Cc\Dc\Rc(\Mc)_{C^I}\bigr)$
is a $\bigl(\Cc\Dc\Rc(\o_X)_{C^I}\bigr)$-module.
\endproclaim

\noindent{\sl Proof:}
A choice of an element $i_0\in I$ defines a morphism
of schemes $p_{i_0,I}: \Lc^0(X)_{C^I}\to X$ as follows.
Recall that $\Lc^0(X)_{C^I}$
represents the functor $\lambda^0_{X, C^I}$
which takes a scheme $S$ into the set of pairs
$(f_I, \rho)$ where
$f_I: S\to C^I$ is a morphism of schemes and $\rho$
is a morphism of $\widehat{\Gamma(f_I)}$, the formal
neighborhood of $\Gamma(f_I)\subset S\times C$,
into $X$. Now, restricting $\rho$ onto the graph
of $f_{i_0}$, which is a subscheme in
$\widehat{\Gamma(f_I)}$ isomorphic to $S$,
we get a natural transformation from $\lambda^0_{X, C^I}$
into the functor represented by $X$, so a morphism
$p_{i_0,I}$.

Denote by $\iota_{\bar I}: \Lc^0(X)_{C^{\bar I}}\to \Lc(X)_{C^{\bar I}}$ 
the embedding.
For any $\Mc\in\Db_X$, we form the object
$(\iota_{\bar I})_\bullet(p_{\heartsuit,\bar I})^\bullet\Mc\in
\Db_{\Lc(X)_{C^{\bar I}}/C^{\bar I}}$. 
The general construction of (5.3), applied
to the restrictions of 
$(\iota_{\bar I})_\bullet(p_{\heartsuit,\bar I})^\bullet\Mc$
onto open subsets in $X\times C^{\bar I}$, gives then
a complex of sheaves on $X\times C^{\bar I}$ which we denote
$\Cc\Dc\Rc(\Mc)_{C^{\bar I}}$.

Notice that in the particular case where $\Mc=\o_X$ the object
$(p_{i_0,I})^\bullet\Mc$   of the category
$\Db_{\Lc^0(X)_{C^I}/C^I}$ is independent (up to a unique isomorphism)
of the choice of $i_0\in I$. Indeed, objects of the latter
category are, by definition, pairs $(n, \Nc)$
where $\Nc$ is a right $\Dc$-module on $\Lc^0_n(X)_{C^I}$
and two such pairs  $(n, \Nc)$ and $(n', \Nc'$)
are isomorphic, if the pullbacks of $\Nc$ and $\Nc'$ to
$\Lc^0_m(X)_{C^I}$, $m\geq n, n',$ are isomorphic
as right $\Dc$-modules. Since the
pullback for right $\Dc$-modules is just the $\Oc$-module
pullback tensored with the relative canonical class,
$(p_{i_0,I})^\bullet\o_X$ is represented by 
$(n, \o_{\Lc^0_n(X)_{C^I}})$ for any $n$, and thus is clearly
independent on $i_0$.

Then, the general construction in (5.3)
gives a complex of sheaves on $X\times C^I$,
denoted by $\Cc\Dc\Rc(\o_X)_{C^I}$.

To prove Claim $(a)$ it is sufficient to observe that the fiber of
$\Lc(X)_C$ at a point $0\in C$ is isomorphic to $\Lc(X)$.
Recall that $\Lc(X)_C$ represents the contravariant functor
$\l_{X,C}\,:\,\Schb\to\Setsb$ such that $\l_{X,C}(S)$ is the set of pairs
$(f,\rho)$ such that
$$f\in\Hom_{\Schb}(S,C)\and
\rho\in\Hom_{\Lrsb}\bigl((\Gamma(f),\Kc_{f}^\sqr),\,X\bigr).$$
Thus the fiber at $0$ represents the subfunctor
$$S\mapsto\bigl\{(f,\rho)\in\l_{X,C}(S)\,\big|\,f(S)_\red=\{0\}\bigr\}.$$
Let $t$ be a local coordinate on $C$ centered at 0.
For any $f$ as above we have
$\bigl(\Gamma(f),\Kc_{f}^\sqr\bigr)=\bigl(S,\Oc_S((t))^\sqr\bigr)$
and this proves (a).
Note that the isomorphism of $\Lc(X)$ and the fiber of
$\Lc(X)_C$ at 0 is compatible with the ind-pro-systems in (2.7), (3.9).

(b) Both $\Lc^0(X)_{C^I}$ and $\Lc(X)_{C^I}$ form factorization monoids
in the categories of ind-schemes.
Since passing to the De Rham complex takes Cartesian
products of (ind-)schemes to tensor products of
vector spaces, we see that $(\Cc\Dc\Rc(\o_X)_{C^I})$ 
form a factorization algebra.

Next, given any surjection $J\to I$ and the corresponding surjection 
$\bar J\to \bar I$, we have 
$$\prod_{\bar i\in\bar I}\Lc(X)_{C^{\bar J}_{\bar i}}=
\Lc(X)_C\times\prod_{i\in I}\Lc(X)_{C^{J_i}},$$
the first factor in the RHS corresponding to $\bar i=\heartsuit$.
Let us use the notation $\kappa$ for the factorization monoid structure
of $\Lc(X)_{C^I}$ as in (3.2.1).
Then, with respect to the identification (6.1.3), we have an isomorphism of
$\Dc$-modules
$$(\kappa^{(\bar J/\bar I)})^\bullet
(\iota_{\bar I})_\bullet(p_{\heartsuit,\bar I})^\bullet\Mc=
(p_{\heartsuit,\{\heartsuit\}})^\bullet\Mc\otimes\bigotimes_{i\in I}
(\iota_{I_i})_\bullet(p_{\heartsuit,I_i})^\bullet\o_X$$
over $U^{\bar J/\bar I}\to C^{\bar I}$.
Using again the fact that passing to the De Rham complexes takes
Cartesian products to tensor products, we conclude that
$(\Cc\Dc\Rc(\Mc)_{C^{\bar I}})$  is a factorization module over
$(\Cc\Dc\Rc(\o_X)_{C^I})$.
\qed

\vskip3mm

\noindent{\bf (6.2)
 Reminder on chiral and vertex algebras.}
Let us recall the basic facts on chiral and vertex algebras.
See \cite{BD1, \S 3}, \cite{K} and \cite{FLM} for more details.
Let $C$ be a smooth curve, as before.
For any right $\Dc_C$-module $\Mc$ the projection formula yields
an isomorphism of right $\Dc_{C^2}$-modules
$\Delta_\bullet\Delta^\bullet(\o_C\boxtimes\Mc)\simto\Delta_\bullet\Mc.$
Let
$$\eps_\Mc\,:\,j_\bullet j^\bullet(\o_C\boxtimes\Mc)\to\Delta_\bullet\Mc$$
be the composition of the projection
$j_\bullet j^\bullet(\o_C\boxtimes\Mc)\to
j_\bullet j^\bullet(\o_C\boxtimes\Mc)/(\o_C\boxtimes\Mc)$
and of the isomorphism
$\bigl(j_\bullet j^\bullet(\o_C\boxtimes\Mc)\bigr)/(\o_C\boxtimes\Mc)\simto
\Delta_\bullet\Delta^\bullet(\o_C\boxtimes\Mc)
\simto\Delta_\bullet\Mc.$

\proclaim{(6.2.1) Definition}
(a) A chiral algebra over $C$ is a right
$\ZZ/2\ZZ$-graded $\Dc_C$-module $\Ac=\Ac^0\oplus\Ac^1$
with two even maps $\mu_\Ac\in\Hom_{\Db_{C^2}}
\bigl(j_\bullet j^\bullet(\Ac\boxtimes\Ac),\Delta_\bullet\Ac\bigr)$
and $1_\Ac\in\Hom_{\Db_C}(\o_C,\Ac^0)$ such that

-- the map $\mu_\Ac(1_\Ac,id_\Ac)$ coincides with $\eps_\Ac$,

-- the map $\mu_\Ac$ is antisymmetric, and it satisfies the Jacobi identity.

(b) A module over a chiral algebra $\Ac$ over $C$ is a right
$\ZZ/2\ZZ$-graded $\Dc_C$-module $\Mc_C$ with an even map
$\mu_\Mc\in\Hom_{\Db_{C^2}}
\bigl(j_\bullet j^\bullet(\Ac\boxtimes\Mc),\Delta_\bullet\Mc\bigr)$
such that

-- the map $\mu_\Mc(1_\Ac,id_\Mc)$ coincides with $\eps_\Mc$,

-- the map $\mu_\Mc$ is compatible with $\mu_\Ac$.
\endproclaim

\noindent
For any factorization algebra on $\Ec$,
each sheaf $\Ec_I$ has a canonical left $\Dc_{C^I}$-module structure,
compatible with the factorization structure,
such that the section $1_\Ec$ is a horizontal,
see \cite{BD1, Proposition 3.4.8}.
It is proved in \cite{BD1, \S 3.4.9} that the
right $\Dc_C$-module $\Ec^r:=\Ec\otimes_{\Oc_C}\o_C$
is a chiral algebra over $C$.
The map $\mu_\Ac$ is the composition of the chain of maps
$$j_\bullet j^\bullet(\Ec^r\boxtimes\Ec^r)=
j_\bullet j^\bullet\o_{C^2}
\otimes_{\Oc_{C^2}}
\Ec_{\{1,2\}}
\to
\Delta_\bullet\o_C
\otimes_{\Oc_{C^2}}
\Ec_{\{1,2\}}
=\Delta_\bullet\Ec^r.$$
Here the 1-st equality is the 2-nd isomorphism in $6.1.1.(a)$,
the second arrow is $\eps_{\o_C}$, and
the last equality results from the 1-st isomorphism in $6.1.1.(a)$
and the projection formula for $\Delta$.

\proclaim{(6.2.2) Definition}
(a) A vertex algebra is a $k$-supervector space $V$
with an even vector $1_V\in V$,
an even endomorphism $\partial_V\in\End(V)$,
and an even  linear map
$V\to\End(V)[[z,z^{-1}]],\,a\mapsto a(z)=\sum_na_nz^{-n-1}$.
These data satisfy the following axioms:

-- $\partial_V(1_V)=0$, $1_V(z)=id_V$,
$a_n(1_V)=0$ if $n\geq 0$, $a_{-1}(1_V)=a$,

-- $[\partial_V,a(z)]=\partial_za(z)$,

-- we have $(z-w)^N[a(z),b(w)]=0$ for $N>>0$.

\noindent
We will also assume that for all elements $a,b\in V$ we have
$a_n(b)=0$ for $n>>0$.

(b) A module over a vertex algebra $V$ is a
$k$-supervector space $W$ with
an even endomorphism $\partial^W\in\End(W)$,
and an even linear map
$V\to\End(W)[[z,z^{-1}]],\,a\mapsto a^W(z)=\sum_na^W_nz^{-n-1}$.
These data satisfy the following axioms :

--  $1^W_V(z)=id_W$,

-- $[\partial^W,a^W(z)]=\partial_za^W(z)=(\partial_V a)^W(z)$,

--  (Borcherds identity)
$$z_0^{-1}\delta\Bigl({z_1-z_2\over z_0}\Bigr)a^W(z_1)b^W(z_2)-
z_0^{-1}\delta\Bigl({-z_2+z_1\over z_0}\Bigr)b^W(z_2)a^W(z_1)=$$
$$=z_2^{-1}\delta\Bigl({z_1-z_0\over z_2}\Bigr)\bigl(a(z_0)(b)\bigr)^W(z_2),$$
where
$$z_0^{-1}\delta\Bigl({z_1-z_2\over z_0}\Bigr)=
\sum_{m\in\NN}\sum_{n\in\ZZ}(-1)^m\Bigl(\matrix n\cr m\endmatrix\Bigr)
z_0^{-n-1}z_1^{n-m}z_2^m.$$
\endproclaim

\noindent
Assume that $C$ is the formal disk $\Spec k[[t]]$.
Let $0$ be the closed point of $C$.
Let $t_1=t\otimes 1$, $t_2=1\otimes t$ be the coordinates on $C^2$.
We have the following basic fact, due to Beilinson-Drinfeld
(see \cite{HL}, \cite{B} for details).
Fix a vertex algebra $V$.
The $k[[t]]$-module $V[[t]]$ has a unique structure
of vertex algebra such that
$$\partial_{V[[t]]}=\partial_V+\partial_t,\quad
1_{V[[t]]}=1_V,\quad
(at^n)(z)=(t+z)^na(z),$$
for any elements $a\in V$, $n\in\ZZ$.
Let $\Ac_V$ be the sheaf on $C$ associated to the
$k[[t]]$-module $V[[t]]\cdot dt$.
The sheaf $\Ac_V$ has a unique structure of a chiral algebra over $C$
such that the field $\partial_t$ acts on $\Ac_V$ as the operator
$\partial_{V[[t]]}$, and such that the chiral product is induced by the map
$$V\otimes V[[t_1,t_2]][(t_1-t_2)^{-1}]\to
V[[t_1,t_2]][(t_1-t_2)^{-1}]/V[[t_1,t_2]]$$
which takes the element $f(t_1,t_2)a\boxtimes b$,
with $f(t_1,t_2)\in k(t_1-t_2)$ and $a,b\in V$,
to the element $f(t_1,t_2)a(t_1-t_2)(b)+V[[t_1,t_2]].$
Similarly, if $W$ is a $V$-module
then the $k[[t]]$-module $W[[t]]$ has a natural
structure of a $V[[t]]$-module,
and the corresponding sheaf $\Mc_W$ on $C$
has a natural structure of a $\Ac_V$-module.
Conversely, we have the following.

\proclaim{(6.2.3) Lemma}
Assume that $C$ is a smooth curve.
Fix a point $0\in C$ and a formal coordinate $t$ at 0.

(a) Let $\Ac$ be a chiral algebra on $C$.
Assume that $\Ac$ is a locally free $\Oc_C$-module.
The fiber, $V$, of $\Ac$ at 0
has a unique structure of a vertex algebra
such that the  chiral algebra $\Ac_V$
is isomorphic to $\Ac|_{\Spec k[[t]]}$.

(b) Let $\Mc$ be a module over a chiral algebra $\Ac$ on $C$.
Assume that $\Ac$, $\Mc$ are locally free $\Oc_C$-modules.
Let $V$, $W$ be the fibers of $\Ac$, $\Mc$ at 0.
The space $W$ has a unique structure of a module over $V$,
see Part $(a)$, such that the $\Ac_V$-module $\Mc_W$
is isomorphic to $\Mc|_{\Spec k[[t]]}$.

\endproclaim

\noindent Thus, Proposition 6.1.2 gives the following.

\proclaim{(6.2.4) Theorem}
The De Rham complex $\Cc\Dc\Rc(\o_X)$ is a sheaf of
vertex algebras on $X$.
For any right $\Dc_X$-module $\Mc$ the
De Rham complex $\Cc\Dc\Rc(\Mc)$ is a sheaf
of $\Cc\Dc\Rc(\o_X)$-modules.
\endproclaim

\vskip.3cm

\noindent {\bf (6.3) The vacuum module and
the chiral de Rham complex.} Here we recall
the original construction of the chiral  de Rham complex
of $X$ as given in [MSV]. One first considers
the case $X = \AA^d$. Similarly to Example 5.4.2, let
$CD_M^N$ be the $\ZZ/2\ZZ$-graded $k$-algebra generated
by even elements $a_{in},a^*_{in}$ and
odd elements $b_{in},b^*_{in}$,
with  $i=1,...,d$ and $-N\leq n\leq M$ for $a_{in}, b_{in}$
and $-M\leq n\leq N$ for $a_{in}^*, b_{in}^*$,
modulo the relations
$$\matrix
a^*_{im}a_{jn}-a_{jn}a^*_{im}
=\delta_{ij}\delta_{m,-n},\quad\hfill&
a_{im}a_{jn}-a_{jn}a_{im}=
a^*_{im}a^*_{jn}-a^*_{jn}a^*_{im}=0,\hfill\cr\cr
b_{im}^*b_{jn}+b_{jn}b_{im}^*=
\delta_{ij}\delta_{m,-n},\quad\hfill&
b_{im}b_{jn}+b_{jn}b_{im}=
b^*_{im}b^*_{jn}+b^*_{jn}b^*_{im}=0.\hfill
\endmatrix$$
We further require that the letters $a$ and $b$ commute
in all cases. Let $\tilde{CD} = \ind_{M,N}CD_M^N$.
Consider the super vector space $k^{d|d}$ with basis consisting of
even vectors $v_1,...v_d$ and odd vectors $v_{d+1},...v_{2d}$.
The space
$$\tilde\hen=\bigl(k[t,t^{-1}]\oplus k[t,t^{-1}]dt\bigr)\otimes k^{d|d}
\oplus k\gamma$$
is then a Lie super algebra over $k$ with respect to the brackets given by
$$[f\otimes v_i,\omega\otimes v_j]=\delta_{ij}\Res(f\cdot\omega)\gamma,$$
all other brackets being zero.
It is clear that
$$\tilde{CD}=U(\tilde\hen)/(\gamma-1),$$
with $a_{in}\mapsto v_it^n$,
$a^*_{in}\mapsto v_it^{n-1}dt$,
$b_{in}\mapsto v_{i+d}t^n$,
$b^*_{in}\mapsto v_{i+d}t^{n-1}dt$.
Let also
$$\hen=\bigl(k((t))\oplus k((t))dt\bigr)\otimes k^{d|d}
\oplus k\gamma$$
with the bracket defined in the same way, and
$$CD=U(\hen)/(\gamma-1).$$
Obviously $\tilde{CD}\subset CD$.
Let
$$\tilde\hen^+=\bigl(k[t]\oplus k[t]dt\bigr)\otimes k^{d|d},\quad
\hen^+=\bigl(k[[t]]\oplus k[[t]]dt\bigr)\otimes k^{d|d}.$$
These are Abelian subalgebras in $\tilde\hen,\hen$.
Set $\tilde{CD}^+=U(\tilde\hen^+)\tilde{CD}$,
$CD^+=U(\hen^+)CD$.
The vacuum module $\tilde{Vac}=\tilde{CD}/\tilde{CD}^+$, $Vac=CD/CD^+$
are identified as in (5.4.3)(a).
We denote by $1\in Vac$ (the vacuum vector)
the image of $1\in CD$.
As well-known, $Vac$ has a structure of a vertex algebra
such that the generating series associated to the (-1)- and
0-modes of $a_i, b_i$ are given by:
$$(a_{i,-1}1)(z)=\sum_{n\in\ZZ}a_{in}z^{-n-1},\quad
(a^*_{i0}1)(z)=\sum_{n\in\ZZ}a^*_{in}z^{-n},$$
and similarly with $(b_{i,-1}1)(z)$, $(b^*_{i0}1)(z)$.
The generating series associated to other modes
are obtained by differentiation, using
the action of $\partial$ given by
$$\partial(a_{in})=na_{i,n-1},\quad
\partial(a^*_{in})=(n-1)a^*_{i,n-1},$$
and similarly for $b_{in}^*$, $b_{in}$. The map
$$\delta=\Sum_{i,n}a^*_{in}b_{i,-n}\,:\,Vac\to Vac$$
is a derivation of vertex algebras with zero square.

Setting $x_i = a_{i0}$ makes $Vac$ into a module
over $k[x_1, ..., x_n]$, the coordinate ring of
$\AA^d$. We denote this ring shortly by $k[x]$.
In [MSV] the authors consider the quasicoherent
sheaf $\Omega^{ch}_{\AA^d} := Vac\otimes_{k[x]} \Oc_{\AA^d}$
corresponding to $k[x]$-module $Vac$ and
extend the vertex algebra structure to it.
One also has a vertex algebra structure on
$$Vac^\wedge:=Vac\otimes_{k[x_i]}k[[x_i]].$$

Now let $X$ be a smooth algebraic variety,
$U\subset X$ be an
open subset and $\phi: U\to \AA^d$ be an \'etale map.
Let $0\in X$ be a point such that $\phi(0)=0$
and let $X^\wedge$ be the formal neighborhood of $0$ in $X$.
Then
$x'_i=\phi^*x_i$ are the coordinates on
$X^\wedge$.

In \cite{MSV} the authors construct
a sheaf $\Omega^{ch}_X$ of differential vertex algebras on $X$ as the unique such sheaf satisfying the following
condition.
For any $\phi$ as above, there is an isomorphism of vertex algebras
$$\phi_{ch}\,:\,\Omega^{ch}_{X}\otimes_{\Oc_X}\Oc_{X^\wedge}
\to\phi^*\Omega^{ch}_{\AA^{d}} \otimes_{\Oc_{\AA^d}}
\Oc_{\AA^{d\wedge}} = Vac^\wedge$$
which coincides, for
$U\subset\AA^d$, with the automorphism of
$Vac^\wedge$
introduced in \cite{MSV, Theorem 3.7}.
Our aim in the rest of this paper is to prove the following
fact.

\proclaim{(6.3.1) Theorem}
There is an isomorphism of
sheaves of differential vertex algebras $\Omega^{ch}_X\simeq\Cc\Dc\Rc(\o_X)$.
\endproclaim

\vskip .3cm

\noindent {\bf (6.4) The factorization algebra
associated to the vacuum module.}
As the first step in  proving Theorem 6.3.1, let us
describe
 the factorization algebra corresponding to $Vac =
\Gamma(\AA^d, \Omega^{ch}_{\AA^d})$.
In this, we follow [BD1], [G] :
the constructions below are  a particular instance
of the general concept of the chiral enveloping
algebra of a Lie*-algebra. For the convenience
of the reader we give a self-contained presentation.

Fix a smooth projective curve $C$. Recall
that $\o_C$ is the sheaf of 1-forms on $C$
(in the Zariski topology). Let $I$ be a finite set.
Consider the product $C^I\times C$ and its projections
$p, q$ to $C^I$ and $C$. Let us specialize the notation
of (3.3) to the case when $S=C^I$ and
$f_I = \op{Id}: C^I\to C^I$. We denote
the subvariety $\Gamma(f_I)$ simply by
$$\Gamma_I = \bigl\{ \bigl( (c_i), x\bigr) \in C^I\times C\bigl| x\in \{c_i\}\bigr\}.$$
Similarly we write $\Oc_I^\wedge, \Kc_I^\wedge$
for $\Oc_{f_I}^\wedge$, $\Kc_{f_I}^\wedge$.
Set
$$\Oc_{[[I]]} = p_*\Oc_I^\wedge, \quad \Oc_{((I))} = p_*\Kc_I^\wedge,$$
$$\o_{[[I]]} = p_*(\Oc_I^\wedge \otimes q^*\o_C), \quad
\o_{((I))} = p_*(\Kc_I^\wedge\otimes q^*\o_C).$$
These are (non-quasicoherent) sheaves of $\Oc$-modules
on $C^I$. Informally, the ``fiber" of, say, $\o_{((I))}$ at a
point $(c_i)\in C^I$ is the space of sections of
$\o_C$ on the punctured formal neighborhood of the set
$\{c_i\}$, and similarly in the other cases. Note
that the sum of residues defines a morphism
$$\Res_{(c_i)}\,:\,\o_{((I))}\to\Oc_{C^I},$$
trivial on $\o_{[[I]]}$.

Consider the super-vector space $k^{d|d}$ as in (6.3). 
The sheaf
$$\hen_I = \bigl(\Oc_{((I))}\oplus\o_{((I))}\bigr)\otimes
k^{d|d}\,\,\oplus \Oc_{C^I}\cdot\gamma$$
is then a Lie superalgebra in the  category of left
$\Dc_{C^I}$-modules, with respect to the super-bracket
$$\matrix
[\Oc_C^{d|d},\Oc_C^{d|d}]=[\o_C^{d|d},\o_C^{d|d}]=0,\hfill\cr\cr
[v_i\otimes f, v_j\otimes\o]=\delta_{ij}\Res_{(c_i)}(f\o)\cdot\gamma,\quad
\forall f\in\Oc_C,\forall\o\in\o_C,\hfill\endmatrix$$
and with $\gamma$ being a central element.
Similarly, let
$\hen^+_I\subset\hen_I$ be the super-Lie subalgebra
$(\Oc_{[[I]]}\oplus\o_{[[I]]})\otimes k^{d|d}$.
For any surjective map $J\twoheadrightarrow I$
there are obvious isomorphisms
$$\Delta^{(J/I)*}\hen_J\simto\hen_I,\qquad
j^{(J/I)*}\bigl(\Prod_I\hen_{J_i}\bigr)\simto j^{(J/I)*}\hen_J.
$$
Let $\Uc_{C^I}$ be the quotient of the associative enveloping algebra
of $\hen_I$, in the category of left $\Dc_{C^I}$-modules,
by the right ideal generated by $\gamma-1$.
Consider the sheaf
$$\Vac_{C^I}=\Uc_{C^I}\big/\Uc^+_{C^I},$$
where $\Uc^+_{C^I}\subset\Uc_{C^I}$
is the right ideal generated by $\hen^+_I.$
The collection $(\Vac_{C^I})$ is clearly
a factorization algebra.
To simplify we may omit the subscript $C$,
writing $\Vac$ instead of $\Vac_C$.

\proclaim{(6.4.1) Lemma}
$Vac$ is isomorphic, as a vertex algebra,
to the fiber of the chiral
algebra $\Vac$ at any point of $C$.
\endproclaim

\noindent {\sl Proof:}
Fix a point $0\in C$ and a formal coordinate $t$ at 0.
Let $D=\Spec k[[t]]$ be the formal neighborhood of 0 in $C$.
We compute the chiral product, $\mu$, on the right $\Dc_C$-module
$\Vac^r:=\Vac\otimes_{\Oc_C}\o_C$.
The scheme $C\times C$ is equipped with the coordinates
$t:=t\boxtimes 1$, $z:=1\boxtimes t$.
The $\Oc_C$-module $\hen_{\{1\}}$ is locally free, such that
$$\Gamma(D,\hen_{\{1\}})=
(k^{d|d}\oplus k^{d|d}dz)[[z-t,t]][(z-t)^{-1}]\oplus
k[[t]]\cdot\gamma.$$
Thus the map
$$(z-t)^{m} v_i\mapsto a_{im},\quad
(z-t)^{m} v_{i+d}\mapsto b_{im},$$
$$(z-t)^{m-1} v_idz\mapsto a^*_{im},\quad
(z-t)^{m-1} v_{i+d}dz\mapsto b^*_{im}$$
extends uniquely to an isomorphism of $k[[t]]$-vector spaces
$\Gamma(D,\Uc_C)\simto CD[[t]]$, where $CD$ was defined in (6.3).
Let $1,1_C,1_{C^2}$ be the vacuum elements of
$Vac$, $\Vac$, $\Vac_{C^2}$.
We consider the unique isomorphism of $CD[[t]]$-modules
$\Gamma(D,\Vac)\simto Vac[[t]]$ such that $1_C\mapsto 1$.

The scheme $C^2\times C$ is equipped with the local coordinates
$t_i:=t_i\boxtimes 1, z:=1\boxtimes t$, $i=1,2$.
Put $R=k[[t_1,t_2]]$.
The $\Oc_{C^2}$-module $\hen_{\{1,2\}}$ is locally free, such that
$$\Gamma(D^2,\hen_{\{1,2\}})=
(R^{d|d}\oplus R^{d|d}dz)[[z-t_1,z-t_2]]
[(z-t_1)^{-1},(z-t_2)^{-1}]\oplus R\cdot\gamma.$$
Let $T_1$ (resp. $T_2$) be the Taylor expansion
$$R[[z-t_1,z-t_2]][(z-t_1)^{-1},(z-t_2)^{-1},(t_1-t_2)^{-1}]
\to R((z-t_2))((z-t_1))$$
(resp. $R((z-t_1))((z-t_2))$).
The factorization map
$$j_\bullet j^\bullet(\hen_{\{1,2\}})
\simto
j_\bullet j^\bullet(\hen_{\{1\}}\times\hen_{\{1\}})$$
takes an element
$a\in\Gamma\bigl(D^2,j_\bullet j^\bullet(\hen_{\{1,2\}})\bigr)$
to $\bigl(T_1(a),T_2(a)\bigr)$.
It induces an action of the
sheaf of Lie algebras $j_\bullet j^\bullet(\hen_{\{1,2\}})$
on $j_\bullet j^\bullet(\Vac^r\boxtimes\Vac^r)$.
The factorization map
$$j_\bullet j^\bullet(\Vac^r\boxtimes\Vac^r)\simto
j_\bullet j^\bullet(\Vac^r_{C^2})$$
is the unique morphism of sheaves of
$j_\bullet j^\bullet(\hen_{\{1,2\}})$-modules
taking $1_C\boxtimes 1_C$ to $1_{C^2}$.
The chiral product $\mu$ is the composition of the chain of maps
$$j_\bullet j^\bullet\bigl(\Vac^r\boxtimes\Vac^r\bigr)\simto
j_\bullet j^\bullet\bigl(\Vac_{C^2}^r\bigr)\to
j_\bullet j^\bullet\bigl(\Vac_{C^2}^r\bigr)/\Vac_{C^2}^r=
\Delta_\bullet\Delta^\bullet\bigl(\Vac_{C^2}^r\bigr)=
\Delta_\bullet\bigl(\Vac^r\bigr).$$
The right $\Gamma(D^2,\Dc_{C^2})$-module
$\Gamma\bigl(D^2,\Delta_\bullet(\Vac^r)\bigr)$
is spanned by the symbols $a(t)\delta(t_1-t_2)$,
for any $a(t)\in\Gamma(D,\Vac^r)$, modulo the relations
$$\bigl(a(t)\delta(t_1-t_2)\bigr)(\partial_{t_1}+\partial_{t_2})=
\bigl(a(t)\partial_t\bigr)\delta(t_1-t_2),$$
$$(a(t)\delta(t_1-t_2))f(t_1,t_2)=\bigl(a(t)f(t,t)\bigr)\delta(t_1-t_2),$$
for any $f(t_1,t_2)\in k[[t_1,t_2]].$
Fix $b\in\Gamma(D,\Vac^r)$.
Note that
$$T_2\bigl((z-t_1)^{-1}\bigr)=
-\sum_{m\geq 0}(t_1-t_2)^{-m-1}(z-t_2)^{m}.$$
Hence,
$$\matrix
\mu\bigl((t_1-t_2)^{n}(a_{i,-1}1_C)\boxtimes b\bigr)
&=\mu\bigl((a_{i,-1}1_C)\boxtimes b\bigr)(t_1-t_2)^{n}
\hfill\cr\cr
&=\sum_{m\in\ZZ}a_{im}b\,\delta(t_1-t_2)\partial_{t_2}^{(m-n)}.\hfill
\endmatrix$$
Similarly, we get
$$\mu\bigl((t_1-t_2)^{n}(b_{i,-1}1_C)\boxtimes b\bigr)
=\sum_{m\in\ZZ}b_{im}b\,\delta(t_1-t_2)\partial_{t_2}^{(m-n)},$$
$$\mu\bigl((t_1-t_2)^{n}(a^*_{i0}1_C)\boxtimes b\bigr)
=\sum_{m\in\ZZ}a^*_{im}b\,\delta(t_1-t_2)\partial_{t_2}^{(m-n)},$$
$$\mu\bigl((t_1-t_2)^{n}(b^*_{i0}1_C)\boxtimes b\bigr)
=\sum_{m\in\ZZ}b^*_{im}b\,\delta(t_1-t_2)\partial_{t_2}^{(m-n)}.$$
On the other hand, the chiral product
associated to the vertex algebra $V[[t]]$ is the map
$$Vac\otimes Vac[[t_1,t_2]][(t_1-t_2)^{-1}]\to
Vac[[t_1,t_2]][(t_1-t_2)^{-1}]/Vac[[t_1,t_2]]$$
taking
$(t_1-t_2)^{n}(a_{i,-1}1)\boxtimes b$
to
$$\sum_{m\in\ZZ}a_{im}b\,\partial_{t_2}^{(m-n)}\delta(t_1-t_2),$$
where $\partial_{t_2}^{(m)}\delta(t_1-t_2)$ stands for the element
$(t_1-t_2)^{-m-1}+k[[t_1,t_2]]$, and similarly for $b_{i,-1}$,
$a^*_{i0}$, $b^*_{i0}$.
Thus, to prove that the chiral algebra
$\Vac^r$ is isomorphic to the chiral algebra on $D$
built from $Vac$ as in (6.2)
it is sufficient to check
that the corresponding right $\Dc_C$-modules coincide.
See \cite{BD1, Remark 3.4.8.$(i)$} for an elementary definition of the
the canonical left $\Dc_C$-module structure on $\Vac$.
By construction we have $\partial_t(1_C)=0$.
It is easy to see that
$\partial_t\bigl(a_{im}1_C\bigr)=ma_{i,m-1}1_C$
for all $m<0$.
Hence, the operators $\partial$ on $\Gamma(D,\Vac^r)$ and $Vac[[t]]$
coincide on $a_{im}$.
The case of $a^*_{im}$, $b^*_{im}$, $b_{im}$ is similar.
\qed

\vskip .3cm

\noindent {\bf (6.5) The action of \'etale morphisms I.}
To prove Theorem 6.3.1 in full generality,
it suffices to establish the following lemma.
Let $U\subset X$ be any affine open set
and $\phi\,:\,U\to\AA^d$ be any \'etale map.
Fix a point $0\in U$ such that $\phi(0)=0$.
Let $X^\wedge$ be the formal neighborhood of 0.
In particular we write $\AA^{d,\wedge}$ for $(\AA^d)^\wedge$.

\proclaim{(6.5.1) Lemma}
(a) 
There is an isomorphism of differential vertex algebras
$F_\phi\,:\,\phi^*\Omega^{ch}_{\AA^{d,\wedge}}
\simto k[X^\wedge]\otimes_{\Oc_X}\Cc\Dc\Rc(\o_X)$.

(b) If $X=\AA^d$ then $F_\phi^{-1}\circ F_{\Id}
=\phi_{ch}$ is the isomorphism
constructed in [MSV].
\endproclaim

\noindent The plan of the proof is as follows.
We will construct an isomorpism $(F_I)$ of factorization
algebras and obtain $F_\phi$ as the fiber of $F_{\{1\}}$
at a point of $C$.
It is enough to
assume that $C=\AA^1$. Set
$$CD_I=\Gamma(\AA^I,\Uc_{C^I}),\quad
CD^+_I=\Gamma(\AA^I,\Uc^+_{C^I}),$$
$$DR(\o_U)_I=
\Gamma\bigl(\Lc(U)_{\AA^I},\Cc\Dc\Rc(\o_U)_{\AA^I}\bigr),
\quad
Vac_I=\Gamma(\AA^I,\Vac_{C^I}).$$
Thus, $Vac_I=CD_I/CD^+_I.$
If $X=\AA^d$ we introduce the algebra 
$$\tilde{CD}_I=\ind_{n,N}CD_{X_n^N},$$
where
$$A=k[\AA^I],\quad A_n^N=A[a^{(j)}_{l\nu};-N\leq l\leq n],\quad
X_n^N=\Spec A_n^N,\leqno(6.5.2)$$
see Example 5.4.2.
Note that $\tilde{CD_I}$ is a subalgebra of $CD_I$.
Let
$$\tilde{CD}^+_I=\ind_{n,N}CD_{X_n^N}^+,$$
where $CD_{X_n^N}^+$ is the right ideal generated by 
$a^{(j)}_{l\nu}$,
$a^{(j)*}_{l+1,\nu}$,
$b^{(j)}_{l\nu}$,
$b^{(j)*}_{l+1,\nu}$ 
with $l\geq 0$ and $\tilde{Vac}_I=\tilde{CD}_I/\tilde{CD}^+_I.$ 

\proclaim{(6.5.3) Lemma}
$\tilde{Vac}_I$ is an irreducible $\tilde{CD}_I$-module
and the natural map of vector spaces $\tilde{Vac}_I\to Vac_I$
is an isomorphism.
\endproclaim

\noindent{\sl Proof:}
Irreducibility follows from the fact that $CD_{X_n^N}/CD^+_{X_n^N}$ 
is irreducible over $CD_{X_n^N}$.
The isomorphism follows from the normal form of elements of 
$Vac_I$ and $\tilde{Vac}_I$.
\qed

\vskip3mm

Now, to prove Lemma (6.5.1) we will construct a right action
$$DR(\o_U)_I\otimes\tilde{CD}_I\to DR(\o_U)_I\leqno(6.5.4)$$
commuting with the factorization maps.
To prove that the factorization algebras $(\Vac_{\AA^I})$,
$(\Cc\Dc\Rc(\o_U)_{\AA^I})$ are isomorphic,
it is then sufficient to check that the right
$\tilde{CD}_I$-module $DR(\o_U)_I$ has a cyclic vector
whose annihilator is $\tilde{CD}^+_I$.
Observe that, since the map $(6.5.3)$ depends on $\phi$,
the resulting isomorphism of sheaves of vertex algebras
$\phi^*\Omega^{ch}_{\AA^{d,\wedge}}\simto k[X^\wedge]\otimes_{\Oc_X}\Cc\Dc\Rc(\o_X)$ will also
depend on $\phi$.

\vskip .2cm

\noindent {\bf (6.5.5) The case $U=\AA^d$.}
First we consider the particular case where $U=\AA^d$ and $\phi=\Id$. 
What we do in this case is to provide an
explicit identification of the right
$\tilde{CD}_I$-module $DR(\o_{\AA^d})_I$.
To simplify we set $\Lc=\Lc(\AA^d)_{\AA^I}$, etc.
Then, in the notations of (6.5.2),
$$\Lc_n=\ind_N\Spf\bigl(A_n^0[[a_{l\nu}^{(j)};-N\leq l<0]]\bigr),\quad
\Lc^\eps_n=\Spec\bigl(A_n^{N_\eps}\big/
(a_{l\nu_1}^{(j)}\cdots a^{(j)}_{l\nu_{1+\eps_l}};l<0)\bigr),$$
see $(3.9.1)$, where we set $N_\eps=\max\{l;\,\eps_{-l}\neq 0\}$
for each $\eps\in\EE$.
To simplify again we set $X_n^\eps=X_n^{N_\eps}$.
Thus $X^\eps_N$ is an affine space of finite dimension.
There are closed embeddings $\Lc^0_n\subset\Lc_n^\eps\subset X_n^\eps$.
Let $i_{n\eps}: \Lc^0_n\hookrightarrow X^\eps_n$ be the composite embedding.
We write $\o_{\Lc^0_n, X^\eps_n}$ for the right
$\Dc$-module $i_{n\eps\bullet}\o_{\Lc^0_n}$ on $X^\eps_n$.
Let $\Dc\Rc(\o^\eps_n)\in\Ob_{\Lc^\eps_n}$
be the subsheaf of
$\Dc\Rc(\o_{\Lc^0_n,X^\eps_n})$
consisting of the sections supported
(scheme-theoretically) on $\Lc^\eps_n$.
By definition,
$$
\ind_\eps\Gamma\bigl(\Lc^\eps_n,\Dc\Rc(\o_n^\eps)\bigr)=
\ind_N\Gamma\bigl(\Lc^0_n,\Dc\Rc(\o_{\Lc^0_n,X^N_n})\bigr)=
DR_{\Lc^0_n,X^\infty_n}.$$
Let us denote this space by $DR(\o_n)_I$.
Hence,
$$DR(\o_{\AA^d})_I=\ind_nDR(\o_n)_I[nd],\leqno(6.5.6)$$
and there is a right $\ind_N CD_{X_n^N}$-action on $DR(\o_n)_I$,
such that $DR(\o_n)_I$ is the quotient of $\ind_NCD_{X^N_n}$ by
$\ind_NCD^+_{X^N_n}.$
Using (6.5.4) we get a right action of $\tilde{CD}_I$
on $DR(\o_{\AA^d})_I$ such that
$$DR(\o_{\AA^d})_I\simeq\tilde{CD}_I/\tilde{CD}^+_I\simeq \tilde{Vac}_I=Vac_I,$$
thus achieving our goal in the case $X=\AA^d$.

\proclaim{(6.5.7) Corollary}
Theorem 6.3.1 is true for $X=\AA^d$.
\endproclaim

\noindent {\bf (6.6) Etale change of coordinates
in Clifford algebras.} In order to prove Lemma 6.5.1
for general $\phi: U\to \AA^d$ we need
some elementary observations about Clifford
algebras.

If $X$ is a smooth algebraic variety, we denote
by $\Cc\Dc_X$ the sheaf of differential operators
in $\Omega^\bullet_X$, the commutative superalgebra
of differential forms.

In particular, if $X=\AA^d$ with coordinates $x_1, ..., x_d$,
then $x_i, dx_i$ are  free generators of
(the algebra of global sections of)
$\Omega^\bullet_X$ and we denote
$\partial_i = \partial/\partial x_i$ and
$\xi_i = ``\partial/\partial\, dx_i"$
the corresponding derivations, which are thus global
sections of $\Cc\Dc_{\AA^d}$.

Let now $U$ be an affine open subset of a smooth variety $X$, and
let $\phi\,:\,U\to\AA^d$ be an \'etale map.
Let $x'_i=\phi^*x_i$ be the coordinate on $X^\wedge$.
There are then uniquely determined derivations
$\partial'_i,\xi'_i$ of $\Omega_{X^\wedge}$ such that
$$[\partial'_i,x'_j]=[\xi'_i,dx'_j]=\delta_{ij},\quad
[\partial'_i,dx'_j]=[\xi'_i,x'_j]=0.$$
Set $CD_{X^\wedge}=k[X^\wedge]\hat\otimes_{\Oc_X}\Cc\Dc_X$.
The \'etale map $\phi$ gives an isomorphism of formal schemes
$X^\wedge\simto\AA^{d,\wedge}$.
Let $\phi'$ be the inverse isomorphism.

\proclaim{(6.6.1) Lemma}
For any $\phi$ as above, there is a unique $k[X^\wedge]$-algebra isomorphism
$$\phi_\sharp\,:\,CD_{X^\wedge}\simto\phi^*CD_{\AA^{d,\wedge}}$$
such that
$$\phi_\sharp(dx'_i)=\sum_j\partial_j\phi_i(x)dx_j,\quad
\phi_\sharp(\xi'_i)=\sum_j\partial_i'\phi'_j(\phi(x))\xi_j,$$
$$\phi_\sharp(\partial'_i)=\sum_j\partial'_i\phi'_j(\phi(x))\partial_j+
\sum_{j,k,l}\partial'_i\partial'_k\phi'_j(\phi(x))\partial_l\phi_k(x)dx_k\xi_j,
$$
where $\phi=(\phi_1,...,\phi_d)$.
\endproclaim

\noindent{\sl Proof:}
It suffices to observe that, since the map $\phi$ is \'etale,
$CD_{X^\wedge}$ is a free $k[X^\wedge]$-module with basis
$$(\partial'_1)^{r_1}(\partial'_2)^{r_2}\cdots(\partial'_d)^{r_d}\otimes
(\xi'_1)^{m_1}\cdots(\xi'_d)^{m_d}(dx'_1)^{n_1}\cdots(dx'_d)^{n_d},$$
where $r_i,m_i,n_i\in\NN$.
Then use the coordinates change formulas,
see \cite{L, chap. II} for instance.
\qed

\vskip .2cm

\noindent{\bf (6.7) The action of  \'etale morphisms II.}
Let now $\phi: U\to \AA^d$ be a general \'etale morphism, with $U$ affine.
By (3.9.2) we have an isomorphism of schemes
$$\Lc^\eps_n(\phi)_{\AA^I}\simeq
\Lc^\eps_n\times_{\Lc^0_0}\Lc^0_0(U)_{\AA^I},\quad
\roman{where}\ \Lc^\eps_n=\Lc^\eps_n(\AA^d)_{\AA^I}.$$
We will use freely the notations in (6.5), (6.6).
There is an obvious map $X_n^N\to\Lc^0_0$, since $\Lc_0^0=X^0_0=\Spec A_0^0$,
which restricts to the map $\Lc^\eps_n\to\Lc^0_0$ when $N=N_\eps$.
For any $N$, consider the fiber product
$$X_{n\phi}^N=X^N_n\times_{\Lc^0_0}\Lc^0_0(U)_{\AA^I}.$$
Denote by $\phi_n^N: X_{n\phi}^N\to X_n^N$
the projection to the first factor.
Being a base change of an \'etale morphism
$\Lc^0_0(\phi)_\Id\,:\,\Lc^0_0(U)_{\AA^I}\to\Lc^0_0$,
see Proposition 3.9.3,
the map $\phi_n^N$ is \'etale.
Set $0=\{a_{0\nu}^{(j)}=0\}\in\Hom_{\Schb}(\AA^I,X_0^0)$.
Fix $0\in\Hom_{\Schb}(\AA^I,X_{0\phi}^0)$ mapping to $0$ by $\phi^N_n$.
Let $X_0^{0,\wedge},X_{0\phi}^{0,\wedge}$
be the formal neighborhoods of 0, and set
$$X_n^{N,\wedge}=X_0^{0,\wedge}\times_{X^0_0}X_n^N,\quad
X_{\phi}^{N,\wedge}=X_{0\phi}^{0,\wedge}\times_{X^0_0}X_n^N.$$
By Lemma 6.6.1 applied to $\phi_n^N$, we have a ring isomorphism
$$\phi_{n\sharp}^N\,:\,
CD_{X_{n\phi}^{N,\wedge}}\simto
CD_{X_{n}^{N,\wedge}}.$$
Let
$$\tilde{CD}^\wedge_I=\ind_{n,N}CD_{X_n^{N,\wedge}},\ 
\tilde{CD}^\wedge_{I,\phi}=\ind_{n,N}CD_{X_{n\phi}^{N,\wedge}},$$
so that the $\phi^N_{n\sharp}$ give a ring isomorphism
$$\phi^\infty_{\infty\sharp}\,:\,
\tilde{CD}^\wedge_{I,\phi}\to\tilde{CD}^\wedge_I.
\leqno(6.7.1)$$
Let also
$$\tilde{CD}^{\wedge,+}_I=\ind_{n,N}CD^+_{X_n^{N,\wedge}},$$
where 
$$CD_{X_n^{N,\wedge}}^+=k[X_n^{N,\wedge}]\otimes_{k[X_n^N]}CD^+_{X_n^N}$$
and $CD^+_{X_n^N}$ is introduced after (6.5.2).
We have then the vacuum modules 
$$Vac_I^\wedge=\tilde{CD}_I^\wedge/\tilde{CD}_I^{\wedge,+},$$
which form the factorization algebra corresponding to the vertex algebra
$Vac^\wedge$ defined in (6.3).
Note that
$$DR(\o_U)_I=\ind_{n, N} DR (\o_{\Lc^0_n(U)_{\AA^I},X^N_{n\phi}})[nd],$$
and after tensoring with $k[U^\wedge]$ we get a module
$DR(\o_{U^\wedge})_I$ over $\tilde{CD}^\wedge_{I,\phi}$.

\proclaim{(6.7.2) Lemma} With respect to the above
structure of a $\tilde{CD}^\wedge_I$-module, $DR(\o_{U^\wedge})_I$
is isomorphic to $Vac^\wedge_I$.
\endproclaim

\noindent {\sl Proof:} Follows from the fact that
each $DR(\o_{\Lc_n^0(U), X^N_{n\phi}})$
is isomorphic to the vacuum module over
$CD_{X^N_{n\phi}}$.
\qed

\vskip .2cm

Note that both $DR(\o_{U^\wedge})_I$ and $Vac^\wedge_I$ have distinguished
generators. Namely $DR(\o_{U^\wedge})_I$ is the limit of an inductive system
with the first term $DR(\o_{\Lc^0_0(U)_{\AA^I}/\AA^I})$.
But for any smooth morphism $Z\to S$ there is a canonical element 
$1_{Z/S}$ in $DR(\o_{Z/S})$ and we take $1'_I\in DR(\o_{U^\wedge})_I$ to be
the image of $1_{\Lc_0^0(U)_{\AA^I}/\AA^I}$ in the limit.
The generator $1''_I\in Vac^\wedge_I$ is the image of
$1\in\tilde{CD}^\wedge_I$. We denote $F_I: Vac^\wedge_I\to DR(\o_{U^\wedge})_I$
the unique module isomorphism taking $1''_I$ to $1'_I$.
Let $\Fc_I: \Vc ac^\wedge_I\to \Dc\Rc(\o_{U^\wedge})_I$ be the
corresponding morphism of quasicoherent sheaves on $\Lc_0^0(U)_{\AA^I}$.

\proclaim{(6.7.3) Lemma}
The $(\Fc_I)$ commute with factorization maps and thus form
an isomorphism of factorization algebras.

\endproclaim

\noindent {\sl Proof:} Follows from the fact that $(1'_I)$
and $(1''_I)$ are compatible with factorization structures:
in the notation of (6.1.1) we have 
$\varkappa^{J/I}(\boxtimes 1'_{J_i}) = 1'_J$ and
similarly for $1''_I$.
\qed

\vskip .2cm

Set $I=\{1\}$.
Let $F_\phi$ be the fiber of
$F_{\{1\}}$ at the point $0\in\AA^1$.
Lemma 6.7.3 implies that $F_\phi$ is a morphism
of vertex algebras. This establishes part (a)
of Lemma 6.5.1

\vskip .2cm

We now prove Lemma 6.5.1(b). So we assume $U\i\AA^d$ and need
to compare two automorphisms of the vertex algebra
$\Gamma(U, \Vc ac)$, namely $F_\phi^{-1}\circ F_{\op{Id}}$
and $\phi_{ch}$.
Notice that the vertex algebra $Vac$ is
strongly generated by the fields
$(a_{i,-1}1)(z)$, $(a^*_{i0}1)(z)$, $(b_{i,-1}1)(z)$, $(b^*_{i0}1)(z)$,
see \cite{K}, and that $\Gamma(U,\Vac)$ is obtained by localization.
Thus it is enough to compare the two automorpisms on the elements
$$(a_{i,-1}1),\,(a^*_{i0}1),\,(b_{i,-1}1),\,(b^*_{i0}1)\in Vac.
\leqno(6.7.4)$$
Since $U\subset\AA^d$, the sheaf of algebras $\tilde{CD}^\wedge_{I,\phi}$ is
identified with $\tilde{CD}^\wedge_I$, so $\phi^\infty_{\infty,\sharp}$ is an automorphism of the latter.

\proclaim{Lemma}
(a) The morphism $\phi^\infty_{\infty\sharp}$ preserves 
$\tilde{CD}_I^{\wedge,+}$ and thus induces an automorphism
$\phi_{\sharp,I}\,:\,Vac^\wedge_I\to Vac^\wedge_I$.

(b) For $I=\{\bullet\}$, the morphism of vertex algebras
$F_\phi^{-1}\circ F_\Id\,:\,Vac^\wedge\to Vac^\wedge$ 
is equal to $\phi_{\sharp,0}$,
which is the fiber over $0\in\AA^1$ of the morphism $\phi_{\sharp,I}$.
\endproclaim

\noindent{\sl Proof:}
(a) is enough to verify for each 
$\phi^N_{n\sharp}\,:\,CD_{X^{N,\wedge}_{n\phi}}\to
CD_{X^{N,\wedge}_{n\phi}}$, in which case it follows from Lemma 6.6.1.
Claim (b) follows from construction of $F_\phi$.
\qed

\vskip3mm

To prove $(6.5.1)(b)$ it suffices therefore to check that $\phi_{\sharp0}=\phi_{ch}$ on elements (6.7.4).
Recall that $\Lc^0_0(X)_{\AA^1}=X\times\AA^1$ for any $X$.
Hence, there is a commutative diagram
$$\matrix
X_{n\phi}^N&\twoheadrightarrow&X^0_{0\phi}&=&U\times\AA^1\cr
{\ss\phi_n^N}\downarrow&&\downarrow&&\downarrow{\ss\phi}\cr
X_n^N&\twoheadrightarrow&X_0^0&=&\AA^d\times\AA^1.
\endmatrix$$
This diagram induces a diagram of $A$-algebra homomorphisms
$$\matrix
CD_{Y_\phi^\wedge}&\hookleftarrow&CD_{X^{0\wedge}_{0\phi}}\cr
{\ss\phi_\sharp}\downarrow&&\downarrow{\ss\phi^0_{0\sharp}}\cr
CD_{Y^\wedge}&\hookleftarrow&CD_{X^{0\wedge}_0}.
\endmatrix$$
Note that the images of the elements
(6.7.4) by $\phi_\sharp$ and $\phi^0_{0\sharp}$ coincide,
modulo the identification
$$a_{i0}\mapsto x_i,\quad b_{i0}\mapsto dx_i,\quad
b^*_{i,-1}\mapsto\xi_i,\quad a^*_{i,-1}\mapsto\partial_i.$$
On the other hand  the images of the elements (6.7.4)
by $\phi_{ch}$ and $\phi^0_{0\sharp}$
coincide, see the formulas \cite{MSV, (3.17)}
for $\phi_{ch}$, and Lemma 6.6.1 for $\phi_{0\sharp}^0$.
We are done.
\qed

\vfill\eject

\vskip3cm
{\eightpoint{
$$\matrix\format\l&\l&\l&\l\\
\phantom{.} & {\text{Mikhail Kapranov}}\phantom{xxxxxxxxxxxxx} &
{\text{Eric Vasserot}}\\
\phantom{.}&{\text{Department of Mathematics}}\phantom{xxxxxxxxxxxxx} &
{\text{D\'epartement de Math\'ematiques}}\\
\phantom{.}&{\text{University of Toronto}}\phantom{xxxxxxxxxxxxx} &
{\text{Universit\'e de Cergy-Pontoise}}\\
\phantom{.}&{\text{100 St. George St. }}\phantom{xxxxxxxxxxxxx} &
{\text{2 Av. A. Chauvin}}\\
\phantom{.}&{\text{Toronto, Ontario M5S 3G3}}\phantom{xxxxxxxxxxxxx} &
{\text{95302 Cergy-Pontoise Cedex}}\\
\phantom{.}&{\text{Canada}}\phantom{xxxxxxxxxxxxx} &
{\roman{France}}\\
&{\text{email: kapranov\@math.toronto.edu}}\phantom{xxxxxxxxxxxxx} &
{\text{email: eric.vasserot\@math.u-cergy.fr}}
\endmatrix$$
}}
\enddocument